\documentclass{clsrbt}
\usepackage{blkarray}
\usepackage{array}

\usepackage{amsmath}
\usepackage{amssymb}
\usepackage[all]{xy}
\newcommand{\xyR}[1]{\xydef@\xymatrixrowsep@{#1}}
\newcommand{\xyC}[1]{\xydef@\xymatrixcolsep@{#1}}
\newtheorem{theorem}{Theorem}[section]
\newtheorem{proposition}[theorem]{Proposition}
\newtheorem{lemma}[theorem]{Lemma}
\newtheorem{corollary}[theorem]{Corollary}
\newnumbered{definition}[theorem]{Definition} 
\newnumbered{example}[theorem]{Example}
\newnumbered{remark}[theorem]{Remark}
\newnumbered{assumption}[theorem]{Assumption}
\newnumbered{notation}[theorem]{Notation}

\usepackage[all]{xy}
\usepackage{url}

\def\rad{\mathop{\mathrm{rad}}}

\def\im{\mathop{\mathrm{im}}}
\def\ker{\mathop{\mathrm{ker}}}

\title[Functorial filtrations for homotopy categories]{Functorial filtrations for homotopy categories\\ of some generalisations of gentle algebras} 

\author{R. Bennett-Tennenhaus}

\dedication{}

\classno{16D70 (primary), 16G20, 16E35 (secondary)}

\extraline{The author was supported by the School of Mathematics at the University of Leeds. The author was also supported by the Alexander von Humboldt Foundation 
in the framework of an Alexander von Humboldt Professorship 
endowed by the German Federal Ministry of Education and Research.}

\begin{document}
\maketitle

\begin{abstract}
We consider algebras defined  over a complete, local and noetherian ground ring. They are gentle algebras in case the ground ring is a field. The unbounded homotopy category of complexes of projective modules is considered. Complexes with finitely-generated homogeneous components are shown to be isomorphic to direct sums of indecomposable \textit{string} and \textit{band}
complexes. The corresponding isoclasses are described, and the Krull-Remak-Schmidt-Azumaya property is verified. This classification problem is solved using the idea of \textit{functorial filtrations}.
\end{abstract}
\section{Introduction.}\label{intro}
 Gentle algebras, as introduced by Assem and Skowro{\'n}ski \cite{AssSko1987}, have the form $\Gamma = kQ/\mathcal{J}$ where $k$ is a field, $Q$ is a quiver, $kQ$ is the path algebra and the ideal $\mathcal{J}\triangleleft kQ$ is generated by length $2$ paths. The combinatorial restrictions \cite[\S1.1, Proposition, R1, R2]{AssSko1987} on $\mathcal{J}$ and $Q$ ensure $\Gamma$ is a finite-dimensional \textit{special biserial algebra}, in the sense of Skowro{\'n}ski and Waschb{\"u}sch \cite[\S1, (SP)]{SkoWas1983}. By a result of Wald and Waschb{\"u}sch \cite[Proposition 2.3]{WalWas1985}, if $\Omega$ is a special biserial algebra, then indecomposables in the category $\Omega\text{-}\mathbf{mod}$ of finite-dimensional modules are \textit{string modules}, \textit{band modules} or non-uniserial projective-injectives. Since $\mathcal{J}$ is generated by paths, indecomposable projective-injective $\Gamma$-modules are uniserial.

Various authors have also capitalised on the aforementioned restrictions of $Q$ and $\mathcal{J}$ to study derived categories of gentle algebras.  Schr{\"o}er and Zimmermann \cite[Corollary 1.2]{SchZim2003} have shown gentle algebras are closed under derived equivalence. Bekkert and Merklen \cite[Theorem 3]{BekMer2003} have shown indecomposable objects in the bounded derived category $\mathcal{D}^{b}(\Gamma\text{-}\mathbf{mod})$ are (what we call) \textit{string complexes} or \textit{band complexes}. The repetitive algebra $\Omega$ of $\Gamma$, as introduced by Hughes and Waschb{\"u}sch \cite{HugWas1983}, was shown by Pogorza{\l}y and Skowro{\'n}ski \cite[Lemma 8]{PogSko1991} to be special biserial. Bobi\'{n}ski \cite[\S 5, \S 6]{Bob2011} has studied connections
between string and band modules over $\Omega$ and string and band complexes over $\Gamma$, describing almost split triangles in the category of perfect complexes. Gei{\ss} and Reiten \cite[Theorem 3.4]{GeiRei2005} have shown gentle algebras are Iwanaga-Gorenstein rings. Kalck \cite[Theorem 2.5(b)]{Kal2015} has since described the singularity category. Arnesen, Laking and Pauksztello \cite[Theorem A]{ArnLakPau2016} have described all morphisms between string and band complexes, and these authors together with Prest \cite[Theorem B]{ArnLakPauPre2017} described all indecomposable pure-injective objects in the unbounded homotopy category $\mathcal{K}(\Gamma\text{-}\mathbf{Proj})$. 
 
In this article we consider more generally \textit{complete gentle algebras}, defined using a quiver $Q$ with relations, and a complete, local and noetherian ground ring $(R,\mathfrak{m},k)$. See \S \ref{2} for details.

For some prime number $p>0$, let $R$ be the ring $\widehat{\mathbb{Z}}_{p}$  of $p$-adic integers. A notable mixed-characteristic example of a complete gentle algebra is the $R$-algebra
\begin{equation}\label{eqn1.1}
\Upsilon_{p}=\left\{\begin{pmatrix}\gamma_{11} & \gamma_{12}\\
\gamma_{21} & \gamma_{22}
\end{pmatrix}\in\mathbb{M}_{2}(\widehat{\mathbb{Z}}_{p})\,\mid\,\gamma_{11}-\gamma_{22},\,\gamma_{12}\in p\widehat{\mathbb{Z}}_{p}\,\right\}
\end{equation} 
If $R$ is a field, then by \cite[Corollary 1.2.11]{Ben2018} any complete gentle $R$-algebra is a finite-dimensional gentle algebra in the sense of \cite{AssSko1987}. Cases where $R$ strictly contains a field give interesting infinite-dimensional examples of complete gentle algebras. To see an example, consider the \textit{completed string algebras} in the sense of Ricke \cite{Ric2017}, which were defined as quotients $\Theta=\overline{kQ}/\overline{(\rho)}$ of the completed path algebra of a quiver $Q$. If the pair $(Q,\rho)$ satisfies \textit{gentle conditions} (see Definition \ref{definition.1.16}) then $\Theta$ is a complete gentle $k[[t]]$-algebra (see Remark \ref{rem}). For example, if $Q$ consists of two loops $a$ and $b$ and $\rho=\{ab,ba\}$ then $\Theta$ is the local $k[[t]]$-algebra $k[[x,y]]/(xy)$, where $t$ acts as $x+y+(xy)$ (see \cite[Example 1.2.30]{Ben2018} for details). The complex of free $k[[x,y]]/(xy)$-modules
\begin{equation}\label{eqn1.2}
P=\xymatrix@C=1.45em{\cdots\ar[r]^{\small{\begin{pmatrix}x^{2}\end{pmatrix}}} & \Lambda\ar[r]^{\small{\begin{pmatrix}y\end{pmatrix}}} & \Lambda\ar[r]^{\small{\begin{pmatrix}x^{2}\end{pmatrix}}} & \Lambda\ar[r]^{\small{\begin{pmatrix}0\\
y
\end{pmatrix}}}  & \Lambda^{2}\ar[rr]^{\small{\begin{pmatrix}x^{2} & 0\\
y^{3} & 0\\
0 & x^{2}
\end{pmatrix}}} &  & \Lambda^{3}\ar[rr]^{\small{\begin{pmatrix}y & 0 & 0\\
0 & x & y
\end{pmatrix}}} &  & \Lambda^{2}\ar[rr]^{\small{\begin{pmatrix}x^{2} & 0\end{pmatrix}}} &  & \Lambda\ar[r] & 0\ar[r] & \cdots}
\end{equation}
is an example of a \textit{string complex}. The name refers to the depiction of $P$ by the schema 
\[
\xymatrix@C=1.65em@R=1.65em{ \ddots\,\,\,\,\,\ar[dr]\sp(0.6){{x}^{2}} &  &  &  &  &  &  &  &  &  & \vdots\ar[d]\\
 &  \Lambda\ar[dr]\sp(0.6){y}\ar@{..}[l]  &  &  &  &  &  &  &  &  & P^{-4}\ar@{..}[lllllllll]\ar[d]^{d_{P}^{-4}}\\
 &  & \Lambda\ar[dr]\sp(0.6){{x}^{2}}\ar@{..}[ll] & &  &  &  \Lambda\ar[dl]\sp(0.45){{y}^{3}}\ar[dr]\sp(0.6){{x}^{2}}\ar@{..}[llll] &  &  &  &  P^{-3}\ar@{..}[llll]\ar[d]^{d_{P}^{-3}}\\
 &  &  & \Lambda\ar[dr]\sp(0.6){y}\ar@{..}[lll] &  &  \Lambda\ar[dl]\sp(0.45){x}\ar@{..}[ll]  &  & \Lambda\ar[dr]\sp(0.6){y}\ar@{..}[ll] &  &  &  P^{-2}\ar@{..}[lll]\ar[d]^{d_{P}^{-2}}\\
 &  &  &  &  \Lambda\ar@{..}[llll]  &  &  &  &  \Lambda\ar[dr]\sp(0.6){{x}^{2}}\ar@{..}[llll]  &  & P^{-1}\ar@{..}[ll]\ar[d]^{d_{P}^{-1}}\\
&  &  &  &  &  &  &  &  &  \Lambda\ar@{..}[lllllllll] & P^{0}\ar@{..}[l]\ar[d]^{d_{P}^{0}}\\
&  &  &  &  &  &  &  &  & & \vdots
}\]
\textit{Band complexes} are defined similarly, but the corresponding schema is finite and joins up with itself. String complexes are indexed by a \textit{homotopy word} $C$ which is \textit{aperiodic}, and band complexes by a \textit{periodic homotopy word} together with an indeceomposable $R[T,T^{-1}]$-module $V$ which is free over $R$. See \S\ref{3} for details. Given a complete gentle algebra $\Lambda$ let $\mathcal{K}(\Lambda\text{-}\boldsymbol{\mathrm{Proj}})$ (respectively $\mathcal{K}(\Lambda\text{-}\boldsymbol{\mathrm{proj}})$) be the unbounded homotopy
category of complexes of (finitely
generated) projective modules. Theorem \ref{theorem.1.1} describes objects in $\mathcal{K}(\Lambda\text{-}\boldsymbol{\mathrm{proj}})$. Theorems \ref{theorem.1.2} and \ref{theorem.1.3} explain the way in which this description is unique.
\begin{theorem}
\label{theorem.1.1}\emph{\cite[Theorem 2.0.1]{Ben2018}} Let $\Lambda$ be a complete gentle algebra.  
\begin{enumerate}
\item Every object in $\mathcal{K}(\Lambda\text{-}\boldsymbol{\mathrm{proj}})$ is isomorphic to a (possibly infinite) direct sum of shifts of string complexes $P(C)$ and shifts of band complexes $P(C,V)$.
\item Each shift of a string or band complex is an indecomposable object in $\mathcal{K}(\Lambda\text{-}\boldsymbol{\mathrm{Proj}})$.
\end{enumerate}
\end{theorem}
A well-known result of Kaplansky \cite[Theorem 2]{Kap1958} says that any projective $R$-module is free. Consider  the category $R[T,T^{-1}]\text{-}\boldsymbol{\mathrm{Mod}}$ of $R[T,T^{-1}]$-modules, and the full subcategory $R[T,T^{-1}]\text{-}\boldsymbol{\mathrm{Mod}}_{R\text{-}\boldsymbol{\mathrm{Proj}}}$ of $R[T,T^{-1}]$-modules which are free as $R$-modules. Note that there is a natural isomorphism $k[T,T^{-1}]\otimes_{R[T,T^{-1}]}-\simeq k\otimes_{R}-$ of functors $R[T,T^{-1}]\text{-}\boldsymbol{\mathrm{Mod}}_{R\text{-}\boldsymbol{\mathrm{Proj}}}\to k[T,T^{-1}]\text{-}\boldsymbol{\mathrm{Mod}}$. Let $\iota$ define the $R$-algebra involution of $R[T,T^{-1}]$
which exchanges $T$ and $T^{-1}$. Define a functor $\mathrm{res}_{\iota,R}:R[T,T^{-1}]\text{-}\boldsymbol{\mathrm{Mod}}_{R\text{-}\boldsymbol{\mathrm{Proj}}}\to R[T,T^{-1}]\text{-}\boldsymbol{\mathrm{Mod}}_{R\text{-}\boldsymbol{\mathrm{Proj}}}$
by setting $\mathrm{res}_{\iota,R}(V)$ to have underlying $R$-module
structure $V$ but where $T$ acts on $v\in\mathrm{res}_{\iota,R}(V)$ by $v\mapsto T^{-1}v$.

The \textit{homotopy letters} which make up a homotopy word $C$ are indexed using a  subset $I_{C}$ of $\mathbb{Z}$. The \textit{inverse} $C^{-1}$ defines a new homtopy word by inverting the letters and reversing their order. Similarly the \textit{shift} $C[t]$ ($t\in\mathbb{Z}$) is defined by reindexing the letters by $i\mapsto i-t$. Any $i\in I_{C}$ defines a vertex $v_{C}(i)$ of the quiver $Q$ defining $\Lambda$, and in this way the set $I_{C}$ defines the projective indecomposable summands of the underlying module of the string or band complex defined by $C$.  Furthermore there is a function $\mu_{C}:I_{C}\to\mathbb{Z}$ which detects the homogeneous degree of each aformentioned summand. See Definitions \ref{definition2} and \ref{def.3.3} for details. Theorem \ref{theorem.1.2} characterises when two shifts of string or band complexes are isomorphic. 
\begin{theorem}
\label{theorem.1.2}\emph{\cite[Theorem 2.0.4]{Ben2018}} Let $\Lambda$ be a complete gentle algebra. Let $C$ and $E$ be homotopy words, let $V$ and $W$ be objects of  $R[T,T^{-1}]\text{-}\boldsymbol{\mathrm{Mod}}_{R\text{-}\boldsymbol{\mathrm{Proj}}}$ and let $n\in\mathbb{Z}$. 
\begin{enumerate}
\item If $C$ and $E$ are aperiodic, then $P(C)[n]\simeq P(E)$ in $\mathcal{K}(\Lambda\text{-}\boldsymbol{\mathrm{Proj}})$ if and only if:
\begin{enumerate}
\item we have  $I_{C}=\{0,\dots,m\}$ and $(I_{E},E,n)=(I_{C},C,0)\text{ or } (I_{C},C^{-1},\mu_{C}(m))$; or
\item we have $I_{C}=\pm\mathbb{N}$ and $(I_{E},E,n)=(\pm\mathbb{N},C,0)\text{ or }(\mp\mathbb{N},C^{-1},0)$; or
\item we have $I_{C}=\mathbb{Z}$ and  $(I_{E},E,n)=(\mathbb{Z}, C^{\pm1}[t],\mu_{C}(\pm t))$ for some $t\in\mathbb{Z}$.
\end{enumerate}
\item If $C$ and $E$ are periodic, then $P(C,V)[n]\simeq P(E,W)$ in $\mathcal{K}(\Lambda\text{-}\boldsymbol{\mathrm{Proj}})$ if and only if:
\begin{enumerate}
\item we have $E=C[t]$, $k\otimes_{R}V\simeq k\otimes_{R}W$ and $n=\mu_{C}(t)$ for some $t\in\mathbb{Z}$; or
\item we have $E=C^{-1}[t]$, $k\otimes_{R}V\simeq k\otimes_{R}\mathrm{res}_{\iota,R} \,W$ and $n=\mu_{C}(-t)$ for some $t\in\mathbb{Z}$.
\end{enumerate}
\item If $C$ is aperiodic and $E$ is periodic, then $P(C)[n]\not\simeq P(E,V)$ in $\mathcal{K}(\Lambda\text{-}\boldsymbol{\mathrm{Proj}})$.
\end{enumerate}
\end{theorem}
\begin{theorem}
\label{theorem.1.3}\emph{\cite[Theorem 2.0.5]{Ben2018}} Let $\Lambda$ be a complete gentle algebra. If two direct sums of shifts of string
and band complexes are isomorphic in $\mathcal{K}(\Lambda\text{-}\boldsymbol{\mathrm{proj}})$, then there is an isoclass
preserving bijection between the summands.
\end{theorem}
For a finite-dimensional gentle algebra $\Gamma$, the classification of objects in $\mathcal{D}^{b}(\Gamma\text{-}\mathbf{mod})$ is due to Bekkert and Merklen \cite{BekMer2003}. These authors constructed a functor to a category of square-zero block matrices which preserves and respects indecomposables. Objects in this category were classified by Bondarenko using the \textit{matrix reductions} method  \cite[Theorem 1]{Bon1975}. 
 Bekkert, Drozd and Furtorny \cite[Theorem 2.7]{BekDroFut2009} note that the technique of reducing to a matrix problem used in \cite{BekMer2003} may be employed to the class of complete gentle $k[[t]]$-algebras $\Theta=\overline{kQ}/\overline{(\rho)}$ discussed before Equation \ref{eqn1.2} (see \cite[Definition 2.5]{BekDroFut2009} and Remark \ref{rem}). Theorems \ref{theorem.1.1}, \ref{theorem.1.2} and \ref{theorem.1.3}, together with Remark \ref{rem}, give a new approach to the aforementioned classification problems solved by Bekkert and Merklen and Bekkert, Drozd and Furtorny (see \cite[Corollaries 2.7.2 and 2.7.6]{Ben2018}).

In our proof of Theorems \ref{theorem.1.1}, \ref{theorem.1.2} and \ref{theorem.1.3} we work directly with the category $\mathcal{K}(\Lambda\text{-}\boldsymbol{\mathrm{Proj}})$ and avoid reducing to a matrix problem. The classification method we adapt is called the \textit{functorial filtrations} method. This approach yields results which appear to be new: our results hold for mixed-characteristic complete gentle algebras, and so Theorems \ref{theorem.1.1}, \ref{theorem.1.2} and \ref{theorem.1.3} strictly generalise \cite[Theorem 2.7]{BekDroFut2009}; string and band complexes whose homogeneous components may be infinitely-generated are indecomposable by Theorem \ref{theorem.1.1}(ii), which generalises the indecomposability statement in \cite[Corollary 3.5]{ArnLakPauPre2017}; in Theorem \ref{theorem.1.1}(i) we classify complexes in $\mathcal{K}(\Lambda\text{-}\boldsymbol{\mathrm{proj}})$ which need not satisfy any bounded-ness conditions; and in Theorem \ref{theorem.7.1} the summands of any direct sum of shifts of string and band complexes have been identified.

Functorial filtrations have been written in Mac Lane's
language of linear relations \cite{Mac1961}, and were used in the past
to classify modules (with certain finiteness conditions) up to isomorphism.
Gel'fand and Ponomarev \cite{GelPon1968} seem to be the first to use this method, during a classification Harish-Chandra modules for the Lorentz group. Their work was interpreted in
the language of functors by Gabriel \cite{Gab1987}. Ringel \cite{Rin1975} then used this approach to describe indecomposable representations of the dihederal 2-group. Since then several classes of rings have had their modules classified by authors adapting the method from \cite{GelPon1968} and \cite{Rin1975}: Brauer graph algebras by Donovan and Freislich \cite{DonFre1978}; locally bounded string algebras by Butler and Ringel \cite{ButRin1987}; clannish algebras \cite{Cra1989}, semidihedral algebras \cite{Cra19892} and string algebras which may be both infinite-dimensional and unital \cite{Cra2018} by Crawley-Boevey; and completed string algebras by Ricke \cite{Ric2017}. Functorial filtrations have also been used to classify representations of significance outside representation theory. Prest and Puninski \cite{PrePun2016} classified pure-injective indecomposable modules over domestic string algebras using linear relations that correspond to subgroups of finite definition. Such objects arise from the model theory of modules. In joint work \cite{BenCra2018} with Crawley-Boevey we employed similar ideas, adapting the focus of the functorial filtrations method in \cite{Cra2018} to classify $\Sigma$-pure-injective modules. 

The results in this article form part of the author's PhD thesis \cite{Ben2018}. The article is organised as follows. In \S\ref{2} we define complete gentle algebras. In \S\ref{3} we define string and band complexes. In \S\ref{structuretheorem} we explain how our classification works by providing a categorical blueprint. In \S\ref{sec:Linear-Relations-and} we consider linear relations over the ring $R$, and introduce the notion of a \textit{reduction} which \textit{meets in} the maximal ideal $\mathfrak{m}$. In \S\ref{sec:Some-Linear-Relations} we explain how each homotopy word gives rise to a functor on the category of complexes, each of which is defined by such linear relations. In the remaining sections of the article we: check compatibility conditions between (string and band complexes) and (linear relations given by homotopy words); we verify that our setting fits the blueprint discussed in \S\ref{structuretheorem}; and we provide proofs of Theorems \ref{theorem.1.1}, \ref{theorem.1.2} and \ref{theorem.1.3}.
\section{Some generalisations of gentle algebras.}\label{2}
 
\begin{assumption}For the remainder of the article we fix:
\begin{enumerate}
\item a commutative,
noetherian and local ring $R$ with maximal ideal $\mathfrak{m}$ and $k=R/\mathfrak{m}$;
\item a finite quiver $Q$ with path algebra $RQ$ and an ideal $\mathcal{J}\triangleleft RQ$ generated by a set paths in $Q$ of length at least $2$; and
\item a surjective $R$-algebra homomorphism $\vartheta:RQ\rightarrow\Lambda$ where $\mathcal{J}\subseteq\ker(\vartheta)$ and $\vartheta(p)\neq0$ for any path $p\notin\mathcal{J}$. Notation is abused by writing $p$ for $\vartheta(p)$.
\end{enumerate}
\end{assumption}
\begin{notation} Let $\mathbf{P}$ be the set of non-trivial paths $p\notin \mathcal{J}$ with head $h(p)$ and tail $t(p)$. For each $t>0$ and each vertex $v$ let $\mathbf{P}(t,v\rightarrow)$
(respectively $\mathbf{P}(t,\rightarrow v)$) be the set of paths $p\in\mathbf{P}$ of length $t$ with $t(p)=v$ (respectively $h(p)=v$). Let $\mathbf{A}$ be the set of arrows in $Q$, $\mathbf{A}(v\rightarrow)=\mathbf{P}(1,v\rightarrow)$ and  $\mathbf{A}(\rightarrow v)=\mathbf{P}(1,\rightarrow v)$. The composition of $a\in\mathbf{A}(\rightarrow v)$ and $b\in\mathbf{A}(u\rightarrow)$ is $ba$ if $u=v$, and $0$ otherwise. If $X$ is a finite set, let $\# X$ be the cardinality of $X$.
\end{notation}
\begin{definition}
\label{definition.1.16}\cite[\S 1, (SP)]{SkoWas1983} We say that the pair $(Q,\mathcal{J})$ satisfies \textit{special conditions} if:
\begin{enumerate}
\item[(SPI)] if $v$ is a vertex then $\#\mathbf{A}(v\rightarrow)\leq2$
and $\#\mathbf{A}(\rightarrow v)\leq2$; and
\item[(SPII)] if $y\in\mathbf{A}$ then $\# \{ x\in \mathbf{A}(h(y)\rightarrow)\mid xy\in\mathbf{P}\}\leq 1$ and $\# \{ z\in \mathbf{A}(t(y)\rightarrow)\mid yz\in\mathbf{P}\}\leq 1$.
\end{enumerate}
 \cite[p.272, R3, R4]{AssSko1987} We say $(Q,\mathcal{J})$ satisfies \textit{gentle conditions} if it satisfies special conditions and:
\begin{enumerate}
\item[(GI)] any path $p\notin\mathbf{P}$ has a subpath $q\notin\mathbf{P}$ of length 2; and
\item[(GII)] if $y\in\mathbf{A}$ then $\# \{ x\in \mathbf{A}(h(y)\rightarrow)\mid xy\notin\mathbf{P}\}\leq 1$ and $\# \{ z\in \mathbf{A}(t(y)\rightarrow)\mid yz\notin\mathbf{P}\}\leq 1$.
\end{enumerate} 
\end{definition}
\begin{example}
\label{exa:p-adic}\cite[Example 1.1.6]{Ben2018} Let $R=\widehat{\mathbb{Z}}_{p}$ and $\mathfrak{m}=p\widehat{\mathbb{Z}}_{p}$. Let $Q$ be the quiver given by two loops $\alpha$
and $\beta$ at a single vertex $v$, and let $\mathcal{J}=\langle \alpha^{2},\beta^{2}\rangle$.
Let $\Lambda$ be the
subring of $2\times2$ matrices whose $ij$-entry $\gamma_{ij}\in R$ satisfies $p\mid \gamma_{11}-\gamma_{22},\gamma_{12}$, and so $\Lambda=\Upsilon_{p}$ from Equation \ref{eqn1.1}. Define $\vartheta:RQ\rightarrow\Lambda$ by (multiplicatively and $R$-linearly) extending the assignments
\[
\vartheta(\alpha)=\begin{pmatrix}0 & 0\\
1 & 0
\end{pmatrix},\quad\vartheta(\beta)=\begin{pmatrix}0 & p\\
0 & 0
\end{pmatrix}.
\]
That $(Q,\mathcal{J})$ satisfies gentle conditions is clear. Furthermore $\rad(\Lambda)$ is the subset of matrices $(\gamma_{ij})$ with $p\mid\gamma_{11},\gamma_{22}$, which is generated as a two-sided ideal by $\alpha$ and $\beta$. Note $\alpha\beta+\beta\alpha=p$ in $\Lambda$, and so $(\rad(\Lambda))^{3}\subseteq\Lambda p\subseteq\rad(\Lambda)$. In particular the ideal  $\rad(\Lambda)/\Lambda p$
of $\Lambda/\Lambda p$ is nilpotent.
\end{example}
\begin{definition}
\label{def.1.13} \cite[Definition 1.1.19]{Ben2018} We say the ring $\Lambda$
is \textit{rad-nilpotent modulo $\mathfrak{m}$} if $(\rad(\Lambda))^{n}\subseteq\Lambda \mathfrak{m}\subseteq \rad(\Lambda)$ for some $n\geq1$ (and so the ideal  $\rad(\Lambda)/\Lambda\mathfrak{m}$
of $\Lambda/\Lambda\mathfrak{m}$ is nilpotent).

\label{definition.3.4} \cite[Definition 1.1.21]{Ben2018} We say that $\Lambda$ is a \textit{quasi-bounded gentle} $R$-algebra if:
\begin{enumerate}
\item the pair $(Q,\mathcal{J})$ satisfies gentle conditions (SPI, SPII, GI and GII);
\item the ring $\Lambda$ is rad-nilpotent
modulo $\mathfrak{m}$;
\item if $a\in\mathbf{A}$ then the $R$-modules $\Lambda a$ and $a\Lambda$ are finitely generated;
\item the ideal $\rad(\Lambda)$ of $\Lambda$ is generated by $\mathbf{A}$; and
\item if $a,a'\in\mathbf{A}$ and $a\neq a'$ then $\Lambda a'\cap\Lambda a=0=a\Lambda\cap a'\Lambda$.
\end{enumerate}
\end{definition}
It is straightforward to check the mixed-characteristic ring $\Upsilon_{p}$ from Equation \ref{eqn1.1} and Example \ref{exa:p-adic} satisfies conditions (i), (ii), (iii), (iv) and (v) from Definition \ref{definition.3.4}. Hence $\Upsilon_{p}$ is a quasi-bounded gentle algebra over the $p$-adic integers. We now add detail to the discussion between Equations \ref{eqn1.1} and \ref{eqn1.2}. 
\begin{remark}\label{rem}Suppose $k$ is a field and $(Q,\mathcal{J})$ satisfies gentle conditions. The \textit{completed path algebra} $\overline{kQ}$ consists of (possibly infinite) sums $\sum\lambda_{p}p$ where $\lambda_{p}\in k$ and $p$ runs once through each of the distinct paths in $Q$. Note that $\rad(\overline{kQ})$ is generated by $\mathbf{A}$ \cite[Proposition 1.2.23]{Ben2018}. Let $z$ be the sum of the cycles $\sigma\in\mathbf{P}$ which are not non-trivial powers of other cycles, and such that $\sigma^{n}\in\mathbf{P}$ for all integers $n>0$. Note $z$ acts centrally on $kQ/\mathcal{J}$, and $kQ/\mathcal{J}$ is module finite over $k[z]$ \cite[Lemmas 3.1 and 3.2]{Cra2018}. Furthermore there exists $n>0$ with $A^{n}\subseteq(z)$, so $\overline{kQ}/\overline{\mathcal{J}}$ is the $A$-adic completion of the ring $kQ/\mathcal{J}$, where $\overline{\mathcal{J}}$ is the ideal in $\overline{kQ}$ generated by $\mathcal{J}$ \cite[Proposition 3.2.7]{Ric2017} (see also \cite[Proposition 1.2.25, Lemma 1.2.27]{Ben2018}). Altogether the above shows that $\Lambda=\overline{kQ}/\overline{\mathcal{J}}$ is a quasi-bounded gentle $k[[z]]$-algebra \cite[Corollary 1.2.29]{Ben2018}. Note that $\Lambda$ is a finite-dimensional gentle algebra over $k$ if and only if there are no cycles $c\in\mathbf{P}$ as above.  
\end{remark}
In the notation from Remark \ref{rem}, note that if $(Q,\mathcal{J})$ satisfies gentle conditions then Bekkert, Drozd and Furtorny \cite{BekDroFut2009} call $\overline{kQ}/\overline{\mathcal{J}}$ a \textit{gentle algebra}. Ricke \cite{Ric2017}, who uses a more general definition, calls $\overline{kQ}/\overline{\mathcal{J}}$ a \textit{completed string algebra}.
\begin{example}\label{example.2.14-2435}\label{example235463}Let $\mathcal{J}=\langle lr,\,rq,\,tf,\,gs,\,eb,\,dc ,\,ba\rangle$
where $Q$ is the quiver 
\[
\xymatrix@R=.2em{& 0\ar[r]^{p} & 1\ar[r]^{q} & 2\ar[dd]_{r} & 3\ar[dd]_{f} & 4\ar[l]_{e} & 5\ar[l]_{b} & \\
& & & & & & \\
& 6\ar[uu]_{n} & 7\ar[l]^{m} & 8\ar[l]^{l}\ar[r]_{s} & 9\ar[uul]^{t}\ar[r]_{g} & 10\ar[uu]_{d} & 11\ar[l]^{c}\ar[uu]_{a}
}
\]
It is clear that $(Q,\mathcal{J})$ satisfies gentle conditions. Here the primitive cycles are: $tsr$ (at $2$); $edgf$ (at $3$); $dgfe$ (at $4$); $rts$ (at $8$); $srt$ and $fedg$ (at $9$); and $gfed$ (at $10$). This shows that the sum $z$ of these paths annihilates $e_{u}$ provided $u$ is $0$, $1$,  $5$, $6$, $7$ or $11$. 

Note that the $k[z]$-module $(kQ/\mathcal{J})e_{9}$ is generated by $e_{9}$, $t$, $g$, $rt$, $dg$, $srt=ze_{9}-fedg$ and $edg$. Moreover, any path in $Q$ of length $6$ lying outside $\mathcal{J}$ must have the form $\mu z$ where $\mu$ is a nontrivial path of length $2$ or $3$  lying outside $\mathcal{J}$, and so $A^{6}\subseteq (z)$.
\end{example}
\begin{remark}
There is an embedding of the $\widehat{\mathbb{Z}}_{p}$-algebra $\Upsilon_{p}$  (from Equation \ref{eqn1.1}) into the hereditary order of $2\times 2$ matrices $(\lambda_{ij})$ with $p\mid\lambda_{12}$, and so $\Upsilon_{p}$ is a \textit{nodal ring} in the sense of Burban and Drozd \cite[Definition 2.1]{BurDro2006}. Like-wise, there is an embedding of the $k[[z]]$-algebra $\overline{kQ}/\overline{\mathcal{J}}$ from Example \ref{example.2.14-2435} into the hereditary algebra $\overline{k\Gamma}$ where $\Gamma$ is the quiver 
\[
\xymatrix@R=.2em{& 0\ar[r]^{p} & 1\ar[r]^{q} & 2' & 2''\ar[dd]_{r} & & 3\ar[dd]_{f} & 4''\ar[l]_{e} &4' & 5'\ar[l]_{b} & 5'' & \\
& & & & & & \\
& 6\ar[uu]_{n} & 7\ar[l]^{m} & 8'\ar[l]^{l} & 8''\ar[r]_{s} & 9'\ar[uul]^{t} & 9''\ar[r]_{g} & 10''\ar[uu]_{d} & & 10' & 11\ar[l]^{c}\ar[uu]_{a}
}
\]
This embedding is defined by the construction from the proof of \cite[Theorem]{Zem2015}, in which Zembyk proves that any finite-dimensional gentle algebra is a \textit{nodal algebra}. Since $\overline{kQ}/\overline{\mathcal{J}}$ is infinite-dimensional and has minimal ideals, it is not nodal in either sense.
\end{remark}
\begin{remark}\label{rem2.6}
Suppose $\Lambda$ is a quasi-bounded gentle $R$-algebra. By Definition \ref{definition.3.4}(iv), for each vertex $v$ the left $\Lambda$-module $\Lambda e_{v}$ (respectively right $\Lambda$-module $e_{v}\Lambda$) has a unique maximal submodule generated by $\mathbf{A}(v\rightarrow)$ (respectively $\mathbf{A}(\rightarrow v)$). By \cite[Corollary 1.1.17]{Ben2018} this shows, for any $p\in\mathbf{P}$, that: every proper non-zero submodule of $\Lambda p$ (respectively $p\Lambda$) has the form $\Lambda qp$ (respectively $pq\Lambda$) with $q\in\mathbf{P}$; $\Lambda p$ is local with $\rad(\Lambda p)=\rad(\Lambda)p$; if $\rad(\Lambda p)\neq0$ then $\rad(\Lambda p)=\Lambda ap$ for an arrow $a$ with $ap\in\mathbf{P}$; and if $p'\in\mathbf{P}$, $\Lambda p=\Lambda p'$ and $\mathrm{f}(p)=\mathrm{f}(p')$ (respectively $p\Lambda=p'\Lambda$ and $\mathrm{l}(p)=\mathrm{l}(p')$) then $p=p'$. In particular, for any vertex $v$, the  left $\Lambda$-module $\rad(\Lambda e_{v})$ (respectively right $\Lambda$-module $\rad(e_{v}\Lambda)$) is the sum of at most two uniserial modules, and this sum is direct by Definition \ref{definition.3.4}(v).

Furthermore the elements $e_{v}$ form a complete set of orthogonal local idempotents in the ring $\Lambda$. By \cite[Corollary 1.1.25]{Ben2018} this means that the $R$-modules $k$, $\Lambda e_{v}/\rad(\Lambda e_{v})$ and $e_{v}\Lambda/\rad(e_{v}\Lambda)$ are isomorphic: and furthermore $\Lambda$ is a noetherian semilocal ring which is module finite over $R$. By the intersection theorem of Krull (see for example \cite[Exercise 4.23]{Lam1991}), if $M$ is a finitely generated $\Lambda$-module, then the intersection of $(\rad(\Lambda))^{n}M$ over all $n>0$ is trivial.
\end{remark}
\begin{notation} \cite[Definition 1.1.13]{Ben2018} Any non-trivial path $p$ in $Q$ has a \textit{first arrow} $\mathrm{f}(p)$ and a \textit{last arrow} $\mathrm{l}(p)$ satisfying $\mathrm{l}(p)p'=p=p''\mathrm{f}(p)$
for some (possibly trivial) paths $p'$ and $p''$.
\end{notation}
The proof of Corollary \ref{corollary.0.1} uses the results discussed in Remark \ref{rem2.6}.
\begin{corollary}
\emph{\cite[Corollary 1.2.14, 1.2.18 and 1.2.21]{Ben2018}}\label{corollary.0.1}\label{corollary.0.2}\label{corollary.1.20=000023} Let $v$ be a vertex, let $t>0$ be an integer and let $q\in\mathbf{P}$. Let $\Lambda$ be a quasi-bounded gentle $R$-algebra.
\begin{enumerate}
\item If $v=h(q)$, $\lambda\in\Lambda e_{v}$ and $\lambda q=0$, then $\lambda e_{v}\in\bigoplus_{a\in\mathbf{A}(v\rightarrow)}\Lambda a$.
\item If $v=t(q)$, $\lambda\in e_{v}\Lambda$ and $q\lambda=0$, then $e_{v}\lambda\in\bigoplus_{a\in\mathbf{A}(\rightarrow v)}a\Lambda$.
\item We have $\rad(\Lambda)q\cap \mathrm{l}(q)\Lambda\subseteq q\rad(\Lambda)$
and $q\rad(\Lambda)\cap\Lambda \mathrm{f}(q)\subseteq\rad(\Lambda)q$.
\item If $\sum_{p}r_{p}p\neq 0$ where $r_{p}\in R$ as $p$ runs over $\{h\in\mathbf{P}\mid h\mathrm{l}(q)\in\mathbf{P}\}$, then $\sum_{p}r_{p}pq\neq 0$.
\item If $\sum_{p}r_{p}p\neq 0$ where $r_{p}\in R$ as $p$ runs over $\{h\in\mathbf{P}\mid \mathrm{f}(q)h\in\mathbf{P}\}$, then $\sum_{p}r_{p}qp\neq 0$.
\item If $R$ is an $\mathfrak{m}$-adically complete ring then $\Lambda$ is a semiperfect ring.
\end{enumerate}
\end{corollary}
\begin{definition}\label{completegentle}
\cite[Definition 1.2.19]{Ben2018} By a \textit{complete gentle algebra} we mean a quasi-bounded gentle $R$-algebra where $R$ is $\mathfrak{m}$-adically complete.
\end{definition}
\begin{assumption}
For the remainder of the article let $\Lambda$ be a complete gentle $R$-algebra.
\end{assumption}
\section{String and band complexes.}\label{3}
 
In what follows we define a new system of words to adapt the functorial filtrations method, as in \S\ref{intro}, to classify complexes. To do so we modify the alphabet used by Bekkert and Merklen \cite{BekMer2003} to define \textit{generalised strings and bands}.
\begin{definition}
\label{definition2}\cite[Definition 1.3.26]{Ben2018} A \textit{homotopy letter $q$} is one of $\gamma$, $\gamma^{-1}$, $d_{\alpha}$,
or $d_{\alpha}^{-1}$ for $\gamma\in\mathbf{P}$ and an arrow $\alpha$. Those of the form $\gamma$ or $d_{\alpha}$ will be called \textit{direct}, and those of the form $\gamma^{-1}$ or $d_{\alpha}^{-1}$ will be called \textit{inverse}. The \textit{inverse $q^{-1}$} of a homotopy letter \textit{$q$} is defined by setting $(\gamma)^{-1}=\gamma^{-1}$, $(\gamma^{-1})^{-1}=\gamma$, $(d_{\alpha})^{-1}=d_{\alpha}^{-1}$ and $(d_{\alpha}^{-1})^{-1}=d_{\alpha}$. 

Let $I$ be one of the
sets $\{0,\dots,m\}$ (for some $m\geq0$), $\mathbb{N}$, $-\mathbb{N}=\{-n\mid n\in\mathbb{N}\}$,
or $\mathbb{Z}$. For $I\neq\{0\}$ a \textit{homotopy $I$-word}
is a sequence of homotopy letters 
\[
C=\begin{cases}
l_{1}^{-1}r_{1}\dots l_{m}^{-1}r_{m} & (\mbox{if }I=\{0,\dots,m\})\\
l_{1}^{-1}r_{1}l_{2}^{-1}r_{2}\dots & (\mbox{if }I=\mathbb{N})\\
\dots l_{-1}^{-1}r_{-1}l_{0}^{-1}r_{0} & (\mbox{if }I=-\mathbb{N})\\
\dots l_{-1}^{-1}r_{-1}l_{0}^{-1}r_{0}\mid l_{1}^{-1}r_{1}l_{2}^{-1}r_{2}\dots & (\mbox{if }I=\mathbb{Z})
\end{cases}
\]
(which will be written as $C=\dots l_{i}^{-1}r_{i}\dots$ to save
space) such that: 
\begin{enumerate}
\item any homotopy letter in $C$ of the form $l_{i}^{-1}$ (respectively $r_i$) is inverse (respectively direct);
\item any sequence of 2 consecutive letters in $C$, which is of the form $l_{i}^{-1}r_{i}$, is one of
$\gamma^{-1}d_{\mathrm{l}(\gamma)}$ or $d_{\mathrm{l}(\gamma)}^{-1}\gamma$ for some $\gamma\in\mathbf{\mathbf{P}}$; and 
\item any sequence of 4 consecutive letters in $C$ of the form $l_{i}^{-1}r_{i}l_{i+1}^{-1}r_{i+1}$ is one of 
\begin{enumerate}
\item $\gamma^{-1}d_{\mathrm{l}(\gamma)}d_{\mathrm{l}(\lambda)}^{-1}\lambda$ where $h(\gamma)=h(\lambda)$ and $\mathrm{l}(\gamma)\neq \mathrm{l}(\lambda)$;
\item $d_{\mathrm{l}(\gamma)}^{-1}\gamma d_{\mathrm{l}(\lambda)}^{-1}\lambda$ where $t(\gamma)=h(\lambda)$ and $\mathrm{f}(\gamma)\mathrm{l}(\lambda)\in\mathcal{J}$;
\item $d_{\mathrm{l}(\gamma)}^{-1}\gamma\lambda^{-1}d_{\mathrm{l}(\lambda)}$ where $t(\gamma)=t(\lambda)$ and $\mathrm{f}(\gamma)\neq \mathrm{f}(\lambda)$;
\item $\gamma^{-1}d_{\mathrm{l}(\gamma)}\lambda^{-1}d_{\mathrm{l}(\lambda)}$ where $h(\gamma)=t(\lambda)$ and $\mathrm{f}(\lambda)\mathrm{l}(\gamma)\in\mathcal{J}$.
\end{enumerate}
\end{enumerate}
For $I=\{0\}$ there are \textit{trivial homotopy words}
$ 1 _{v,1}$ and $ 1 _{v,-1}$ for each vertex
$v$.

\cite[Definition 1.3.29]{Ben2018} The head and tail of any path $\gamma\in\mathbf{P}$
are already defined and we extend this by setting
$h(d_{a}^{\pm1})=h(a)$ for any arrow $a$ and $h(q^{-1})=t(q)$ for all homotopy letters $q$. For each $i\in I$ there is an \textit{associated vertex} $v_{C}(i)$
defined by: $v_{C}(i)=t(l_{i+1})$ for $i\leq0$ and $v_{C}(i)=t(r_{i})$
for $i>0$ provided $C=\dots l_{i}^{-1}r_{i}\dots$ is non-trivial;
and $v_{ 1 _{v,\pm1}}(0)=v$ otherwise.

If $\gamma\in\mathbf{P}$ and $a=\mathrm{l}(\gamma)$ let $H(\gamma^{-1}d_{a})=-1$
and $H(d_{a}^{-1}\gamma)=1$. Let $\mu_{C}(0)=0$ and
\[
\mu_{C}(i)=\begin{cases}
H(l{}_{1}^{-1}r{}_{1})+\dots+H(l{}_{i}^{-1}r{}_{i}) & (\text{if } 0<i\in I)\\
-(H(l_{0}^{-1}r{}_{0})+\dots+H(l_{i+1}^{-1}r{}_{i+1}))) & (\text{if }0>i\in I)
\end{cases}
\]
\cite[Definition 2]{BekMer2003} (see also \cite[Definition 1.3.34]{Ben2018}). For $n\in\mathbb{Z}$ let $P^{n}(C)$ be the sum $\bigoplus\Lambda e_{v_{C}(i)}$ over $i\in \mu_{C}^{-1}(n)$. For each $i\in I$ let $ b _{i,C}$ denote the coset of $e_{v_{C}(i)}$
in $P(C)$ (in degree $\mu_{C}(i)$). If the dependency on $C$ is irrelevant let $b_{i,C}=b_{i}$. We
define the complex $P(C)$ by extending the assignment $d_{P(C)}( b _{i})= b _{i}^{-}+ b _{i}^{+}$
linearly over $\Lambda$ for each $i\in I$, where 
\[
 \begin{array}{c}
 b _{i}^{+}=\begin{Bmatrix}\alpha b _{i+1} & (\mbox{if }i+1\in I,\,l_{i+1}^{-1}r_{i+1}=d_{\mathrm{l}(\alpha)}^{-1}\alpha)\\
0 & (\mbox{otherwise})
\end{Bmatrix}\\
\\
 b _{i}^{-}=\begin{Bmatrix}\beta b _{i-1} & (\mbox{if }i-1\in I,\,l_{i}^{-1}r_{i}=\beta^{-1}d_{\mathrm{l}(\beta)})\\
0 & (\mbox{otherwise})
\end{Bmatrix}
\end{array}
\]
\end{definition}
Let $[C]_{i}=[\gamma^{-1}]$ if $l_{i}^{-1}r_{i}=d_{\mathrm{l}(\gamma)}^{-1}\gamma$
and $[C]_{i}=[\gamma]$ if $l_{i}^{-1}r_{i}=\gamma^{-1}d_{\mathrm{l}(\gamma)}$. Then $[C]=\dots [C]_{i}\dots$ defines a \textit{generalised string} or a \textit{generalised band} as in Bekkert and Merklen \cite[\S 4.1]{BekMer2003}.
\begin{example}\label{example.2.14-1}\label{example8}(See also \cite[Examples 1.3.27 and 1.3.35]{Ben2018}). Consider the finite-dimensional gentle algebra $\Lambda=kQ/\mathcal{J}$ given by $\mathcal{J}=\langle gf,\,hg,\,fh,\,sr,\,ts,\,rt\rangle$
where $Q$ is the quiver 
\[
\xymatrix@R=.2em{0 & 1\ar[l]_{w}\ar[dd]_{g} &  & 2\ar[dl]_{t} & 3\ar[l]_{x}\\
& & 4\ar[ul]_{f}\ar[dr]_{r} & & \\
5 & 6\ar[l]_{y}\ar[ur]_{h} &  & 7\ar[uu]_{s} & 8\ar[l]_{z}
}
\]
Let $C=d_{r}^{-1}rh d_{g}^{-1}g d_{f}^{-1}f r^{-1}d_{r} s^{-1}d_{s} d_{x}^{-1} x$. We may depict the string complex $P(C)$ by
\[
\xymatrix@C=1.65em@R=1.65em{ 
 &  \Lambda e_{7}\ar[dr]^{rh}\ar@{..}[l]  &  &  &  &  &  &  &  &  & P^{0}(C)\ar@{..}[lllllllll]\ar[d]^{d_{P(C)}^{0}}\\
 &  & \Lambda e_{6}\ar[dr]^{g}\ar@{..}[ll] & &  &  &  \Lambda e_{2}\ar[dl]^{s}\ar[dr]^{x}\ar@{..}[llll] &  &  &  &  P^{1}(C)\ar@{..}[llll]\ar[d]^{d_{P(C)}^{1}}\\
 &  &  & \Lambda e_{1}\ar[dr]^{f}\ar@{..}[lll] &  &  \Lambda e_{7}\ar[dl]^{r}\ar@{..}[ll]  &  & \Lambda e_{3}\ar@{..}[ll] &  &  &  P^{2}(C)\ar@{..}[lll]\ar[d]^{d_{P(C)}^{2}}\\
&  &  &  &  \Lambda e_{4}\ar@{..}[llll]  &  &  &  &  &  & P^{3}(C)\ar@{..}[llllll]
}\]
where an arrow $\Lambda e_{v}\to \Lambda e_{u}$ labelled by a path $p$ with head $v$ and tail $v$ indicates right-multiplication by $p$. The corresponding generalised string is $[C]=[(rh)^{-1}][g^{-1}][f^{-1}][r][s][x^{-1}]$.
\end{example}
\begin{definition}\cite[Definitions 1.3.26, 1.3.32 and 1.3.42]{Ben2018}\label{def.3.3}
Let $C$ be a homotopy word. Write $I_{C}$ for the subset of $\mathbb{Z}$ where $C$ is a homotopy $I_{C}$-word. Let $t\in\mathbb{Z}$. We say $C$
has \textit{controlled homogeny} if the preimage $\mu_{C}^{-1}(n)=\{i\in I_{C}\mid \mu_{C}(i)=n\}$
is a finite set for each $n\in\mathbb{Z}$.

Let $t\in\mathbb{Z}$. If $I_{C}=\mathbb{Z}$ we let $C[t]=\dots l_{t}^{-1}r_{t}\mid l_{t+1}^{-1}r_{t+1}\dots$. That is, in the language of generalised strings and bands, if $I_{C}=\mathbb{Z}$ let $[C[t]]_{i}=[C]_{i+t}$. If instead $I_{C}\neq\mathbb{Z}$ we let $C=C[t]$. 

The \textit{ inverse} $C^{-1}$ of $C$ is defined
by $( 1 _{v,\delta})^{-1}= 1 _{v,-\delta}$ if $I=\{0\}$,
and otherwise inverting the homotopy letters and reversing their order. Note the homotopy $\mathbb{Z}$-words are indexed so
that
\[
\left(\dots l_{-1}^{-1}r_{-1}l_{0}^{-1}r_{0}\mid l_{1}^{-1}r_{1}l_{2}^{-1}r_{2}\dots\right)^{-1}=\dots r_{2}^{-1}l_{2}r_{1}^{-1}l_{1}\mid r_{0}^{-1}l_{0}r_{-1}^{-1}l_{-1}\dots
\]
\end{definition}
\begin{lemma}\emph{\cite[Lemma 1.3.33, Corollary 1.3.43]{Ben2018}}
\label{lemma.4.1}\label{corollary.2.1}Let $C$ be a homotopy $I$-word and $i\in I$.
\begin{enumerate}
\item If $I=\{0,\dots,m\}$ then $v_{C^{-1}}(i)=v_{C}(m-i)$, $\mu_{C^{-1}}(i)=\mu_{C}(m-i)-\mu_{C}(m)$ and there is an isomorphism of complexes
$P(C^{-1})\rightarrow P(C)[\mu_{C}(m)]$ given by $b_{i,C}\mapsto b_{C^{-1},m-i}$.
\item If $I$ is infinite then $v_{C^{-1}}(i)=v_{C}(-i)$, $\mu_{C^{-1}}(i)=\mu_{C}(-i)$ and there is an isomorphism of complexes
$P(C^{-1})\rightarrow P(C)$ given by $b_{i,C}\mapsto b_{-i,C^{-1}}$.
\item If $I=\mathbb{Z}$ and $t\in\mathbb{Z}$ then $v_{C}(i+t)=v_{C[t]}(i)$, $\mu_{C}(i+t)=\mu_{C[t]}(i)+\mu_{C}(t)$ and there is an isomorphism of complexes
$P(C[t])\rightarrow P(C)[\mu_{C}(t)]$ given by $b_{i,C}\mapsto b_{C[t],i-t}$.
\end{enumerate}
\end{lemma}
\begin{definition}
\cite[Definitions 1.3.42 and 1.3.45]{Ben2018} We say $C$ is \textit{periodic}
if $I_{C}=\mathbb{Z}$, $C=C\left[p\right]$ and $\mu_{C}(p)=0$ for some $p>0$.  In this case the minimal such $p$ is the \textit{period} of $C$, and we say $C$ is $p$-\textit{periodic}. We say $C$ is \textit{aperiodic} if $C$ is not periodic. 

If $C$ is periodic of period $p$ then by Lemma \ref{corollary.2.1} $P^{n}(C)$ is a $\Lambda\text{-}R[T,T^{-1}]$-bimodule
where $T$ acts on the right by $ b _{i}\mapsto b _{i-p}$. By translational symmetry the map $d_{P(C)}^{n}:P^{n}(C)\rightarrow P^{n+1}(C)$ is $\Lambda\otimes_{R}R[T,T^{-1}]$-linear. For an $R[T,T^{-1}]$-module $V$ we define $P(C,V)$ by $P^{n}(C,V)=P^{n}(C)\otimes_{R[T,T^{-1}]}V$
and $d_{P(C,V)}^{n}=d_{P(C)}^{n}\otimes\mathrm{id}_{V}$ for each $n\in\mathbb{Z}$. 
\end{definition}
\begin{lemma}\label{lemma.2.3}\emph{\cite[Lemma 1.3.47]{Ben2018}} Let $n\in\mathbb{Z}$ and $C$ be a $p$-periodic homotopy word, and let $V$ be an $R[T,T^{-1}]$-module
which is free as an $R$-module.
\begin{enumerate}
\item Letting $\left\langle n,p\right\rangle =\mu_{C}^{-1}(n)\cap[0,p-1]$ gives $\mu_{C}^{-1}(n)=\{j+ps\mid j\in\left\langle n,p\right\rangle,s\in\mathbb{Z}\}$.
\item Let $L$ be the $\Lambda$-module $\bigoplus_{j\in\left\langle n,p\right\rangle }\Lambda b_{j}\otimes_{R}V$. The map $P^{n}(C)\times V\to L$ given by
\[
\begin{array}{c}
\bigg(\sum_{i\in\mu_{C}^{-1}(n)}\lambda_{i}b_{i}, v\bigg)\mapsto \sum_{j\in\left\langle n,p\right\rangle }\bigg(\sum_{s\in\mathbb{Z}}\lambda_{j+ps}b_{j}\otimes T^{-s}v\bigg)
\end{array}
\]
is $R[T,T^{-1}]$-balanced, and defines a $\Lambda$-module isomorphism $\kappa_{n}:P^{n}(C,V)\rightarrow L$.
\item The $\Lambda$-module $P^{n}(C,V)$ is projective. Moreover, if $V$ has an $R$-basis $(v_{\lambda}\mid\lambda\in\Omega)$ then the underlying $\Lambda$-module of $P(C,V)$ is generated by $\{b_{i}\otimes v_{\lambda}\mid 0\leq i\leq p-1, \lambda\in\Omega\}$.
\end{enumerate}
\end{lemma}
\begin{definition}\label{bandpres}\cite[Definition 2]{BekMer2003} If $C$ is a $p$-periodic homotopy word then $C=\dots EE\mid EE\dots$ for some unique homotopy word $E=l_{1}^{-1}r_{1}\dots l_{p}^{-1}r_{p}$. As in Lemma \ref{lemma.2.3}, let $V$  be an $R[T,T^{-1}]$-module which is a free $R$-module with basis $(v_{\lambda}\mid\lambda\in\Omega)$, and for each $n\in\mathbb{Z}$ let $\left\langle n,p\right\rangle =\mu_{C}^{-1}(n)\cap[0,p-1]$. Let $T^{\pm1}v_{\lambda}=\sum_{\mu}a^{\pm}_{\mu\lambda}v_{\mu}$ for some $a^{\pm}_{\mu\lambda}\in R$. 

Fix $n\in\mathbb{Z}$. Let $P^{n}(E,V)$ be the $\Lambda$-module $\bigoplus_{\Omega}\bigoplus_{j}\Lambda b_{j,C}$ where $j$ runs through $\left\langle n,p\right\rangle$. Let  $c_{j,\lambda}$ denote the copy of $b_{j,C}$ in the summand indexed by ($j$ and) $\lambda\in\Omega$. Define the $\Lambda$-module map $d_{P(E,V)}^{n}:P^{n}(E,V)\to P^{n+1}(E,V)$ by sending $c_{j,\lambda}\mapsto c^{+}_{j,\lambda}+c^{-}_{j,\lambda}$ where 
\[
 \begin{array}{c}
 c_{j,\lambda}^{+}=\begin{Bmatrix}\alpha c _{j+1,\lambda} & (\mbox{if }j<p-1,\,l_{j+1}^{-1}r_{j+1}=d_{\mathrm{l}(\alpha)}^{-1}\alpha)\\
\alpha(\sum_{\mu}a^{-}_{\mu\lambda} c _{0,\mu}) &  (\mbox{if }j=p-1,\,l_{p}^{-1}r_{p}=d_{\mathrm{l}(\alpha)}^{-1}\alpha)\\
0 & (\mbox{otherwise})
\end{Bmatrix}\\
\\
  c_{j,\lambda}^{-}=\begin{Bmatrix}\beta c _{j-1,\lambda} & (\mbox{if }j>0,\,l_{j}^{-1}r_{j}=\beta^{-1}d_{\mathrm{l}(\beta)})\\
\beta(\sum_{\mu}a^{+}_{\mu\lambda} c _{p-1,\mu}) &  (\mbox{if }j=0,\,l_{p}^{-1}r_{p}=\beta^{-1}d_{\mathrm{l}(\beta)})\\
0 & (\mbox{otherwise})
\end{Bmatrix}
\end{array}
\]
\end{definition}
\begin{remark}\label{indbasis}Fixing the notation from Definition \ref{bandpres}, it is straightforward to check that $P(E,V)$ defines a complex of projective modules. 
Furthermore, it is straightforward to check that the composition of the isomorphism $\kappa_{n}$ from Lemma \ref{lemma.2.3}(ii) together with the canonical isomorphism $\bigoplus\Lambda b_{j,C}\otimes_{R}V\to P^{n}(E,V)$ of $\Lambda$-modules (where $j$ runs through $\left\langle n,p\right\rangle $) defines an isomorphism of complexes $P(C,V)\to P(E,V)$. 
In particular, this means that the complex $P(E,V)$ is defined, up to isomorphism, independently of the choice of an $R$-basis for $V$. 
\end{remark}
\begin{example}\label{example11}\label{example 10}(See also \cite[Examples 1.3.41 and 1.3.46]{Ben2018}). Let $\Lambda=k[[x,y]]/(xy)$. Note that 
\[
C=\dots y^{-3}d_{y}d_{x}^{-1}x y^{-2}d_{y} d_{x}^{-1}x^{5}y^{-3}d_{y}d_{x}^{-1}x\mid y^{-2}d_{y} d_{x}^{-1}x^{5}y^{-3}d_{y}d_{x}^{-1}xy^{-2}d_{y} d_{x}^{-1}x^{5}\dots
\]
is 4-periodic. The action of $T$ on $P^{0}(C)$ and $P^{1}(C)$ may depicted by dashed arrows such as
\[
\xymatrix@R=1.2em@C=2em{&  &  &  &  &  &  &  &  &\\
\cdots \, \, \,  & \Lambda\ar[dl]_{\,y^{3}}\ar[dr]^{x} &  &
\Lambda\ar[dl]_{y^{2}}\ar[dr]^{x^{5}}\ar@{-->}@/_{1pc}/[lll] &  &
\Lambda\ar[dr]^{x}\ar[dl]_{y^{3}}\ar@{-->}@/_{1pc}/[llll] &  & \Lambda\ar[dr]^{x^{5}}\ar[dl]_{y^{2}}\ar@{-->}@/_{1pc}/[llll] & \ar@{-->}@/_{1pc}/[lll] \, \, \, \cdots & P^{0}(C)\ar[d]^{d_{P(C)}^{0}}\\
\cdots \, \, \,  &  & \Lambda &  & \Lambda\ar@{-->}@/^{1pc}/[llll]  &  & \Lambda\ar@{-->}@/^{1pc}/[llll] &  & \ar@{-->}@/^{1pc}/[llll]\, \, \,\cdots &  P^{1}(C)\\
&  &  &  &  &  &  &  &  &\\}
\]
Define the $k[[t]][T,T^{-1}]$-module $V=k[[t]]\oplus k[[t]]$ by $
T(f(t),g(t))=(f(t)-tg(t),tf(t)+\lambda g(t))$ where $0\neq\lambda\in k$. 
Note that $P(E,V)\simeq P(C,V)$ where $P(E,V)$ is depicted by
\[
\xymatrix@R=0.5em@C=0.5em{&
\Lambda\oplus\Lambda\ar@/_{1pc}/[ddl]^{x}_(0.7){A=\small{\begin{pmatrix}1 & -t\\
t & \lambda
\end{pmatrix}}^{-1}}\ar[ddrrrr]|>>>>>>>>>>>>>>>>{\hole}_(0.55){y^{3}} & & &
\Lambda\oplus\Lambda\ar@/^{0.55pc}/[ddllll]^(0.5){y^{2}}\ar@/^{1pc}/[ddr]^{x^{5}} &  & & & P^{-1}(E,V)\ar[dd]^{d_{P(E,V)}^{-1}}\\
&  &  &  &  &  &  & & &\\
\Lambda\oplus\Lambda &  &  &  &  &  \Lambda\oplus\Lambda & & & P^{0}(E,V) 
}
\]
where $E= y^{-2}d_{y} d_{x}^{-1}x^{5}y^{-3}d_{y}d_{x}^{-1}x$ (so that $C={}^{\infty}E{}^{\,\infty}$). Note that the matrix $A^{-1}$ defines the action of $T$ on $V$. Furthermore $A^{-1}$ has determinant  $\lambda+t^{2}$, and the inverse of  $\lambda+t^{2}$ is $\sum_{n=0}^{\infty}(-1)^{n}\lambda^{-(n+1)}t^{2n}\in k[[t]]$.
\end{example}
\begin{lemma}\label{bandsgood}Let $C$ be a $p$-periodic homotopy word of the form $C={}^{\infty}E{}^{\,\infty}$ where $E=l_{1}^{-1}r_{1}\dots l_{p}^{-1}r_{p}$ is a homotopy word. Suppose $l_{p}^{-1}r_{p}l_{1}^{-1}r_{1}=d_{\mathrm{l}(\alpha)}^{-1}\alpha\beta^{-1} d_{\mathrm{l}(\beta)}$ for some $\alpha,\beta\in\mathbf{P}$. Let $U$ and $V$ be $R[T,T^{-1}]$-modules which are free as $R$-modules. If $k\otimes_{R} U\simeq k\otimes V$ as $k[T,T^{-1}]$-modules then $P(E,U)\simeq P(E,V)$ as complexes.
\end{lemma}
Before giving the proof of Lemma \ref{bandsgood} we recall why, over a semiperfect ring, any homomorphism between projective modules, which becomes an isomorphism upon factoring the radicals, must have been an isomorphism. 
\begin{remark}\label{projiso} Let $\Gamma$ be a semiperfect ring, and so $\Gamma=\bigoplus_{s}\Gamma e_{s}=\bigoplus_{s}e_{s}\Gamma$  for some finite set $E=\{e_{s}\mid s\in S\}$
of idempotents where the left (respectively right) $\Gamma$-module $\Gamma e$ (respectively $e\Gamma$) is local with unique maximal submodule $\mathrm{rad}(\Gamma)e$ (respectively $e\mathrm{rad}(\Gamma)$). A $\Gamma$-module is said to be \textit{quasi-free} if it is a direct sum of modules of the form $\Gamma e_{s}$. We write $\Gamma\text{-}\mathbf{quas}$ for the full subcategory of $\Gamma\text{-}\mathbf{Mod}$ consisting of finitely generated quasi-free modules. We say a full subcategory $\mathcal{A}$ of $\Gamma\text{-}\mathbf{Mod}$ of all $\Gamma$-modules \textit{reflects} \textit{monomorphisms} (respectively \textit{epimorphisms}, respectively \textit{isomorphisms}) \textit{modulo the radical} provided that for any homomorphism $f:N\rightarrow L$ in $\mathcal{A}$, if the induced morphism $\bar{f}:N/\mathrm{rad}(N)\to L/\mathrm{rad}(L)$ is a monomorphism (respectively epimorphism, respectively isomorphism) then $f$ is a monomorphism (respectively epimorphism, respectively isomorphism). Note that every quasi-free $\Gamma$-module is a union of its finitely-generated quasi-free summands \cite[Lemma 3.1.40]{Ben2018}. So, by \cite[Lemma 3.1.39]{Ben2018}, if $\Gamma\text{-}\mathbf{quas}$ reflects monomorphisms (respectively epimorphisms, respectively isomorphisms) modulo the radical, then so does the full subcategory $\Gamma\text{-}\mathbf{Proj}$ of $\Gamma\text{-}\mathbf{Mod}$ consisting of projective modules.

We now explain why $\Gamma\text{-}\mathbf{Proj}$ reflects isomorphisms modulo the radical. By the above, assuming $f:N\rightarrow L$ is a homomorphism in $\Gamma\text{-}\mathbf{quas}$ such that $\bar{f}:N/\mathrm{rad}(N)\to L/\mathrm{rad}(L)$ is an isomorphism, it suffices to show $f$ is an isomorphism. Write $h$ for the inverse of $\bar{f}$. Since $N$ and $L$ are quasi-free the canonical epimorphisms $N\to N/\mathrm{rad}(N)$ and $L\to L/\mathrm{rad}(L)$ are projective covers. Since $\bar{f}$ is a section, and since $N$ and $L$ are (finitely generated and) projective, $f$ is a section, by (the proof of) \cite[Lemma 2.2]{Jon1976}. Write $g$ for the corresponding retraction. Since $gf$ is the identity on $N$, $\bar{g}\bar{f}$ is the identity on $N/\mathrm{rad}(N)$, and so $h=\bar{g}$ (note inverses of isomorphisms are unique). This means $\bar{g}$ is an isomorphism, and so as above $g$ is a section. Again, since inverses are unique, this means $f$ is an isomorphism. 
\end{remark}
\begin{proof}[of Lemma \ref{bandsgood}] Let $\bar{\tau} :k\otimes_{R} U\to k\otimes V$ be the said isomorphism of $k[T,T^{-1}]$-modules, that is, an isomorphism of $k$-vector spaces which respects the action of $T$. For any $R$-basis $(w_{\lambda}\mid \lambda\in\Omega)$ of a free $R$-module $W$ let $\bar{w}_{\lambda}=1\otimes_{R} w_{\lambda}$ for each $\lambda$, which defines a $k$-basis $(\bar{w}_{\lambda}\mid \lambda\in\Omega)$ of $k\otimes_{R} W$. If $W$ is also an $R[T,T^{-1}]$-module, then by Remark \ref{indbasis} any choice of an $R$-basis defines $P(E,W)$ up to isomorphism. 
Altogether we can choose $R$-bases $(u_{\lambda}\mid\lambda\in\Omega)$ and $(v_{\lambda}\mid\lambda\in\Omega)$ of $U$ and $V$ respectively such that $\bar{\tau} (T^{-1}\bar{u}_{\lambda})=T^{-1}\bar{v}_{\lambda}$ for each $\lambda$. 

 Let $c_{j,\lambda}$ (respectively $e_{j,\lambda}$) denote the copy of $b_{j,C}$ in the summand of $P^{n}(E,U)$ (respectively $P^{n}(E,V)$) indexed by $j$ and $\lambda$. Let $\phi(c_{j,\lambda})=e_{j,\lambda}$ for all $\lambda$ and all $j$ with $0<j\leq p-1$.  To define a morphism  $\phi:P(E,U)\to P(E,V)$ of complexes it suffices to define $\phi(c_{0,\lambda})\in P^{0}(E,V)$ for all $\lambda$, and explain why extending these assignments (linearly over $\Lambda$) satisfies $\phi(d_{P(E,U)}(c_{j,\lambda}))=d_{P(E,V)}(\phi(c_{j,\lambda}))$ for all $j$ and $\lambda$. We begin by defining $\phi(c_{0,\lambda})$.

For each $\lambda$ let $\bar{\tau} (T^{\pm1}\bar{u}_{\lambda})=\sum_{\mu}s^{\pm}_{\mu\lambda}\bar{v}_{\mu}$ and $T^{\pm1}\bar{v}_{\lambda}=\sum_{\mu}r^{\pm}_{\mu\lambda}\bar{v}_{\mu}$ for some $s^{\pm}_{\mu\lambda},r^{\pm}_{\mu\lambda}\in R$. Since $\bar{\tau} (T^{-1}\bar{u}_{\lambda})=T^{-1}\bar{v}_{\lambda}$ for each $\lambda$, this means that $\sum_{\mu}(s^{-}_{\mu\lambda}-r^{-} _{\mu\lambda})v_{\mu}=\sum_{\mu}z _{\mu\lambda}v_{\mu}$ for some $z _{\mu\lambda}\in\mathfrak{m}$. Since $V$ has $R$-basis $(v_{\lambda}\mid\lambda\in\Omega)$ this means $s^{-}_{\mu\lambda}-r^{-}_{\mu\lambda}=z _{\mu\lambda}$ for all $\mu,\lambda$. Since $EE$ is a homotopy word, $t(E)=h(E)$ and $-s(E^{-1})=s(E)$, which means $t(\alpha)=t(\beta)$ and $\mathrm{f}(\alpha)\neq\mathrm{f}(\beta)$. Let $w=t(\alpha)$. Since $\Lambda$ is rad-nilpotent modulo $\mathfrak{m}$ we have $e_{w}z_{\mu\lambda} \in \sum a\Lambda$ where the sum runs over $a\in\mathbf{A}(\rightarrow w)$. Without loss of generality we assume $\mathbf{A}(\rightarrow w)=\{x,y\}$ where $x\neq y$. 

Hence for all $\mu,\lambda\in\Omega$ there exists some $\gamma_{\mu\lambda,x},\gamma_{\mu\lambda,y}\in\Lambda$ such that $e_{w}z_{\mu\lambda} =x\gamma_{\mu\lambda,x} +y\gamma_{\mu\lambda,y} $. Without loss of generality $\alpha x,\beta y\notin\mathbf{P}\ni \beta x,\alpha y$.  For all $\lambda\in\Omega$ let $\delta_{\lambda\lambda}=1$ and $\delta_{\mu\lambda}=0$ for all $\mu\in\Omega$ with $\mu\neq\lambda$. For each $\mu\in\Omega$ extend the assingment $\phi(c_{0,\mu})=\sum_{\eta}(\delta_{\mu\eta} +r^{+}_{\eta\mu}y\gamma_{\mu\eta,y})e_{0,\mu}$ linearly over $\Lambda$, where the sum runs over $\eta\in\Omega$. We now check $\phi(d_{P(E,U)}(c_{j,\lambda}))=d_{P(E,V)}(\phi(c_{j,\lambda}))$ for all $j$ and all $\lambda$. By our assumption on the form of the homotopy word $E$ we have that: $c_{0,\lambda}^{-}=0$ and $e_{0,\lambda}^{-}=0$;  if $j>1$ then $\phi(c^{-}_{j,\lambda})=e^{-}_{j,\lambda}$; and if $j< p-1$ then $\phi(c^{+}_{j,\lambda})=e^{+}_{j,\lambda}$. In case $j=1$ then
\[\begin{array}{c}
\phi(c^{-}_{1,\lambda})=\phi(\beta c _{0,\lambda}))=\beta\phi( c _{0,\lambda})=\beta\phi( c _{0,\lambda})=\beta(\sum_{\eta\in\Omega}(\delta_{\lambda\eta} +r^{+}_{\eta\lambda}y\gamma_{\lambda\eta,y})e_{0,\lambda})\\
=\sum_{\eta}\beta\delta_{\lambda\eta}e_{0,\lambda} +\sum_{\eta}\beta r^{+}_{\eta\lambda}y\gamma_{\lambda\eta,y}e_{0,\lambda}=\beta\delta_{\lambda\lambda}e_{0,\lambda} + \sum_{\eta}r^{+}_{\eta\lambda}\beta y\gamma_{\lambda\eta,y}e_{0,\lambda}=\beta e_{0,\lambda}=e^{-}_{1,\lambda}
\end{array}
\]
In case $j=p-1$ then
\[\begin{array}{c}
\phi(c^{+}_{p-1,\lambda})=\phi(\alpha(\sum_{\mu}r^{-}_{\mu\lambda} c _{0,\mu}))=\alpha\sum_{\mu}r^{-}_{\mu\lambda} \phi(c _{0,\mu})\\
=\alpha(\sum_{\mu}r^{-}_{\mu\lambda} (\sum_{\eta} (\delta_{\mu\eta} +r^{+}_{\eta\mu}y\gamma_{\mu\eta,y})e_{0,\mu})=\sum_{\mu,\eta}\alpha r^{-}_{\mu\lambda}\delta_{\mu\eta}e_{0,\mu}+\sum_{\mu,\eta}\alpha r^{-}_{\mu\lambda}r^{+}_{\eta\mu}y\gamma_{\mu\eta,y}e_{0,\mu}\\
=\sum_{\mu}\alpha r^{-}_{\mu\lambda}e_{0,\mu}+\sum_{\mu,\eta}\alpha \delta_{\eta\lambda}y\gamma_{\mu\eta,y}e_{0,\mu}=\sum_{\mu}(\alpha r^{-}_{\mu\lambda}e_{0,\mu}+\alpha y\gamma_{\mu\lambda,y}e_{0,\mu})\\
=\sum_{\mu}(\alpha r^{-}_{\mu\lambda}e_{0,\mu}+\alpha (y\gamma_{\mu\lambda,y}+x\gamma_{\mu\lambda,x})e_{0,\mu})=\sum_{\mu}(\alpha r^{-}_{\mu\lambda}e_{0,\mu}+\alpha z_{\mu\lambda}e_{0,\mu})\\
=\alpha(\sum_{\mu}( r^{-}_{\mu\lambda}+ z_{\mu\lambda})e_{0,\mu})
=\alpha(\sum_{\mu}s^{-}_{\mu\lambda}e_{0,\mu})=e^{+}_{p-1,\lambda}.
\end{array}
\]
Altogether we have shown that $\phi(c^{\pm}_{j,\lambda})=e^{\pm}_{j,\lambda}$, and so $\phi$ defines a morphism of complexes. It suffices to explain why each map $\phi^{n}$ is an isomorphism of $\Lambda$-modules. By construction, and by Remark \ref{projiso}, to show each $\phi^{n}$ is an isomorphism it will be enough to explain why the induced morphism $\bar{\phi}^{n}:P^{n}(E,U)/\mathrm{rad}(P^{n}(E,U)) \to P^{n}(E,V)/\mathrm{rad}(P^{n}(E,V))$ is an isomorphism. By construction, $\bar{\phi}^{n}$ maps $c_{j,\lambda}+\mathrm{rad}(P^{n}(E,U))$ to $e_{j,\lambda}+\mathrm{rad}(P^{n}(E,U))$ (for all appropriate $j$ and $n$), which is clearly an isomorphism.
\end{proof}
\begin{definition}\cite[Definition 1.3.48]{Ben2018} A \textit{string complex} has the form $P(C)$ where $C$ is aperiodic. If $V$ is an $R[T,T^{-1}]$-module we call $P(C,V)$ a \textit{band complex} provided $C$ is a periodic homotopy $\mathbb{Z}$-word, $V$ is an indecomposable
$R[T,T^{-1}]$-module and $V$ is free as an $R$-module.
\end{definition}
\section{How the classification works.}\label{structuretheorem}
\begin{notation}
If $\mathcal{A}$ is an abelian category let $\mathcal{P}_{\mathcal{A}}$ be the full subcategory of $\mathcal{A}$ consisting of the projective objects and let $\mathcal{C}(\mathcal{P}_{\mathcal{A}})$ be the category of complexes in $\mathcal{P}_{\mathcal{A}}$. We say $\mathcal{A}$ has  \textit{enough radicals} if any object $X$ of $\mathcal{A}$ has a set of maximal subobjects,  whose infimum, the \textit{radical} $\rad(X)$ of $X$, exists in $\mathcal{A}$. If $\mathcal{A}$ is an abelian category with enough radicals let $\mathcal{C}_{\mathrm{min}}(\mathcal{P}_{\mathcal{A}})$
and $\mathcal{K}_{\mathrm{min}}(\mathcal{P}_{\mathcal{A}})$ be the full subcategories
of $\mathcal{C}(\mathcal{P}_{\mathcal{A}})$ and $\mathcal{K}(\mathcal{P}_{\mathcal{A}})$ consisting of \textit{homotopically minimal} complexes: objects $M$ in $\mathcal{C}(\mathcal{P}_{\mathcal{A}})$ such that $\im(d_{M}^{n})\subseteq\rad(M^{n+1})$
for all $n\in\mathbb{Z}$. Here let $\Xi_{\mathcal{A}}:\mathcal{C}_{\mathrm{min}}(\mathcal{P}_{\mathcal{A}})\rightarrow\mathcal{K}_{\mathrm{min}}(\mathcal{P}_{\mathcal{A}})$ be the restriction of the quotient functor $\mathcal{C}(\mathcal{P}_{\mathcal{A}})\rightarrow\mathcal{K}(\mathcal{P}_{\mathcal{A}})$.
\end{notation}
\begin{assumption} In \S\ref{structuretheorem} we let: $\mathcal{M}$ be an abelian category with enough radicals and small coproducts; $\mathcal{N}$ be an abelian subcategory of $\mathcal{M}$; $\mathfrak{I}$ be an index set; and $S_{i}:\mathcal{A}_{i}\rightarrow\mathcal{C}_{\mathrm{min}}(\mathcal{P}_{\mathcal{M}})$
and $F_{i}:\mathcal{K}_{\mathrm{min}}(\mathcal{P}_{\mathcal{M}})\rightarrow\mathcal{X}_{i}$ be additive functors for all $i\in\mathfrak{I}$. \end{assumption}
To prove Theorems \ref{theorem.1.1}, \ref{theorem.1.2} and \ref{theorem.1.3} we choose $\mathcal{M}=\Lambda\text{-}\mathbf{Mod}$ and $\mathcal{N}=\Lambda\text{-}\mathbf{mod}$. Note $\mathcal{N}$ has enough projective covers by Corollary \ref{corollary.0.2}(vi). In Lemma \ref{lemma} we adapt an interpretation by Ringel \cite[\S 3, p. 22, Lemma]{Rin1975} of a result by Gabriel \cite[\S 4, Structure theorem]{Gab1987}. 
\begin{corollary}
\label{corollary.9.1}\emph{\cite[Corollary 3.2.25]{Ben2018} (}see also \emph{\cite[Proposition B.2]{Kra2005})}. If $\mathcal{N}$ has enough projective covers then the subcategory $\mathcal{K}_{\mathrm{min}}(\mathcal{P}_{\mathcal{N}})$ of $\mathcal{K}(\mathcal{P}_{\mathcal{N}})$ is dense.
\end{corollary}
We now provide the categorical blueprint that was discussed at the end of \S\ref{intro}.
\begin{definition}\label{detect}\cite[Definition 2.6.4]{Ben2018} 
Suppose $\mathcal{N}$ has enough projective covers. We say that the collection  $\{(S_{i},F_{i})\mid i\in\mathfrak{I}\}$ \textit{detects the
objects} in $\mathcal{K}(\mathcal{P}_{\mathcal{N}})$
if the following statements hold. 
\begin{enumerate}
\item For any $i\in\mathfrak{I}$:
\begin{enumerate}
\item[(FFI)] the functor $F_{i}\Xi_\mathcal{M} S_{i}$ is dense and reflects isomorphisms;
\item[(FFII)] $F_{j}\Xi _\mathcal{M}S_{i}\simeq0$ for each $j\in\mathfrak{I}$
with $j\neq i$;
\item[(FFIII)] $F_{i}$ preserves small coproducts; and
\item[(FFIV)] for each object $M$ in $\mathcal{K}_{\mathrm{min}}(\mathcal{P}_{\mathcal{N}})$ there exists an object $A_{i,M}$ in $\mathcal{A}_{i}$ and a morphism $\gamma_{i,M}:\Xi_\mathcal{N}(S_{i}(A_{i,M}))\rightarrow M$ in $\mathcal{K}_{\mathrm{min}}(\mathcal{P}_{\mathcal{M}})$ such that $F_{i}(\gamma_{i,M})$ is an isomorphism.
\end{enumerate}
\item For all morphisms $\theta:N\rightarrow M$ in $\mathcal{C}_{\mathrm{min}}(\mathcal{P}_{\mathcal{M}})$:
\begin{enumerate}
\item[(FFV)] if $M$ lies in $\mathcal{C}_{\mathrm{min}}(\mathcal{P}_{\mathcal{N}})$
and $F_{i}(\Xi_\mathcal{M}(\theta))$ is epic for all $i\in\mathcal{I}$ then each
$\theta^{n}$ is epic; and
\item[(FFVI)] if $N=\bigoplus_{i\in\mathfrak{I}}S_{i}(A_{i})$ and $F_{i}(\Xi_\mathcal{M}(\theta))$ is monic
for each $i\in\mathcal{I}$ then each $\theta^{n}$ is monic.
\end{enumerate}
\end{enumerate}
\end{definition}
\begin{lemma}
\label{lemma}\emph{\cite[Lemma 2.6.5]{Ben2018}} Suppose that $\mathcal{N}$ has enough projective covers, and that $\{(S_{i},F_{i})\mid i\in\mathfrak{I}\}$ detects
the objects in $\mathcal{K}(\mathcal{P}_{\mathcal{N}})$.
\begin{enumerate}
\item Any object $N$ of $\mathcal{K}(\mathcal{P}_{\mathcal{N}})$ is isomorphic to $\bigoplus_{i\in\mathfrak{I}}\Xi_\mathcal{N}(S_{i}(A_{i,M}))$ for some  $M$ in $\mathcal{K}_{\mathrm{min}}(\mathcal{P}_{\mathcal{N}})$.
\item If $M$ is an indecomposable object of $\mathcal{K}(\mathcal{P}_{\mathcal{N}})$ then there is some $i\in\mathfrak{I}$ and some indecomposable object $A$ of $\mathcal{A}_{i}$ such that $M\simeq\Xi_\mathcal{N}(S_{i}(A))$.
\item If $i\in\mathfrak{I}$ and $A$ is an indecomposable object of $\mathcal{A}_{i}$ then $\Xi_\mathcal{M}(S_{i}(A))$ is indecomposable.
\end{enumerate}
\end{lemma}
\begin{proof}(i) Since $\mathcal{N}$ has projective covers we have $N\simeq M$ for some object $M$ of $\mathcal{K}_{\mathrm{min}}(\mathcal{P}_{\mathcal{N}})$ by Corollary \ref{corollary.9.1}. Let $N_{i}=\Xi_{\mathcal{M}}(S_{i}(A_{i,M}))$
for each $i\in\mathfrak{I}$. Fix $l\in\mathfrak{I}$ and let $\iota_{l}:N_{l}\rightarrow\bigoplus_{i\in\mathfrak{I}}N_{i}$
be the monic in the coproduct. For each $i$ let $g_{i}$ be the morphism defined by the identity on $N_{l}$ if $i=l$, and $N_{i}\rightarrow 0$ otherwise. By the universal property of the coproduct there exists $\pi_{l}:\bigoplus_{i\in\mathfrak{I}}N_{i}\rightarrow N_{l}$ with $\pi_{l}\iota_{l}=\mathrm{id}$. 

By FFIII, for any set $A$ and any collection $X=\{X_{a}\mid a\in A\}$ of objects in $\mathcal{K}_{\mathrm{min}}(\mathcal{P}_{\mathcal{M}})$ there is an isomorphism $\sigma_{X}:\bigoplus_{a}F_{l}(X_{a})\to F(\bigoplus_{a}X_{a})$ such that $F_{l}(\bigoplus _{a}f_{a})\sigma_{X}=\sigma_{Y}\bigoplus _{a}F_{l}(f_{a})$ for any collection of morphisms $\{f_{a}:X_{a}\to Y_{a}\mid a\in A\}$. Let $A=\mathfrak{I}$, and for each $i\in A$ let $f_{i}$ be the morphism defined by the identity on $N_{l}$ if $i=l$, and $0\rightarrow N_{i}$ otherwise.  
By FFII and the above there are isomorphisms $\sigma_{X}:F_{l}(N_{l})\rightarrow F_{l}(N_{l})$
and $\sigma_{Y}:F_{l}(\bigoplus_{i\in\mathfrak{I}}N_{i})\rightarrow F_{l}(N_{l})$
with $\sigma_{Y}F_{l}(\iota_{l})=F_{l}(\mathrm{id})\sigma_{X}$, so $F_{l}(\iota_{l})F_{l}(\pi_{l})=\mathrm{id}$.

By FFIV and the coproduct universal property there exists $\theta:\Xi_{\mathcal{M}}(\bigoplus_{i\in\mathfrak{I}}S_{i}(A_{i,M}))\rightarrow M$ in $\mathcal{K}_{\mathrm{min}}(\mathcal{P}_{\mathcal{M}})$ 
satisfying $\theta\iota_{i}=\gamma_{i,M}$ for each $i\in\mathfrak{I}$.
Since $\Xi_{\mathcal{N}}$ is dense there is an object $L$ in $\mathcal{C}_{\mathrm{min}}(\mathcal{P}_{\mathcal{N}})$
and an isomorphism $\psi:\Xi_{\mathcal{N}}(L)\rightarrow M$ in $\mathcal{K}_{\mathrm{min}}(\mathcal{P}_{\mathcal{N}})$. Since $\Xi_{\mathcal{M}}$ is full we have a morphism $\varphi:\bigoplus_{i}S_{i}(A_{i,M})\rightarrow L$
in $\mathcal{C}_{\mathrm{min}}(\mathcal{P}_{\mathcal{M}})$ with
$\Xi(\varphi)=\psi^{-1}\theta$.  

By the above we have $F_{i}(\theta)=F_{i}(\gamma_{i,M})F_{i}(\pi_{i})$
which is an isomorphism for each $i$ by previous paragraph, and so $F_{i}(\Xi_{\mathcal{M}}(\varphi))$
is an isomorphism for each $i\in\mathfrak{I}$. By FFV and FFVI this shows $\varphi^{n}$ is an isomorphism for each $n\in\mathbb{Z}$, so $\theta$ is an isomorphism. As $\Xi_{\mathcal{M}}$ preserves small
coproducts, altogether we have $M\simeq\bigoplus_{i\in\mathfrak{I}}\Xi_{\mathcal{M}}(S_{i}(A_{i,M}))$.

(ii) By (i) $M\simeq\bigoplus_{i\in\mathfrak{I}}\Xi_{\mathcal{N}}(S_{i}(A_{i,M}))$
and so $\Xi_{\mathcal{N}}(S_{i}(A_{i,M}))=0$ apart from when $i=t$ for some $t\in\mathfrak{I}$.
Hence $M\simeq\Xi_{\mathcal{N}}(S_{t}(A_{t,M}))$. Given objects $X$ and $Y$ of $\mathcal{A}_{t}$
for which $A_{t,M}=X\oplus Y$ we have $\Xi_{\mathcal{M}}(S_{t}(X))=0$ without loss of generality. This means $F_{t}(\Xi_{\mathcal{M}}(S_{t}(\mathbf{0})))$
is an isomorphism where $\mathbf{0}:X\rightarrow0$
in $\mathcal{A}_{t}$. Since $F_{t}\Xi_{\mathcal{M}} S_{t}$ reflects isomorphisms
by FFI, $X=0$.

(iii) If $F_{i}(\Xi_{\mathcal{M}}(S_{i}(A))=0$ then $A=0$ since $F_{i}\Xi_{\mathcal{M}} S_{i}$
is dense and reflects isomorphisms by FFI. Hence $\Xi_{\mathcal{M}}(S_{i}(A))\neq0$,
and we suppose $\Xi_{\mathcal{M}}(S_{i}(A))=X'\oplus Y'$ for objects $X'$ and
$Y'$ of $\mathcal{K}_{\mathrm{min}}(\mathcal{P}_{\mathcal{M}})$. So there are objects $X$ and $Y$
in $\mathcal{C}_{\mathrm{min}}(\mathcal{P}_{\mathcal{M}})$ with
$S_{i}(A)=X\oplus Y$, $\Xi(X)=X'$ and $\Xi(Y)=Y'$. Hence, without loss of generality, $F_{i}(X')=0$.

Label the monics of the coproduct $\iota_{X}:X\rightarrow X\oplus Y$ and $\iota_{Y}:Y\rightarrow X\oplus Y$,
and epics of the product $\pi_{X}:X\oplus Y\rightarrow X$ and $\pi_{Y}:X\oplus Y\rightarrow Y$. For $j\in\mathfrak{I}$
with $j\neq i$ we have $F_{j}(X')\oplus F_{j}(Y')=F_{j}(\Xi(S_{i}(A)))=0$
by FFII. Hence for any $j$ we have $F_{j}(\Xi_{\mathcal{M}}(\iota_{X}))=0$ and so $F_{j}(\Xi_{\mathcal{M}}(\pi_{Y}))$ is an isomorphism with inverse $F_{j}(\Xi_{\mathcal{M}}(\iota_{Y}))$.
Since $\pi_{Y}$ is an morphism in $\mathcal{C}_{\mathrm{min}}(\mathcal{P}_{\mathcal{M}})$ of the form $S_{i}(A)\rightarrow Y$ this means $\pi_{Y}^{n}$ is monic for each $n\in\mathbb{Z}$
by FFVI. Since each $\pi_{Y}^{n}$ is also a split epic, we have $X'=0$.
\end{proof}
In \S\ref{Refined Functors for Complexes.} and \S\ref{subsec:Natural-Isomorphisms.} we define functors $S_{B,D,n}$ and $F_{B,D,n}$ for each integer $n$ and certain pairs $B,D$ of homotopy words. This is done in a way so that the collection of $(S_{B,D,n},F_{B,D,n})$ detects the objects in $\mathcal{K}(\Lambda\text{-}\mathbf{proj})$ (see Proposition \ref{proposition6}).
\section{\label{sec:Linear-Relations-and}Linear relations and reductions}
 In previous applications of the functorial filtrations method to module classifications, so called \textit{refined functors} were defined as quotients of subfunctors of the forgetful
functor (sending a module over a $k$-algebra to its underlying vector space). These
subfunctors were constructed using the language of \textit{linear relations} studied by Mac Lane \cite[\S 2]{Mac1961}. To adapt the functorial filtrations method to classify complexes we
use a similar approach.

Given $R$-modules $L$ and $M$ an \emph{$R$-linear relation from} $L$ \emph{to} $M$ (or \emph{on} $M$ if $L=M$) is an $R$-submodule $C$ of the direct sum $L\oplus M$. We will often write \emph{relation} to mean $R$-linear relation. This notion generalises the graph of an $R$-linear map $L\rightarrow M$. Let $Cm = \{ n\in M \colon (m,n)\in C \}$ for any $m\in M$, and for a subset $S\subseteq M$ let $CS$ be the union $ \bigcup Cs$ over $s\in S$. When $C$ is the graph of a map $f$ then $CS$ is the image of $S$ under $f$. If $W$ is a relation from a module $L$
to $M$ the \textit{composition} $CW$ is the relation from $L$ to
$N$ given by pairs $(l,n)\in L\oplus N$ such that $(m,n)\in C$
and $(l,m)\in W$ for some $m\in M$. Now suppose $C$ is a relation on $M$. For any integer $n>0$ let $C^{n}$ be the $n$-fold composition
of $C$ with itself (so $C^{1}=C$ and $C^{2}=CC$). Let $C^{\sharp}=C''\cap(C^{-1})'', C^{\flat}=C''\cap(C^{-1})'+C'\cap(C^{-1})''$ where, for a relation $D$ on $M$,
\[
\begin{array}{cc}
D'=\bigcup_{n>0}U^{n}0, & D''=\{m\in M\mid\exists m_{0},m_{1},m_{2},\dots\in M\,:\,m_{0}=m,\,m_{i}\in Dm_{i+1} \, \, \, \forall i\}\end{array}
\]
The following result from \cite{Cra2018} was written only in the context where $R$ is a field. The proof does not make use of this assumption, and generalises with no complication.
\begin{lemma}
\label{lemma.4.4.}\label{lemma.2.10}
\emph{\cite[Lemmas 4.4 and 4.5]{Cra2018}} For any relation $C$ on $M$ we have: $C^{\sharp}\subseteq CC^{\sharp}$; $C^{\flat}=C^{\sharp}\cap CC^{\flat}$; $C^{\sharp}\subseteq C^{-1}C^{\sharp}$; $C^{\flat}=C^{\sharp}\cap C^{-1}C^{\flat}$. Consequently there is an $R$-module automorphism $\theta$ on $C^{\sharp}/C{}^{\flat}$
defined by $\theta(m+C{}^{\flat})=m'+C{}^{\flat}$ if and only if $m'\in C^{\sharp}\cap(C^{\flat}+Cm)$.
\end{lemma}

\begin{definition}\label{definition.1.4.32}\cite[Definition 1.4.32]{Ben2018} A \textit{reduction
of} a relation $C$ on an $R$-module $M$ is a pair $(U,f)$ where: $U$ is an $R[T,T^{-1}]$-module;
$U$ is free over $R$; $f:U\rightarrow M$ is $R$-linear; $C^{\sharp}=\im(f)+C^{\flat}$; and $f(Tu)\in Cf(u)$
for all $u\in U$. We say a reduction $(U,f)$ of $C$ \textit{meets in
}$\mathfrak{m}$ if $\{u\in U:f(u)\in C^{\flat}\}\subseteq\mathfrak{m}U$. 
\end{definition}
When $R$ is a field, if $(U,f)$ is a reduction of a relation $C$ on $M$ which meets in $0$, then by \cite[Corollary 1.4.33]{Ben2018} $C$ is \textit{split} in the sense of \cite[p. 9]{Cra2018}: that is, there is an $R$-linear subspace $W$ of $M$ such that $C^{\sharp}=W\oplus C^{\flat}$ and $\# Cm\cap W=1$ for each $m\in W$. In our setting $C^{\flat}$ need not have an complement
as an $R$-submodule of $C^{\sharp}$ (see \cite[Example 1.4.35]{Ben2018}). So, the following is a generalisation of \cite[Lemma 4.6]{Cra2018}. 
\begin{lemma}
\label{lemma.2.11}\emph{\cite[Lemma 1.4.34]{Ben2018}} Let
$M$ be an $R$-module and $C$ a relation on $M$ such that $C^{\sharp}/C^{\flat}$
is a finite-dimensional $R/\mathfrak{m}=k$-vector space. Then there
is a reduction $(U,f)$  of $C$ which meets in $\mathfrak{m}$ where $U$ has free $R$-rank $\mathrm{dim}_{k}(C^{\sharp}/C^{\flat})$.
\end{lemma}
\begin{proof} Let $\theta$ denote the induced $R$-module automorphism of $C^{\sharp}/C^{\flat}$
from Lemma \ref{lemma.2.10}. Let $\overline{A}=(\overline{a_{ij}})$ be the
matrix of $\theta$ (with entries from $k$) with respect to a $k$-basis
$\overline{v_{1}},\dots,\overline{v_{d}}$ of $C^{\sharp}/C^{\flat}$.
For each $i$ choose $v_{i}\in C^{\sharp}$ such that $\overline{v_{i}}=v_{i}+C^{\flat}$
and for each $j$ choose $a_{ij}\in R$ such that $\overline{a_{ij}}=a_{ij}+\mathfrak{m}$. As $\overline{A}\in\mbox{GL}_{d}(k)$, $\mbox{det}(\overline{A})\neq0$
and so if we let $A$ be the matrix $(a_{ij})$
we have $\mbox{det}(A)\notin\mathfrak{m}$ and so $A\in\mbox{GL}{}_{d}(R)$ as $R$ is local. Now
fix $j\in\{1,\dots,d\}$. 

We have $\theta(\overline{v_{j}})=\sum_{i=1}^{d}a_{ij}v_{i}+C^{\flat}$ as $\mathfrak{m}C^{\sharp}\subseteq C^{\flat}$, and so by definition there is some $w_{j}\in Cv_{j}$ for which $\sum_{i=1}^{d}a_{ij}v_{i}-w_{j}\in C^{\flat}$. Let $z_{j}=w_{j}-\sum_{i=1}^{d}a_{ij}v_{i}$.
Write $z_{j}=z_{j}^{+}+z_{j}^{-}$ for elements $z_{j}^{+}\in(C^{-1})'\cap C''$
and $z_{j}^{-}\in(C^{-1})''\cap C'$. Hence there are some integers
$n_{-}$ and $n_{+}$, and a collection $\{z_{j,n}^{\pm}\mid n\in\mathbb{Z}\}\subseteq M$
for which: $z_{j,n}^{\pm}\in Cz_{j,n+1}^{\pm}$ for each $n\in\mathbb{Z}$;
$z_{j,n}^{-}=0$ for each $n>n_{-}$; and $z_{j,n}^{+}=0$ for each
$n<n_{+}$. Now for each $n\in\mathbb{Z}$ define the matrices $L^{\pm,n}$ by
\[
L^{+,n}=\begin{cases}
0 & (\mbox{if }n>0)\\
(A^{-1})^{1-n} & (\mbox{otherwise)}
\end{cases} \, \, \, L^{-,n}=\begin{cases}
-A^{n-1} & (\mbox{if }n>0)\\
0 & (\mbox{otheriwse})
\end{cases}
\]
Write $L^{\pm,n}=(m_{ij}^{\pm,n})_{i,j}$ for elements $m_{ij}^{\pm,n}\in R$. Note that (if $m_{ij}^{+,n}z_{i,n}^{+}\neq 0$ then $n_{+}\leq n\leq 0$) and (if $m_{ij}^{-,n}z_{i,n}^{-}\neq 0$ then $n_{-}\geq n\leq 1$). This means the sums $\sum_{n\in\mathbb{Z}}\sum_{i=1}^{d}m_{ij}^{\pm,n}z_{i,n}^{\pm}$ are finite. As in the proof of \cite[Lemma 4.6]{Cra2018} this gives $\sum_{j=1}^{d}a_{ij}u_{i}\in Cu_{j}$ where for each $i$ we let 
\[
\begin{array}{c}
u_{i}=v_{i}+\sum_{n\in\mathbb{Z}}(\sum_{k=1}^{d}m_{ki}^{+,n}z_{k,n}^{+}))+\sum_{n\in\mathbb{Z}}(\sum_{k=1}^{d}m_{ki}^{-,n}z_{k,n}^{-}).
\end{array}
\] 
Let $U=\bigoplus_{i=1}^{d}R$, $T(r_{i})_{i}=(\sum_{j=1}^{d}a_{ij}r_{j})_{i}$ and $f((r_{i}))=\sum_{t=1}^{d}r_{t}u_{t}$. Since $C^{\sharp}/C^{\flat}$ has $k$-basis $\overline{v_{1}},\dots,\overline{v_{d}}$,
for any $m\in C^{\sharp}$ there are elements $s_{1},\dots,s_{d}\in R$
such that writing $\overline{s_{i}}=s_{i}+\mathfrak{m}$ for each
$i$ gives $m+C^{\flat}=\sum_{i}\overline{s_{i}}(v_{i}+C^{\flat})$
which equals $\sum_{i}s_{i}u_{i}+C^{\flat}$. There is an element
$x=\sum_{i}s_{i}u_{j}=f((s_{j}))\in\im(f)$ and an element
$c\in C^{\flat}$ with $m=x+c$. This shows $C^{\sharp}\subseteq\im(f)+C^{\flat}$
and as $u_{i}\in C^{\sharp}$ for each $i$ we have equality. 

Since $\mathfrak{m}C^{\sharp}\subseteq C^{\flat}$
we have $\mathfrak{m}U\subseteq\{u\in U:f(u)\in C^{\flat}\}$. Conversely
if $f(u)\in C^{\flat}$ where $f(u)=\sum_{t=1}^{d}r_{t}u_{t}$ then
$\bar{0}=\sum_{t=1}^{d}\overline{r_{t}}\overline{u_{t}}=\sum_{t=1}^{d}\overline{r_{t}}\overline{v_{t}}$,
and as $\overline{v_{1}},\dots,\overline{v_{d}}$ is an $R/\mathfrak{m}=k$-basis
for $C^{\sharp}/C^{\flat}$, this means $r_{i}+\mathfrak{m}=0$ in
$k$ (and hence $r_{i}\in\mathfrak{m}$) for each $i$. Hence $\mathfrak{m}U\supseteq\{u\in U:f(u)\in C^{\flat}\}$.
Now fix $u=(r_{i})\in U$. By definition, $Tu=(\sum_{j=1}^{d}a_{ij}r_{j})_{i}$,
and so $f(Tu)=\sum_{j=1}^{d}r_{j}\sum_{t=1}^{d}a_{tj}u_{t}\in\sum_{j=1}^{d}r_{j}Cu_{j}\subseteq Cf(u)$.\end{proof}
In \S\ref{sec:Some-Linear-Relations}  we provide examples of $R$-linear relations.
\section{\label{sec:Some-Linear-Relations}Homotopy words and relations.\label{subsec:Linear-relations-of}}

\begin{assumption}\label{ass6}
In \S\ref{sec:Some-Linear-Relations} we fix a homotopically minimal complex of projective $\Lambda$-modules $M$. Hence $M^{i}$ is a projective $\Lambda$-module and $\im(d^{i}_{M})\subseteq\rad(M^{i+1})$ for each integer $i$.
\end{assumption}
\begin{notation}
We abuse notation writing $M$ for the projective $\Lambda$-module $\bigoplus_{i\in\mathbb{Z}}M^{i}$, and let $d_{M}$ be the $\Lambda$-module endomorphism $\bigoplus_{i\in\mathbb{Z}}d_{M}^{i}$
of $M$ sending $\sum_{i} m_{i}$ to $\sum_{i} d_{M}^{i}(m_{i})$. If $L$ is an $R$-module let $\mathrm{End}_{R}(L)$ be the ring of $R$-module endomorphisms of $L$. For each
vertex $v$ define $d_{v,M}\in\mathrm{End}_{R}(e_{v}M)$ by the restriction of $d_{M}$.
\end{notation} 
\begin{lemma}\label{lemma.3.9}\label{lemma.1.1}Fix $M$ as in Assumption \ref{ass6}.
\begin{enumerate} \item \emph{\cite[Lemma 2.1.1]{Ben2018} (}see also \emph{\cite[Lemma 5]{BekMer2003})}. If $a,b\in\mathbf{A}$:
\begin{enumerate}
\item if $v=h(b)=t(a)$ and $ab\in\mathbf{P}$, $abm=0$ implies
$bm=0$ for all $m\in M$;
\item if $v=t(a)$ then $\{m'\in e_{v}M\mid am'=0\}=\sum_{b' \in\mathbf{A}(\rightarrow v)\, :\, ab'\notin\mathbf{P}} b'M$; and
\item if $v=h(b)=h(a)$ and $a\neq b$ the sum $aM+bM$ is direct.
\end{enumerate}
\item \emph{\cite[Lemma 2.1.2]{Ben2018}} For each $\alpha\in\mathbf{A}$ there exists $d_{\alpha,M}\in \mathrm{End}_{R}(e_{h(\alpha)}M)$ such that $d_{v,M}=\sum_{\beta\in\mathbf{A}(\rightarrow v)}d_{\beta,M}$ for any vertex $v$. Furthermore for all $\tau\in\mathbf{P}$ and all $x\in e_{t(\tau)}M$:
\begin{enumerate}
\item if there exists $\sigma\in\mathbf{A}$ such that $\tau\sigma\in\mathbf{P}$ then $d_{ \mathrm{l}( \tau),M}(\tau x)=\tau d_{\sigma,M}(x)$; 
\item if $\tau\sigma\notin\mathbf{P}$ for all $\sigma\in\mathbf{A}$ then $d_{ \mathrm{l}( \tau),M}(\tau x)=0$; 
\item if $h(\theta)=h(\tau)$ for
some arrow $\theta\neq \mathrm{l}(\tau)$ then $d_{\theta,M}(\tau x)=0$;
\item if $h(\phi)=h(\tau)$ for some arrow $\phi$ then $d_{\phi,M}d_{ \mathrm{l}( \tau),M}=0$; and
\item if $\tau x\in\im(d_{ \mathrm{l}( \tau),M})$ then $d_{\varsigma,M}(x)=0$ for any arrow $\varsigma$ such that $\tau\varsigma\in\mathbf{P}$.
\end{enumerate}
\end{enumerate}
\end{lemma}
The proof of Lemma \ref{lemma.3.9}(i) is straightforward and omitted. The proof Lemma \ref{lemma.3.9}(ii) involves the definition of the $R$-module endomorphisms $d_{\alpha,M}:e_{h(\alpha)}M\to e_{h(\alpha}M$ which motivate the introduction of homotopy letters of the form $d^{\pm}_{\alpha}$.
\begin{proof}[of Lemma \ref{lemma.3.9}(ii)]
Since $M$ is projective we have $\mbox{rad}(M)=\rad(\Lambda)M$ (see for example \cite[Theorem 24.7]{Lam1991}), and by Assumption  \ref{ass6} $\mathrm{im}(d_{v,M})\subseteq e_{v}\mbox{rad}(\Lambda)M$. Note $e_{v}\mbox{rad}(\Lambda)=\bigoplus_{\beta}\beta\Lambda$ where $\beta$ runs through $\mathbf{A}(\rightarrow v)$, and so $e_{v}\mbox{rad}(M)=\sum_{\beta}\beta M$ which is a direct sum by part (ic). For any $\gamma\in\mathbf{A}(\rightarrow v)$ let $\pi_{\gamma}:\bigoplus_{\beta}\beta M\rightarrow\gamma M$
and $\iota_{\gamma}:\gamma M\rightarrow\bigoplus_{\beta}\beta M$
be the natural projections and inclusions of $R$-modules.
Define $d_{\alpha,M}$ by $d_{\alpha,M}(m)=\iota_{\alpha}(\pi_{\alpha}(  d_{v,M}(m)))$. The proof of parts (a), (b), (c), (d) and (e) of (ii) follow from part (i).
\end{proof}
\begin{definition}\cite[Example 1.4.2, Definitions 1.4.9 and 2.1.5]{Ben2018} If $p\in\mathbf{P}$ and $a\in\mathbf{A}$ let $\mathrm{rel}^{p}(M)=\{(m,pm)\mid m\in e_{t(p)}M\}$ and $\mathrm{rel}^{a}_{d}(M)=\{(m,d_{a,M}(m))\mid m\in e_{h(a)}M\}$. 

If $v$ is a vertex and $C= 1 _{v,\pm}$ let $\mathrm{rel}^{C}(M)$ be the relation $\{(m,m)\mid m\in e_{v}M\}$. Now let $C=l_{1}^{-1}r_{1}\dots l_{n}^{-1}r_{n}$ be a non-trivial homotopy word. For each $i$ with $0<i\leq n$ let
\[
\begin{array}{c}\mathrm{rel}^{C}_{\,i}(M)=\begin{cases}(\mathrm{rel}^{\gamma}(M))^{-1}\mathrm{rel}^{\mathrm{l}(\gamma)}_{d}(M) & (\text{if }l_{i}^{-1}r_{i}=\gamma^{-1}d_{\mathrm{l}(\gamma)})\\
(\mathrm{rel}^{\mathrm{l}(\gamma)}_{d}(M))^{-1}\mathrm{rel}^{\gamma}(M) & (\text{if }l_{i}^{-1}r_{i}=d_{\mathrm{l}(\gamma)}^{-1}\gamma)\end{cases}
\end{array}
\]
and let $\mathrm{rel}^{C}(M)=\mathrm{rel}^{C}_{\,1}(M)\dots\mathrm{rel}^{C}_{\,n}(M)$, the $n$-fold composition of these relations.

Let $q$ be a homotopy letter (that is, one of $\gamma$, $\gamma^{-1}$, $d_{\alpha}$ or $d_{\alpha}^{-1}$ for some path $\gamma\in\mathbf{P}$ or some arrow $\alpha$). If $U$ is a subset of $e_{t(q)}M$ then define the subset $qU$ of $e_{h(q)}M$ by
\[
\begin{array}{cc}
\gamma U=\{\gamma m\in e_{h(\gamma)}M\mid m\in U\}, & \gamma^{-1}U=\{m\in e_{t(\gamma)}M\mid\gamma m\in U\},\\
d_{\alpha}U=\{d_{\alpha,M}(m)\in e_{h(\alpha)}M\mid m\in U\}, & d_{\alpha}^{-1}U=\{m\in e_{h(\alpha)}M\mid d_{\alpha,M}(m)\in U\}.
\end{array}
\]
For any vertex $v$ and any subset $U$ of $e_{v}M$ let $ 1 _{v,\pm1}U=U$. When $U=e_{t(q)}M$ we let $qM=qU$. Similarly when $U=e_{t(q)}\mathrm{rad}(M)$ we let $q\mathrm{rad}(M)=qU$.
\end{definition}
\begin{example}\label{example.2.1.6}(See also \cite[Examples 2.1.6 and 2.1.8]{Ben2018}). Consider the finite-dimensional gentle algebra $\Lambda=kQ/\mathcal{J}$ from Example \ref{example8}. Let $C=s^{-1}d_{s}t{}^{-1}d_{t}d_{r}^{-1}rh$. Then
\[
\begin{array}{c}
\mathrm{rel}^{C}(M)=\{(m_{3},m_{0})\mid sm_{0} = d_{s,M}(m_{1}),\,tm_{1} = d_{t,M}(m_{2}) ,\,d_{h,M}(m_{2})=hm_{3}\text{ for some }m_{1},m_{2}\}.
\end{array}
\]
The elements $m_{0},\dots,m_{3}$ may be arranged using schemas from \cite[Remark 1.3.36]{Ben2018}, as follows
\[
\small{\xymatrix@C=0.6em@R=.3em{
 &  &  &  &  &  &  &  & m_{2}\ar[dr]_{r}\ar[dl]^{f}\ar@{.>}@/_{2pc}/[ddddlll]_{d_{t}}\ar@{.>}@/^{2pc}/[ddddrrrrr]^{d_{h}} &\\
 &  &  &  &  &  &  & fm_{2}\ar@{.>}@/^{1pc}/[ddddl]_(0.65){d_{f}}\ar[dl]^{w} &  & rm_{2}\ar@{.>}@/^{1.5pc}/[ddddrrr]_{d_{r}} &\\
 &  &  &  & & &  wfm_{2}\ar@{.>}@/^{1.5pc}/[ddddr]|>>>>>>>>>>>>>{\hole}^(0.7){d_{w}} &  &  &  &  &  &  & & &\\
 &  &  &  & m_{1}\ar[dr]_(0.35){t}\ar@{.>}@/_{1.5pc}/[dddddlll]_{d_{s}} &  &  &  &  &  &  &  & & & m_{3}\ar[dr]_(0.55){y}\ar[dl]^(0.35){h} &\\
 &  &  &  &  & tm_{1}\ar[dr]_(0.55){f} &  &  &  &  &  &  & & hm_{3}\ar[dl]^(0.55){r} & & ym_{3}\\
  &  &  &  &  &  & ftm_{1}\ar[dr]_{w} &  &  &  & & & rhm_{3} & \\
    &  &  &  &  &  & & wftm_{1} &  &  &  & &  & \\
m_{0}\ar[dr]_(0.35){s} &  &  &  &  &  &  &  &  &  & & \\
 & sm_{0}
}}
\]
Only the elements $m_{0}$, $m_{1}$, $m_{2}$ and $m_{3}$ are needed to describe the set $\mathrm{rel}^{C}(M)$. The bold arrows labelled by homotopy letters of the form $\alpha$ (for $\alpha\in\mathbf{A}$) indicate the action of $\Lambda$.  Their shape indicate the indecomposable summands $\Lambda e_{7}$, $\Lambda e_{2}$, $\Lambda e_{4}$ and $\Lambda e_{6}$ of $P(C)$. The dotted arrows indicate the action of $d_{M}$. Note $d_{M}(m_{2})=tm_{1}+hm_{3}$. For a subset $U\subseteq e_{7}M$, $CU$ is the set of $m_{0}\in e_{6}M$ above where $m_{3}\in U$.
\end{example}
\begin{corollary}\emph{\cite[Corollary 2.1.9]{Ben2018}}
\label{corollary.3.1-1}If $a$ is an arrow then \emph{$a^{-1}d_{a}\rad(M)\subseteq e_{t(a)}\rad(M)$}.
Furthermore, given an arrow $b$ with \emph{$ab\in\mathbf{P}$} we
have $(ab)^{-1}ad_{b}M=b{}^{-1}d_{b}M$.
\end{corollary}
\begin{proof}
From the definition of $d_{a,M}$ from Lemma \ref{lemma.1.1}(ii) we have $a^{-1}d_{a}\rad(M)=a^{-1}d_{a}aM$. By parts (a) and (b) of Lemma \ref{lemma.1.1}(ii) we have that if $b$
exists then $a^{-1}d_{a}aM=a^{-1}ad_{b}M$, and otherwise $a^{-1}d_{a}aM=a^{-1}0$. If $b$ exists then by Lemma \ref{lemma.3.9}(iia) any $m\in a^{-1}d_{a}aM$ satisfies $am=ad_{b,M}(m')$
for some $m'\in e_{h(b)}M$, and so $m-d_{b,M}(m')\in a^{-1}0=\sum b'M$ where the sum ranges over all arrows $b'$ with $ab'\notin\mathbf{P}$, by Lemma \ref{lemma.3.9}(ib). As
$\mbox{im}(d_{b,M})\subseteq\rad(M)$ this shows that if $b$ exists then $a^{-1}d_{a}aM\subseteq e_{t(\alpha)}\mbox{rad}(M)$. If $b$ does not exist then $a^{-1}d_{a}aM=a^{-1}d_{a}0=a^{-1}0=\sum b'M$, concluding that $a^{-1}d_{a}\rad(M)\subseteq e_{t(a)}\rad(M)$.

Now assume $b$ exists. If $m\in b{}^{-1}a{}^{-1}ad_{b}M$ there exists $m'\in e_{h(b)}M$ such that $bm-d_{b,M}(m')$ lies in $a^{-1}0$. As $bm-d_{b,M}(m')=bm''$
for some $m''\in M$, we have $abm''=0$ which means $bm''=0$ by Lemma \ref{lemma.3.9}(ia). This gives $b^{-1}a^{-1}ad_{b}M\subseteq b{}^{-1}d_{b}M$. The reverse inclusion is clear.
\end{proof}
\begin{corollary}\emph{\cite[Corollary 2.1.10]{Ben2018}}
\label{corollary.3.2-1}If $\alpha,\beta,\gamma,\sigma,\alpha\beta\in\mathbf{P}$,
$h(\gamma)=h(\sigma)$ and $ \mathrm{l}( \gamma)\neq  \mathrm{l}( \sigma)$ then we have the following inclusions
\[
\begin{array}{c}
\begin{array}{ccc}
\beta^{-1}d_{ \mathrm{l}( \beta)}M\subseteq(\alpha\beta)^{-1}d_{ \mathrm{l}( \alpha)}M, &  & d_{ \mathrm{l}( \alpha)}^{-1}\alpha\beta M\subseteq d_{ \mathrm{l}( \alpha)}^{-1}\alpha M,\end{array}\\
\begin{array}{ccccc}
\alpha^{-1}d_{ \mathrm{l}( \alpha)}M\subseteq d_{ \mathrm{l}( \beta)}^{-1}\beta0, &  & \gamma M\subseteq d_{ \mathrm{l}( \sigma)}^{-1}\sigma0, &  & d_{ \mathrm{l}( \sigma)}M\subseteq d_{ \mathrm{l}( \sigma)}^{-1}\sigma0.\end{array}
\end{array}
\]
\end{corollary}
In Definition \ref{definition.4.14} we define functors $C^{\pm}:\mathcal{C}_{\mathrm{min}}(\Lambda\text{-}\boldsymbol{\mathrm{Proj}})\rightarrow R\text{-}\boldsymbol{\mathrm{Mod}}$ (see Corollary \ref{corollary.3.4}). To do so we adapt ideas used by Ringel \cite{Rin1975} which were developed by Crawley-Boevey \cite{Cra2018}.
\begin{definition}
Choose a \textit{sign} $s(q)\in\{\pm1\}$ for each homotopy letter $q$ in  the set $\mathbf{A}{}^{\pm}$ of homotopy letters of the form $\alpha$ or $\alpha^{-1}$ for some $\alpha\in\mathbf{A}$, such that if distinct $q,q'\in\mathbf{A}{}^{\pm}$ satisfy $h(q)=h(q')$, then $s(q)=s(q')$ if and only if $\{q,q'\}=\{\alpha^{-1},\beta\}$
with $\alpha\beta\notin\mathbf{P}$. Now let $s(\gamma)=s( \mathrm{l}( \gamma))$, $s(\gamma^{-1})=s( \mathrm{f}( \gamma)^{-1})$,
and $s(d_{\alpha}^{\pm1})=-s(\alpha)$ for each $\gamma\in\mathbf{P}$
and each $\alpha$. 

For a (non-trivial finite or $\mathbb{N}$)-homotopy word
$C$ we let $h(C)$ and $s(C)$ be the head and sign of the first
letter of $C$. For any vertex $v$ let $s(1_{v,\pm1})=\pm1$
and $h(1_{v,\pm1})=v$. 

Let $D$ and $E$ be homotopy words where $I_{D^{-1}}\subseteq\mathbb{N}$ and $I_{E}\subseteq\mathbb{N}$.
If $u=h(D^{-1})$ and $\epsilon=-s(D^{-1})$ let $D 1 _{u,\epsilon}=D$.
If $v=h(E)$ and $\delta=s(E)$ we let $ 1 _{v,\delta}E=E$.
The \textit {composition} $DE$ is the concatenation of the homotopy letters in $D$ with those in $E$. By \cite[Proposition 2.1.13]{Ben2018}, $DE$ is a homotopy word
if and only if $h(D^{-1})=h(E)$ and $s(D^{-1})=-s(E)$. If $D=\dots l_{-1}^{-1}r_{-1}l_{0}^{-1}r_{0}$ is a $-\mathbb{N}$-word and $E=l_{1}^{-1}r_{1}l_{2}^{-1}r_{2}\dots$
is an $\mathbb{N}$-word, write $DE=\dots l_{0}^{-1}r_{0}\mid l_{1}^{-1}r_{1}\dots$. 
\end{definition}

\begin{corollary}
\label{corollary.3.5}\emph{\cite[Corollary 2.1.15]{Ben2018}} Let $a,b\in\mathbf{A}$ and $Ca^{-1}d_{a},Cd_{b}^{-1}b$ be homotopy words.
\begin{enumerate}
\item If $\gamma\in\mathbf{P}$ then $C\gamma^{-1}d_{ \mathrm{l}( \gamma)}$  is a homotopy word if and only if $\mathrm{f}( \gamma)=a$.
\item If $\tau\in\mathbf{P}$ then $Cd_{ \mathrm{l}( \tau)}^{-1}\tau$ is a homotopy word if and only if $\mathrm{l}( \tau)=b$.
\item If $\gamma'\in\mathbf{P}$ is longer than $\gamma\in\mathbf{P}$ and  $\mathrm{f}(\gamma')=\mathrm{f}( \gamma)=a$ then $C\gamma^{-1}d_{ \mathrm{l}( \gamma)}M\subseteq C\gamma'^{-1}d_{ \mathrm{l}( \gamma')}M$.
\item If $\tau'\in\mathbf{P}$ is longer than $\tau\in\mathbf{P}$ and   $\mathrm{l}( \tau')=\mathrm{l}( \tau)=b$ then $Cd_{ \mathrm{l}( \tau')}^{-1}\tau'M\subseteq Cd_{ \mathrm{l}( \tau)}^{-1}\tau M$.
\end{enumerate}
\end{corollary}
\begin{example}
\cite[Example 2.1.31]{Ben2018} For the
complete gentle algebra $k[[x,y]]/(xy)$  the iterated inclusions given by Corollaries \ref{corollary.3.2-1} and \ref{corollary.3.5} may be used to construct \textit{intervals} such as
\[
\xymatrix@R=.8em@C=.75em{M\ar@{=}[r]\ar@{-}[dddd] & d_{y}^{-1}yM\ar@{-}[d] &  &  & d_{y}^{-1}y^{2}M\ar@{=}[r]\ar@{-}[dddd] & d_{y}^{-1}y^{2}d_{x}^{-1}xM\ar@{-}[d]\\
 & d_{y}^{-1}y^{2}M\ar@{-}[d]\ar@{--}@(r,dl)[urrr] &  &  &  & d_{y}^{-1}y^{2}d_{x}^{-1}x^{2}M\ar@{-}[d]\\
 & d_{y}^{-1}y^{3}M\ar@{-}[d]\ar@{--}@(r,ul)[ddddddrrr] &  &  &  & d_{y}^{-1}y^{2}d_{x}^{-1}x^{3}M\ar@{-}[d]\\
 & \underset{\,}{\vdots}\ar@{-}[d] &  &  &  & \underset{\,}{\vdots}\ar@{-}[d]\\
( 1 _{1,1})^{+}(M)\ar@{=}[r] & \bigcap_{n>0}d_{y}^{-1}y^{n}M &  &  & (d_{y}^{-1}y^{2})^{+}(M)\ar@{=}[r] & \bigcap_{n>0}d_{y}^{-1}y^{2}d_{x}^{-1}x^{n}M\\
( 1 _{1,1})^{-}(M)\ar@{=}[r]\ar@{-}[ddd] & \bigcup_{n>0}y^{-n}d_{y}M\ar@{-}[d] &  &  & (d_{y}^{-1}y^{2})^{-}(M)\ar@{=}[r]\ar@{-}[ddd] & \bigcup_{n>0}d_{y}^{-1}y^{2}x^{-n}d_{x}M\ar@{-}[d]\\
 & \underset{\,}{\vdots}\ar@{-}[d] &  &  &  & \underset{\,}{\vdots}\ar@{-}[d]\\
 & y^{-1}d_{y}M\ar@{-}[d] &  &  &  & d_{y}^{-1}y^{2}x^{-1}d_{x}M\ar@{-}[d]\\
xM+d_{y}M\ar@{=}[r] & xM+d_{y}M &  &  & d_{y}^{-1}y^{3}M\ar@{=}[r] & d_{y}^{-1}y^{2}(yM+d_{x}M)
}
\]
We define the sets $( 1 _{1,1})^{\pm}(M)$ and $(d^{-1}_{y}y^{2})^{\pm}(M)$ in general in Definition \ref{definition.4.14}. 
\end{example}
\begin{notation}\cite[Definition 2.1.17]{Ben2018} If $C=\dots l_{i}^{-1}r_{i}\dots$
is a homotopy word and $i\in I_{C}$ is arbitrary, let $C_{i}=l_{i}^{-1}r_{i}$
and $C_{\leq i}=\dots l_{i}^{-1}r_{i}$ given $i-1\in I_{C}$, and
otherwise $C_{i}=C_{\leq i}= 1 _{h(C),s(C)}$. Similarly let $C_{>i}=l_{i+1}^{-1}r_{i+1}\dots$ given $i+1\in I_{C}$ and
otherwise $C_{>i}= 1 _{h(C^{-1}),-s(C^{-1})}$. So, $C_{<i}$ and $C_{\geq i}$ are the unique homotopy words with $C_{\leq i}=C_{<i}C_{i}$
and $C_{i}C_{>i}=C_{\geq i}$. For any vertex $v$ and $\delta=\pm1$,
let $\mathcal{W}_{v,\delta}$ be the set of homotopy
$I$-words with $I\subseteq\mathbb{N}$, head $v$ and sign $\delta$. 
\end{notation}
\begin{definition}
\label{definition.4.14}\cite[Definition 2.1.17]{Ben2018} Let $C\in\mathcal{W}_{v,\delta}$. Suppose $I_{C}$ is finite. If $a$ is an arrow and $Cd_{a}^{-1}a$ is
a homotopy word let $C^{+}(M)$ be the intersection $\bigcap_{\beta} Cd_{a}^{-1}\beta\rad(M)$ 
over $\beta\in\mathbf{P}$ with
$\mathrm{l}( \beta)=a$. By \cite[Lemma 2.1.19]{Ben2018}, if there are finitely many such $\beta$ then $C^{+}(M)= Cd_{a}^{-1}0$, and otherwise $C^{+}(M)=\bigcap_{\beta} Cd_{a}^{-1}\beta M$. If there is no arrow $a$ such that $Cd_{a}^{-1}a$ is
a homotopy word, then we let $C^{+}(M)=CM$.

If there exists an arrow $b$ where $Cb^{-1}d_{b}$ is a homotopy word let $C^{-}(M)$ be the union $\bigcup_{\alpha} C\alpha^{-1}d_{ b}M$  over all $\alpha\in\mathbf{P}$ with
$\mathrm{f}(\alpha)=b$. Otherwise let $C^{-}(M)=C(\sum d_{c(+)}M+\sum c(-)M)$
where $c(\pm)$ runs through all arrows with head $h(C^{-1})$ and
sign $\pm s(C^{-1})$.

Suppose instead $I_{C}=\mathbb{N}$. In this case let $C^{+}(M)$ be the set of all $m\in e_{v}M$ with
a sequence of elements $(m_{i})\in\prod_{i\in\mathbb{N}}e_{v_{C}(i)}M$
satisfying $m_{0}=m$ and $m_{i}\in l_{i+1}^{-1}r_{i+1}m_{i+1}$ for
each $i\geq0$, and let $C^{-}(M)$ be the subset of $C^{+}(M)$ where
each sequence $(m_{i})$ is eventually zero.
\end{definition}

\begin{corollary}
\label{corollary.3.4} Let $C\in\mathcal{W}_{v,\delta}$.
\begin{enumerate}
\item The assignments $M\mapsto C^{+}(M)$, $M\mapsto C^{-}(M)$ and ($M\mapsto CM$ for when $C$ is finite) respectively define subfunctors $C^{+}$, $C^{-}$ and ($C$ for when $C$ is finite) of the forgetful functor $\mathcal{C}_{\mathrm{min}}(\Lambda\text{-}\boldsymbol{\mathrm{Proj}})\rightarrow R\text{-}\boldsymbol{\mathrm{Mod}}$. Furthermore $C^{-}\leq C^{+}$.
\item If $I_{C}$ is finite then the functor $C$ preserves small direct sums and small direct products.
\item \emph{\cite[Corollary 2.1.20, Lemma 2.1.21]{Ben2018}} The functors $C^{\pm}$ preserve small direct sums.
\end{enumerate}
\end{corollary}
\begin{proof}
(i) We show that, if $C$ is finite, then $\im(\left.g\right|_{C(M)})\subseteq C(N)$ for any morphism $g:M\rightarrow N$ in $\mathcal{C}_{\mathrm{min}}(\Lambda\text{-}\boldsymbol{\mathrm{Proj}})$. Since $g$ is $\Lambda$-linear, for any $p\in\mathbf{P}$ we have $(g(m),g(m'))\in \mathrm{rel}^{p}(N)$ for all $(m,m')\in\mathrm{rel}^{p}(M)$. Next we show, for any $a\in\mathbf{A}$, that $(g(m),g(m'))\in \mathrm{rel}^{a}_{d}(N)$ assuming $(m,m')\in\mathrm{rel}^{a}_{d}(M)$. For $v=h(a)$ we have $gd_{v,M}=d_{v,N}g$ and so $\sum_{\beta}g(d_{\beta,M}(m))-d_{\beta,N}(g(m))=0$ where $\beta$ runs through $\mathbf{A}(\rightarrow v)$, which is a direct sum by Lemma \ref{lemma.3.9}(ic).

This shows that, if $C$ is finite, $C$ is functorial. From this, and by Definition \ref{definition.4.14}, we have that $\im(\left.g\right|_{C^{\pm}(M)})\subseteq C^{\pm}(N)$ for any $C$ and for any morphism $g:M\rightarrow N$ in $\mathcal{C}_{\mathrm{min}}(\Lambda\text{-}\boldsymbol{\mathrm{Proj}})$. This shows $C^{\pm}$ is functorial. The proof that $C^{-}(M)\subseteq C^{+}(M)$ follows by applying Corollaries \ref{corollary.3.2-1} and \ref{corollary.3.5} to Definition \ref{definition.4.14} (see the proof of \cite[Corollary 2.1.20]{Ben2018}(ii) for details). 

(ii) Fix a set $S$ and a collection $(M(s)\mid s\in S)$ of objects in $\mathcal{C}_{\mathrm{min}}(\Lambda\text{-}\boldsymbol{\mathrm{Proj}})$. Recall the direct product $\prod M(s)$ (respectively direct sum $\bigoplus M(s)$)  over this collection is the complex whose homogeneous part in degree $n\in\mathbb{Z}$ is $\prod _{s}M(s)^{n}$ (respectively $\bigoplus _{s}M(s)^{n}$), and whose differential is defined by $d_{M}((m_{s})_{s})=(d_{M(s)}(m_{s}))_{s}$. In general, considering $M(s)$ as a $\Lambda$-module for each $s$, we have $\bigoplus\mathrm{rad}(M(s)) = \mathrm{rad}(\bigoplus M(s))$. Similarly note $\mathrm{rad}(\prod M(s))\subseteq \prod \mathrm{rad}(M(s))$, which we now prove is an equality. Since $\Lambda$ is semilocal we have $\mathrm{rad}(X)=\mathrm{rad}(\Lambda)X$ for any $\Lambda$-module $X$ (see for example \cite[Proposition 24.4]{Lam1991}), and so $\mathrm{rad}(\prod M(s))=\mathrm{rad}(\Lambda)\prod M(s)$ and $\mathrm{rad}(M(s))=\mathrm{rad}(\Lambda)M(s)$ for each $s$. Hence one can show $\prod \mathrm{rad}(M(s))\subseteq\mathrm{rad}(\Lambda)\prod M(s)$, since $\mathrm{rad}(\Lambda)$ is finitely generated by $\mathbf{A}$.

Fix some $\gamma\in\mathbf{P}$ with $a=\mathrm{l}(\gamma)$ and $u=t(\gamma)$. Now let $N(s)$ be an $R$-submodule of $e_{v}M(s)$ for each $s$. Since $e_{v}\mathrm{rad}(\prod M(s))= \prod e_{v}\mathrm{rad}(M(s))$ we have $d_{a,\prod M(s)}((m_{s}))=(d_{a,M(s)}(m_{s}))$ for any $(m_{s})\in e_{v}\prod M(s)=\prod e_{v}M(s)$. Hence we have that $d_{a}^{-1}\gamma(\prod N(s))=\prod d_{a}^{-1}\gamma N(s)$. Similarly since $e_{v}\mathrm{rad}(\bigoplus M(s))=\bigoplus e_{v}\mathrm{rad}(M(s))$ we have $d_{a}^{-1}\gamma(\bigoplus N(s))=\bigoplus d_{a}^{-1}\gamma N(s)$. It follows, using similar arguments, that $\gamma^{-1}d_{a}(\bigoplus N(s))=\bigoplus \gamma^{-1}d_{a}N(s)$ and $\gamma^{-1}d_{a}(\prod N(s))=\prod \gamma^{-1}d_{a}N(s)$. Altogether we have $\l^{-1}r(\bigoplus N(s))=\bigoplus l^{-1}rN(s)$ and $l^{-1}r(\prod N(s))=\prod l^{-1}rN(s)$ whenever $l^{-1}r$ is a homotopy $\{0,1\}$-word. Let $C=l_{1}^{-1}r_{1}\dots l_{n}^{-1}r_{n}$. Then by iterating the above equalities (in case $n>0$) we have
\[
C(\bigoplus M(s))=l_{1}^{-1}r_{1}\dots l_{n-1}^{-1}r_{n-1}(\bigoplus l_{n}^{-1}r_{n}M(s))=\dots=\bigoplus CM(s).
\]
Similarly one can show $C(\prod M(s))=\prod CM(s)$.

(iii) Suppose firstly that $I_{C}$ is finite. We start by showing $C^{-}(\bigoplus M(s))=\bigoplus C^{-}M(s)$. Let $Cb^{-1}d_{b}$ be a homotopy word for some $b\in\mathbf{A}$. By part (ii) we have 
\[\begin{array}{c}
C^{-}(\bigoplus M(s))=\bigcup_{\alpha} C\alpha^{-1}d_{ b}(\bigoplus_{s} M(s))=\bigcup_{\alpha} \bigoplus_{s} C\alpha^{-1}d_{ b}M(s)\subseteq \bigoplus C^{-}M(s) 
\end{array}
\]
where the union runs over all $\alpha\in\mathbf{P}$ with $\mathrm{f}(\alpha)=b$. Suppose conversely that $(m_{s})\in \bigoplus C^{-}(M(s))$, and so for each $s$ there is some $\alpha(s)\in\mathbf{P}$ with $\mathrm{f}(\alpha(s))=b$ such that $m_{s}\in C\alpha(s)^{-1}d_{a}M(s)$. Since $m_{t}=0$ (and so $\alpha(t)=b$) for all but finitely many $t$, we can choose $\gamma$ of maximal length among $(\alpha(s))$. This gives $(m_{s})\in C\gamma^{-1}d_{ b}(\bigoplus_{s} M(s))\subseteq C^{-}(\bigoplus M(s))$. 

Assume instead $Cb^{-1}d_{b}$ is not a homotopy word for all $b\in\mathbf{A}$. By definition $C^{-}(\bigoplus M(s))=C(\sum d_{c(+)}\bigoplus M(s)+\sum c(+)\bigoplus M(s))$
where $c(\pm)$ runs through all arrows with head $h(C^{-1})$ and
sign $\pm s(C^{-1})$. As above, and by part (ii), we have $\sum d_{c(+)}\bigoplus M(s)=\bigoplus(\sum d_{c(+)} M(s))$ and $\sum c(-)\bigoplus M(s)=\bigoplus(\sum c(-)M(s))$. Altogether this shows  $C^{-}(\bigoplus M(s))=\bigoplus C^{-}M(s)$.

We now show $C^{+}(\bigoplus M(s))=\bigoplus C^{+}M(s)$. If $Cd^{-1}_{a}a$ is not a homotopy word for all $a\in\mathbf{A}$ then $C(\bigoplus M(s))=\bigoplus CM(s)$ and so $C^{+}(\bigoplus M(s))=\bigoplus C^{+}M(s)$. Suppose instead $Cd^{-1}_{a}a$ is a homotopy word for some $a\in\mathbf{A}$. Then by part (ii) we have
\[\begin{array}{c}
C^{+}(\bigoplus M(s))=\bigcap_{\beta} Cd^{-1}_{a}\beta(\bigoplus_{s} M(s))=\bigcap_{\beta} \bigoplus_{s}  Cd^{-1}_{a}\beta M(s)= \bigoplus C^{+}M(s) 
\end{array}
\]
where the intersection runs over all $\beta\in\mathbf{P}$ with $\mathrm{f}(\beta)=a$.

Suppose finally that $C$ is a homotopy $\mathbb{N}$-word. By definition, $(m_{s})\in C^{+}(\bigoplus M(s))$ if and only if (there is a sequence $(m_{0,s}),(m_{1,s}),\dots \in\bigoplus M(s)$ with $m_{s}=m_{0,s}$ for all $s$ and $(m_{n-1,s})\in l^{-1}_{n}r_{n}(m_{n,s})$ for all $n>0$). As in the proof of (ii) we have $d_{a,\bigoplus M(s)}((m_{s}))=(d_{a,M(s)}(m_{s}))$ for any $a\in\mathbf{A}$, and so ($(m_{n-1,s})\in l^{-1}_{n}r_{n}(m_{n,s})$ if and only if $m_{n-1,s}\in l^{-1}_{n}r_{n}m_{n,s}$ for all $s$). Altogether this shows that $C^{+}(\bigoplus M(s))=\bigoplus C^{+}M(s)$. 

For each $s$ let $m_{s}\in C^{-}(M(s))$. So there is a sequence $m_{0,s},m_{1,s},m_{2,s},\dots \in M(s)$ such that $m_{s}=m_{0,s}$, $m_{n-1,s}\in l^{-1}_{n}r_{n}m_{n,s}$ for all $n>0$ and  $m_{n,s}=0$ for all $n\geq l(s)$ for some $l(s)\in\mathbb{N}$. We can assume $l(s)=0$ for all $s$ with $m_{s}=0$, which holds for all but finitely many $s$. Hence there exists $r>0$ with $r>l(s)$ for all $s$. So we have $(m_{0,s}),(m_{1,s}),\dots \in\bigoplus M(s)$ such that $(m_{s})=(m_{0,s})$, $(m_{n-1,s})\in l^{-1}_{n}r_{n}(m_{n,s})$ for all $n>0$ and $(m_{n,s})=0$ for all $n\geq r$. This shows $\bigoplus C^{-}M(s)\subseteq C^{-}(\bigoplus M(s))$. The reverse inclusion is straightforward.
\end{proof}
\section{\label{Refined Functors for Complexes.}Refined functors.}

\begin{definition}
\label{definition.3.3}Fix a vertex $v$ and $\delta=\pm1$. Below we define a total order $<$ on $\mathcal{W}_{v,\delta}$. \begin{enumerate}
\item \cite[Lemma 2.1.22, Definition 2.1.23]{Ben2018} For distinct $l^{-1}r,l'^{-1}r'\in \mathcal{W}_{v,\delta}$ let $l^{-1}r<l'^{-1}r'$ if:
\begin{enumerate}
\item $l^{-1}r=d_{ \mathrm{l}( \gamma)}^{-1}\gamma\nu$ and $l'^{-1}r'=d_{ \mathrm{l}( \gamma)}^{-1}\gamma$ for some $\gamma,\nu\in\mathbf{P}$ such that $\gamma\nu\in\mathbf{P}$; or
\item $l^{-1}r=\mu^{-1}d_{ \mathrm{l}( \mu)}$ and $l'^{-1}r'=d_{ \mathrm{l}( \eta)}^{-1}\eta$ for some $\mu,\eta\in\mathbf{P}$ such that $\mu\eta\in\mathbf{P}$; or
\item $l^{-1}r=\lambda^{-1}d_{ \mathrm{l}( \lambda)}$ and $l'^{-1}r'=\lambda^{-1}\kappa^{-1}d_{ \mathrm{l}( \kappa)}$ for some $\kappa,\lambda\in\mathbf{P}$ such that $\kappa\lambda\in\mathbf{P}$.
\end{enumerate}
\item \cite[Definition 2.1.26, Lemma 2.1.27]{Ben2018} For distinct $C,C'\in\mathcal{W}_{v,\delta}$
let $C<C'$ if:
\begin{enumerate}
\item there are homotopy letters $l$, $l'$, $r$ and $r'$ and homotopy
words $B,D,D'$ for which $C=Bl^{-1}rD$, $C'=Bl'^{-1}r'D'$ and $l^{-1}r<l'^{-1}r'$; or
\item there is some $\beta\in\mathbf{P}$ for which $C'=Cd_{ \mathrm{l}( \beta)}^{-1}\beta E$
for some homotopy word $E$; or
\item there is some $\alpha\in\mathbf{P}$ for which $C=C'\alpha^{-1}d_{ \mathrm{l}( \alpha)}E'$
for some homotopy word $E'$.
\end{enumerate}
\end{enumerate}
\end{definition}
\begin{lemma}\label{lemma.3.6}\label{lemma.3.4} Fix a vertex $v$ and $\delta=\pm1$.
\begin{enumerate}
\item \emph{\cite[Lemma 2.1.29]{Ben2018}} If \emph{(}$\alpha\in\mathbf{P}$, $s(\alpha)=\delta$, $v=h(\alpha)$ and $C\in\mathcal{W}_{v,\delta}$\emph{)} then $ \mathrm{l}( \alpha)M\subseteq C^{-}(M)$.
\item  \emph{\cite[Lemma 2.1.28]{Ben2018}} If $l^{-1}rD,l'^{-1}r'D'\in\mathcal{W}_{v,\delta}$
for homotopy words $D,D'$ and homotopy letters $l,l',r,r'$ where $l^{-1}r<l'^{-1}r'$, then $(l^{-1}rD)^{+}(M)\subseteq(l'^{-1}r'D')^{-}(M)$.
\item \emph{\cite[Proposition 2.1.30]{Ben2018}} For all $C,C'\in\mathcal{W}_{v,\delta}$
with $C<C'$ we have $C^{+}(M)\subseteq C'{}^{-}(M)$.
\end{enumerate}
\end{lemma}
\begin{proof}
(i) Suppose firstly
that $C$ is trivial.
If $C\beta^{-1}d_{ \mathrm{l}( \beta)}$
is a homotopy word for some $\beta\in\mathbf{P}$ then $h(\beta^{-1})=h(C^{-1})=h( 1 _{h(\alpha),s(\alpha)})$
and so $t(\beta)=h(\alpha)$. Similarly $s( \mathrm{f}( \beta)^{-1})=s( \mathrm{l}( \alpha))$ which gives $ \mathrm{f}( \beta) \mathrm{l}( \alpha)\notin\mathbf{P}$, so $ \mathrm{l}( \alpha)M\subseteq\beta^{-1}0$ which lies in $C^{-}(M)$. Otherwise there is no such $\beta$, and so $ \mathrm{l}( \alpha)M\subseteq C^{-}(M)$ by Definition \ref{definition.4.14}.

Suppose $C$ is non-trivial. If $C=\beta^{-1}d_{ \mathrm{l}( \beta)}D$ for
some homotopy word $D$ and some $\beta\in\mathbf{P}$ then $t(\beta)=h(C)=h(\alpha)$
and $s( \mathrm{f}( \beta)^{-1})=s(\alpha)=s( \mathrm{l}( \alpha)).$ So as before $ \mathrm{f}( \beta) \mathrm{l}(\alpha)\notin\mathbf{P}$
and again $ \mathrm{l}( \alpha)M\subseteq\beta^{-1}0\subseteq C^{-}(M)$. Suppose $C=d_{ \mathrm{l}( \gamma)}^{-1}\gamma E$, for some
homotopy word $E$ and $\gamma\in\mathbf{P}$. Here $h(\gamma)=h(\alpha)$
and $s( \mathrm{l}( \alpha))=-s( \mathrm{l}( \gamma)).$ So $ \mathrm{l}( \alpha)\neq  \mathrm{l}( \gamma)$,
and so $ \mathrm{l}( \alpha)M\subseteq d_{ \mathrm{l}( \gamma)}^{-1}0$ by Lemma \ref{lemma.1.1}(iic), which means $ \mathrm{l}( \alpha)M\subseteq d_{ \mathrm{l}( \gamma)}^{-1}\gamma E{}^{-}(M)=C^{-}(M)$. This shows that, in any case, $ \mathrm{l}( \alpha)M\subseteq C^{-}(M)$.

(ii), (iii) Suppose $l^{-1}r=d_{ \mathrm{l}( \gamma)}^{-1}\gamma\nu$ and $l'^{-1}r'=d_{ \mathrm{l}( \gamma)}^{-1}\gamma$ for some $\gamma,\nu\in\mathbf{P}$ such that $\gamma\nu\in\mathbf{P}$. Here $(l^{-1}rD)^{+}(M)\subseteq d_{\mathrm{l}(\gamma)}^{-1}\gamma \mathrm{l}(\nu)M$ and $d_{\mathrm{l}(\gamma)}^{-1}\gamma(D')^{-}(M)=(l'^{-1}r'D')^{-}(M)$. Since $\gamma\nu\in\mathbf{P}$ we have that $\mathrm{f}(\gamma)\mathrm{l}(\nu)\in\mathbf{P}$ and so $s(\mathrm{l}(\nu))=-s((\mathrm{f}(\gamma)^{-1})$. Since $d_{\mathrm{l}(\gamma)}^{-1}\gamma D'$ is a homotopy word we have that $s((d_{\mathrm{l}(\gamma)}^{-1}\gamma)^{-1})=-s(D')$, and so $\mathrm{l}(\nu)=s(D')$. By part (i) this means $\mathrm{l}(\nu)M\subseteq (D')^{-}(M)$. For the proof of (ii), the other cases of $l^{-1}r<l'^{-1}r'$ are similar and omitted. The proof of part (iii) follows from part (ii) by considering cases (iia), (iib) and (iic) of Definition \ref{definition.3.3}.
\end{proof}
\begin{assumption}
In the remainder of \S\ref{Refined Functors for Complexes.} let $M$ and $N$ be objects in
$\mathcal{C}_{\mathrm{min}}(\Lambda\text{-}\boldsymbol{\mathrm{Proj}})$, and let $B$ and $D$ be homotopy words with head $v$ such that $C=B^{-1}D$ is a homotopy word. 
\end{assumption}
\begin{corollary}
\label{corollary.3.1}\emph{\cite[Corollary 2.2.3]{Ben2018}} Let $n\in\mathbb{Z}$, $I_{C}=\{0,\dots,t\}$ and $X^{i}$ and $Y^{i}$
be $R$-submodules of $e_{t(C)}M^{i}$ and $e_{h(C)}M^{i}$ for all $i\in\mathbb{Z}$. Then $Y^{n}\cap C(\bigoplus_{i}X^{i})=Y^{n}\cap CX^{n+\mu_{C}(t)}$. 
\end{corollary}
\begin{proof} Both the map $d_{a,M}$ (for $a\in\mathbf{A}$), and the map given by multiplication by $\gamma\in\mathbf{P}$, are homogeneous: and so $Y^{n}\cap C(\bigoplus_{i}X^{i})\subseteq Y^{n}\cap CX^{n+\mu_{C}(t)}$. The reverse inclusion is clear.
\end{proof}
\begin{lemma}
\label{lemma.6.1}\emph{\cite[Lemma 2.2.4]{Ben2018}} We have the following inclusions:
\begin{enumerate}
\item $B^{+}(M)\cap D^{+}(M)\cap e_{v}\rad(M)\subseteq(B^{+}(M)\cap D^{-}(M))+(B^{-}(M)\cap D^{+}(M))$;
\item $(B^{-}(M)+D^{+}(M)\cap B^{+}(M))\cap e_{v}\rad(M)\subseteq(B^{-}(M)+D^{-}(M)\cap B^{+}(M))$;
and
\item $B^{+}(M)\cap D^{\pm}(M)+e_{v}\rad(M)=(B^{+}(M)+e_{v}\rad(M))\cap(D^{\pm}(M)+e_{v}\rad(M))$.
\end{enumerate}
\end{lemma}
\begin{proof}
Without loss of generality we assume that $s(B)=1$, $s(D)=-1$ and that $\mathbf{A}(\rightarrow v)=\{x(+),x(-)\}$ where $s(x(\pm))=\pm 1$.

(i) Since $\Lambda$ is semilocal we have $\mbox{rad}(M)=\mbox{rad}(\Lambda)M$. So for any $m\in e_{v}\mbox{rad}(M)$ there are some $m_{+},m_{-}\in M$
for which $m=x(-)m_{-}+x(+)m_{+}$. By Lemma \ref{lemma.3.4}(i) we have that $x(+)m_{+}\in B^{-}(M)$ and
$x(-)m_{-}\in D^{-}(M)$. If also $m\in B^{+}(M)\cap D^{+}(M)$
we have $x(+)m_{+}\in D^{+}(M)\cap B^{-}(M)$ as $x(+)m_{+}=m-x(-)m_{-}$
and $D^{-}(M)\subseteq D^{+}(M)$. By symmetry $x(-)m_{-}\in B^{+}(M)\cap D^{-}(M)$. 

(ii) If $m\in(B^{-}(M)+D^{+}(M)\cap B^{+}(M))\cap e_{v}\rad(M)$
then $m=x(+)m_{+}+x(-)m_{-}$ for some $m_{+},m_{-}\in M$
as above. We also have $m=m'+m''$ where $m'\in B^{-}(M)$
and $m''\in D^{+}(M)\cap B^{+}(M)$, and so $x(-)m_{-}=m'+m''-x(+)m_{+}\in B^{+}(M)$.

(iii) Clearly $B^{+}(M)\cap D^{\pm}(M)+e_{v}\rad(M)$ is contained
in the intersection of $B^{+}(M)+e_{v}\rad(M)$ and $D^{\pm}(M)+e_{v}\rad(M)$.
Any $m$ from this intersection may be written as:
$m_{B}+x(-)m_{-}+x(+)m_{+}$ for $m_{B}\in B^{+}(M)$ and $m_{\pm1}\in M$;
and as $m_{D}+x(-)m'_{-}+x(+)m'_{+}$ for $m_{D}\in D^{\pm}(M)$
and $m'_{\pm1}\in M$. So, $x(+)m_{+},x(+)m'_{+}\in B^{-}(M)$
and $x(-)m_{-},x(-)m'_{-}\in D^{-}(M)\subseteq D^{\pm}(M)$.
This means $m_{B}+x(+)m_{+}-x(+)m'_{+}=x(-)m'_{-}-x(-)m_{-}+m_{D}$ lies in $D^{\pm}(M)$. This shows $m=(m_{B}+x(+)m_{+}-x(+)m'_{+})+(x(+)m'_{+}+x(-)m_{-})$ lies in $B^{+}(M)\cap D^{\pm}(M)+e_{v}\rad(M)$.
\end{proof}
\begin{definition}\cite[Definition 2.2.2]{Ben2018} For any $n\in\mathbb{Z}$ consider the $R$-submodules of $e_{v}M^{n}$
\[
\begin{array}{c}
F_{B,D,n}^{+}(M)=M^{n}\cap\left(B^{+}(M)\cap D^{+}(M)\right),\\
F_{B,D,n}^{-}(M)=M^{n}\cap\left(B^{+}(M)\cap D^{-}(M)+B^{-}(M)\cap D^{+}(M)\right),\\
G_{B,D,n}^{\pm}(M)=M^{n}\cap\left(B^{-}(M)+D^{\pm}(M)\cap B^{+}(M)\right).
\end{array}
\]
Define the quotients $F_{B,D,n}(M)$ and $G_{B,D,n}(M)$ by
\[
\begin{array}{cc}
F_{B,D,n}(M)=F_{B,D,n}^{+}(M)/F_{B,D,n}^{-}(M), &  G_{B,D,n}(M)=G_{B,D,n}^{+}(M)/G_{B,D,n}^{-}(M).
\end{array}
\]
\cite[Definition 2.2.5, Lemma 2.2.7]{Ben2018} Let 
\[
\begin{array}{cc}
\bar{F}_{B,D,n}(M)=\bar{F}_{B,D,n}^{+}(M)/\bar{F}_{B,D,n}^{-}(M), & 
\bar{G}_{B,D,n}(M)=\bar{G}_{B,D,n}^{+}(M)/\bar{G}_{B,D,n}^{-}(M)
\end{array}
\]
where
\[
\begin{array}{cc}
\bar{F}_{B,D,n}^{\pm}(M)=F_{B,D,n}^{\pm}(M)+e_{v}\rad(M^{n}), & \bar{G}_{B,D,n}^{\pm}(M)=G_{B,D,n}^{\pm}(M)+e_{v}\rad(M^{n}).
\end{array}
\]
Let $ \bar{A}^{\pm}(M)=A^{\pm}(M)+e_{h(C)}\rad(M)$ for any  homotopy word $A$ with $I_{A}\subseteq\mathbb{N}$. 
\end{definition}
\begin{remark}\label{remreffun} By Lemma \ref{lemma.6.1} we have
\[
\begin{array}{c}
\bar{F}_{B,D,n}^{+}(M)=e_{v}M^{n}\cap\bar{B}^{+}(M)\cap\bar{D}^{+}(M)\text{,}\\
\bar{F}_{B,D,n}^{-}(M)=e_{v}M^{n}\cap((\bar{B}^{+}(M)\cap\bar{D}^{-}(M))+(\bar{B}^{-}(M)\cap\bar{D}^{+}(M)))\text{, and}\\
\text{and }\bar{G}_{B,D,n}^{\pm}(M)=e_{v}M^{n}\cap(\bar{B}^{-}(M)+\bar{D}^{\pm}(M)\cap\bar{B}^{+}(M)).
\end{array}
\]
Consequently $F_{B,D,n}$, $\bar{F}_{B,D,n}$, $G_{B,D,n}$,
and $\bar{G}_{B,D,n}$ all define naturally isomorphic additive functors \cite[Lemma 2.2.7]{Ben2018}. Furthermore note that $\mathrm{im}(n)\subseteq \mathrm{rad}(N)$ for any null-homotopic morphism $n:M\to N$ between homotopically minimal complexes of projectives, and so $F_{B,D,n}$, $\bar{F}_{B,D,n}$, $G_{B,D,n}$,
and $\bar{G}_{B,D,n}$ all define functors $\mathcal{K}_{\mathrm{min}}(\Lambda\text{-}\boldsymbol{\mathrm{Proj}})\rightarrow k\text{-}\boldsymbol{\mathrm{Mod}}$ \cite[Corollary 2.2.8]{Ben2018}. Moreover, if $M$ is a complex of finitely generated
projectives, then since $\Lambda$ is $R$-module finite, $F_{B,D,n}(M)$ and $G_{B,D,n}(M)$ are finite-dimensional \cite[Corollary 2.2.9]{Ben2018}.

\cite[Definition 2.2.10]{Ben2018} Assume $C$ is $p$-periodic, so there is a homotopy word $E=l_{1}^{-1}r_{1}\dots l_{p}^{-1}r_{p}$ with $C={}^{\infty}E{}^{\,\infty}$. For each $n\in\mathbb{Z}$ let $E(n)$ be the relation on $e_{v}M^{n}$ given by all $(m,m')\in e_{v}M^{n}\oplus e_{v}M^{n}$ with $m\in Em'$. Then: $E(n)^{\sharp}=F_{B,D,n}^{+}(M)$
and $E(n)^{\flat}=F_{B,D,n}^{-}(M)$; and by Lemma \ref{lemma.2.10} there is a $k$-vector space automorphism of $E(n)^{\sharp}/E(n)^{\flat}$ defined by sending $m+E(n)^{\flat}$ to $m'+E(n)^{\flat}$
if and only if $m'\in E(n)^{\sharp}\cap(E(n)^{\flat}+E(n)m)$ (see \cite[Lemma 2.2.11]{Ben2018}). 

Consequently, as above, if $C$ is periodic then $F_{B,D,n}$,
$\bar{F}_{B,D,n}$, $G_{B,D,n}$, and $\bar{G}_{B,D,n}$ all define
naturally isomorphic functors $\mathcal{K}_{\mathrm{min}}(\Lambda\text{-}\boldsymbol{\mathrm{Proj}})\rightarrow k[T,T^{-1}]\text{-}\boldsymbol{\mathrm{Mod}}$. Furthermore these functors take objects in $\mathcal{K}_{\mathrm{min}}(\Lambda\text{-}\boldsymbol{\mathrm{proj}})$ to objects in $k[T,T^{-1}]\text{-}\boldsymbol{\mathrm{Mod}}_{k\text{-}\boldsymbol{\mathrm{mod}}}$, the full subcategory of $k[T,T^{-1}]\text{-}\boldsymbol{\mathrm{Mod}}$
consisting of finite-dimensional modules \cite[Corollary 2.2.12]{Ben2018}.
\end{remark}
\section{\label{subsec:Natural-Isomorphisms.}Natural isomorphisms and constructive functors..}
\begin{definition}\label{defintion.7.12}\cite[Definition 2.2.16]{Ben2018} Let $\Sigma$ be the set of all triples $(B,D,n)$ where $B^{-1}D$ is a homotopy word (equivalently $(B,D)\in\mathcal{W}_{v,\pm 1}\times\mathcal{W}_{v,\mp 1}$) and $n$ is an integer. 
\end{definition}
\begin{assumption}\label{ass9}
In \S\ref{subsec:Natural-Isomorphisms.} fix $(B,D,n),(B',D',n')\in\Sigma$ and let $C=B^{-1}D$ and $C'=B'^{-1}D'$.
\end{assumption}
\begin{definition}\label{definition.7.14}\cite[Definition 2.2.14]{Ben2018}  We write $C\sim C'$ 
if and only if $C'=C^{\pm1}[t]$ for some $t\in\mathbb{Z}$. So either ($C'=C^{\pm1}$ and $I_{C}\neq\mathbb{Z}\neq I_{C'}$) or ($C'=C^{\pm1}[t]$ and $I_{C}=I_{C'}=\mathbb{Z}$) \cite[Lemma 2.2.17]{Ben2018} (see also \cite[Lemma 2.1]{Cra2018}). Define the \textit{axis} $a_{B,D}\in\mathbb{Z}$ of $(B,D)$ by $C_{\leq a_{B,D}}=B^{-1}$ and $C_{>a_{B,D}}=D$. If $I_{C}=\{0,\dots,m\}$ then $a_{D,B}=m-a_{B,D}$; if $I_{C}=\pm\mathbb{N}$ then $a_{D,B}=-a_{B,D}$; and if $I_{C}=\mathbb{Z}$ then $a_{B,D}=0$ \cite[Lemma 2.2.15]{Ben2018}. 

\cite[Definition 2.2.18]{Ben2018} If $C\sim C'$ let
\[
r(B,D;B',D')=\begin{cases}
\mu_{C}(a_{B',D'})-\mu_{C}(a_{B,D}) & \mbox{(if }C'=C\mbox{ is not a homotopy }\mathbb{Z}\mbox{-word)}\\
\mu_{C}(a_{D',B'})-\mu_{C}(a_{B,D}) & \mbox{(if }C'=C^{-1}\mbox{ is not a homotopy }\mathbb{Z}\mbox{-word)}\\
\mu_{C}(\pm t) & \mbox{(if }C'=C^{\pm1}[t]\mbox{ is a homotopy }\mathbb{Z}\mbox{-word)}
\end{cases}
\]
We write $(B,D,n)\sim(B',D',n')$ if and only if
\[B^{-1}D\text{ and }B'^{-1}D'\text{ are equivalent and }n'-n=r(B,D;B',D').
\]
By \cite[Lemma 2.2.19]{Ben2018} we have that $r(B,D;B',D')=-r(B',D';B,D)$ and $r(B,D;B'',D'')=r(B,D;B',D')+r(B',D';B'',D'')$ for all $(B'',D'',n'')\in\Sigma$ with $B''^{-1}D''\sim B'^{-1}D'$, and so $\sim$ is an equivalence relation. 
\end{definition}
Recall the involution $\mathrm{res}_{\iota,R}$ of $R[T,T^{-1}]\text{-}\boldsymbol{\mathrm{Mod}}$, which swaps the action of $T$ and $T^{-1}$, and which restricts to an involution of $R[T,T^{-1}]\text{-}\boldsymbol{\mathrm{Mod}}_{R\text{-}\boldsymbol{\mathrm{Proj}}}$.
\begin{lemma}\label{lemma.5.7}\label{corollary.6.6} Fix $(B,D,n)$, $(B',D',n')$, $C$ and $C'$ as in Assumption \ref{ass9}.
\begin{enumerate}
\item \emph{\cite[Lemma 2.2.21]{Ben2018}} \emph{(}see also \emph{\cite[Lemma 7.1]{Cra2018})}. 

\begin{enumerate}
\item  If $C$ is aperiodic then $F_{B,D,n}\simeq F_{D,B,n}$;
and
\item If $C$ is periodic then $F_{B,D,n}\simeq\mathrm{res}_{\iota,k} F_{D,B,n}$.
\end{enumerate}
\item \label{lemma.6.7}
\emph{\cite[Corollary 2.2.24]{Ben2018}} Suppose $(B,D,n)\sim (B',D',n')$
in $\Sigma$.
\begin{enumerate}
\item  If $C$ is aperiodic then $G_{B,D,n}\simeq G_{B',D',n'}$.
\item  If $C$ is periodic and $C'=C[t]$ for some $t\in\mathbb{Z}$ then $G_{B,D,n}\simeq G_{B',D',n'}$. 
\item  If $C$ is periodic and $C'=C^{-1}[t]$ for some $t\in\mathbb{Z}$ then $G_{B,D,n}\simeq\mathrm{res}_{\iota,k} \,G_{B',D',n'}$.
\end{enumerate}
\end{enumerate}
\end{lemma}
\begin{proof}
(i) Let $M$ be an object in $\mathcal{K}_{\mathrm{min}}(\Lambda\text{-}\boldsymbol{\mathrm{Proj}}))$. If $C$ is aperiodic then $F^{\pm}_{B,D,n}(M)=F^{\pm}_{D,B,n}(M)$ as $R$-modules; and if $C$ is periodic the action of $T^{\pm1}$ on $F_{B,D,n}(M)$ depends on the order of $B$ and $D$, and must be exchanged with $T^{\mp 1}$ if and only if $B$ and $D$ are swapped.

(ii) Let $A$, $E$, $(d_{ \mathrm{l}( \gamma)}^{-1}\gamma A)^{-1}E$, $A'=\gamma^{-1}d_{ \mathrm{l}( \gamma)}E$ and  $E'=d_{ \mathrm{l}( \gamma)}^{-1}\gamma A$ be homotopy words. Then $G_{A,A',0}\simeq L\simeq G_{E',E,-1}$ where $L(M)=L^{+}(M)/L^{-}(M)$ and $
L^{\pm}(M)=e_{v}M^{n}\cap(\gamma A^{-}(M)+d_{ \mathrm{l}( \gamma)}E{}^{\pm}(M)\cap\gamma A^{+}(M))$ \cite[Lemma 2.2.22]{Ben2018} (see also the proof of the first lemma on \cite[p. 25]{Rin1975}). As in \cite[Corollary 2.2.23]{Ben2018}, iterating these natural isomorphisms gives the proof of (ii).
\end{proof}
\begin{definition}\label{9.444}
\cite[Definition 2.2.20]{Ben2018} Let $\Sigma(s)$ be the set of $(B,D,n)\in\Sigma$ where $B^{-1}D$ is aperiodic, and $\Sigma(b)$ the set of such $(B,D,n)$ where $B^{-1}D$ is periodic. Note that the relation $\sim$ on $\Sigma$ restricts to an equivalence relation $\sim_{s}$ (respectively $\sim_{b}$) on $\Sigma(s)$ (respectively $\Sigma(b)$). Let $\mathcal{I}(s)\subseteq\Sigma(s)$ (respectively $\mathcal{I}(b)\subseteq\Sigma(b)$) denote a chosen
collection of representatives $(B,D,n)$, one for each
equivalence class of $\Sigma(s)$ (respectively $\Sigma(b)$). Let $\mathcal{I}=\mathcal{I}(s)\sqcup\mathcal{I}(b)$. 

\cite[Definition 2.2.25]{Ben2018} Let $(B,D,n)\in\mathcal{I}$, $C=B^{-1}D$, $P=P(C)[\mu_{C}(a_{B,D})-n]$, $V$ and $V'$ be free $R$-modules with bases $(v_{\lambda}\mid\lambda\in\Omega)$ and $(v'_{\lambda'}\mid\lambda'\in\Omega')$ and let $f:V\rightarrow V'$ be $R$-linear. In this notation define $f_{\lambda',\lambda}\in R$ by $f( v_{\lambda})=\sum_{\lambda'} f_{\lambda',\lambda}v_{\lambda'}$. If $(B,D,n)\in\mathcal{I}(b)$ then $a_{B,D}=0$, $f$ is $R[T,T^{-1}]$-linear and $T$ defines automorphisms
$\varphi_{V}$ of $V$ and $\varphi_{V'}$ of $V'$
with $f\varphi_{V}=\varphi_{V'}f$.

If $(B,D,n)$ lies in $\mathcal{I}(s)$ (respectively $\mathcal{I}(b)$) we use $S_{B,D,n}$ to denote a functor  $R\text{-}\boldsymbol{\mathrm{Proj}}\rightarrow\mathcal{C}_{\mathrm{min}}(\Lambda\text{-}\boldsymbol{\mathrm{Proj}})$ (respectively $R[T,T^{-1}]\text{-}\boldsymbol{\mathrm{Mod}}_{R\text{-}\boldsymbol{\mathrm{Proj}}}\rightarrow\mathcal{C}_{\mathrm{min}}(\Lambda\text{-}\boldsymbol{\mathrm{Proj}})$), defined as follows. On objects $V$ define the homogeneous component $S^{m}_{B,D,n}(V)$ of the complex $S_{B,D,n}(V)$ in degree $m\in\mathbb{Z}$ by $P^{m}\otimes_{R}V$  (respectively $P^{m}\otimes_{R[T,T^{-1}]}V$). Define the corresponding differential $d^{m}_{S_{B,D,n}(V)}$ by $d^{m}_{P}\otimes_{R}\mathrm{id}_{V}$ (respectively $d^{m}_{P}\otimes_{R[T,T^{-1}]}\mathrm{id}_{V}$) in degree $m\in\mathbb{Z}$. Similarly we can define the map $S^{m}_{B,D,n}(f)$ of the image $S_{B,D,n}(f)$ of $S_{B,D,n}$ on a morphism $f$ by $\mathrm{id}_{P}^{m}\otimes_{R}f$ (respectively $\mathrm{id}_{P}^{m}\otimes_{R[T,T^{-1}]}f$).
\end{definition}
\begin{corollary}\label{cor9.6}
\emph{\cite[Corollary 2.2.26]{Ben2018}} 
 Suppose $(B,D,n)\sim (B',D',n')$
in $\Sigma$.
\begin{enumerate}
\item  If $C$ is aperiodic then $S_{B,D,n}\simeq S_{B',D',n'}$.
\item  If $C$ is periodic and $C'=C[t]$
for some $t\in\mathbb{Z}$ then $S_{B,D,n}\simeq S_{B',D',n'}$. 
\item  If $C$ is periodic and $C'=C^{-1}[t]$ for some $t\in\mathbb{Z}$ then $S_{B,D,n}\simeq S_{B',D',n'}\,\mathrm{res}_{\iota,R} $.
\end{enumerate}
\end{corollary}
\begin{proof} Recall $b_{i,C}$ denotes the coset of $e_{v_{C}(i)}$ in $P(C)$. By Lemma \ref{lemma.4.1} there is a canonical bijection $\omega:I_{C'}\rightarrow I_{C}$ such that $b_{i,C}\mapsto b_{\omega(i),C'}$ gives a natural isomorphism $\theta$ of complexes $P(C')[\mu_{C'}(a_{B',D'})-n']\rightarrow P(C)[\mu_{C}(a_{B,D})-n]$. If $C$ is aperiodic this shows $S_{B,D,n}\simeq S_{B',D',n'}$. If instead $C$ is aperiodic we have ($\theta( b _{i,C'}T)= b _{\omega(i),C}T$
when $C'=C[t]$) and ($\theta( b _{i,C'}T^{-1})= b _{\omega(i),C}T$
when $C'=C^{-1}[t]$), as required. 
\end{proof}
\begin{corollary}\label{allbandsgood}Let $(B,D,n)\sim(B',D',n')$, let $C$ be $p$-periodic and let $V$ and $W$ be objects in $R[T,T^{-1}]\text{-}\boldsymbol{\mathrm{Mod}}_{R\text{-}\boldsymbol{\mathrm{Proj}}}$. If $k\otimes_{R} V\simeq k\otimes W$ in $k[T,T^{-1}]\text{-}\boldsymbol{\mathrm{Mod}}$ then \emph{(}$S_{B,D,n}(V)\simeq S_{B',D',n'}(W)$ if $C'=C[t]$\emph{)} and \emph{(}$S_{B,D,n}(V)\simeq S_{B',D',n'}(\mathrm{res}_{\iota,R}(W))$ if $C'=C^{-1}[t]$\emph{)}.
\end{corollary}
\begin{proof}Since $p>0$, and by the definition of $\mu_{C}$, there exists $s\in \mathbb{Z}$ such that
\[
C[s]=\dots l_{-2}^{-1}r_{-2}l_{-1}^{-1}r_{-1}d_{\mathrm{l}(\alpha)}^{-1}\alpha \mid \beta^{-1} d_{\mathrm{l}(\beta)}l_{2}^{-1}r_{2}l_{3}^{-1}r_{3}\dots
\]
for some $\alpha,\beta\in\mathbf{P}$. Note that $C[s]$ is also $p$-periodic, and so  $C[s]={}^{\infty}{E}^{\,\infty}$ where $E=l_{1}^{-1}r_{1}\dots l_{p}^{-1}r_{p}$ is a homotopy word. Since $l_{p}^{-1}r_{p}l_{1}^{-1}r_{1}=d_{\mathrm{l}(\alpha)}^{-1}\alpha\beta^{-1} d_{\mathrm{l}(\beta)}$, and since $k\otimes_{R} V\simeq k\otimes W$, by Lemma \ref{bandsgood} we have $P(E,V)\simeq P(E,W)$ as complexes. By Remark \ref{indbasis} we have $P(C[s],V)\simeq P(E,V)$ and $P(C[s],W)\simeq P(E,W)$ as complexes. Setting $H=(E^{-1})^{\infty}$ and $L={E}^{\,\infty}$ gives $(B,D,n)\sim(H,L,n+\mu_{C}(s))$ (since $H^{-1}L=C[s]$), and so we have $S_{B,D,n}\simeq S_{H,L,n+\mu_{C}(s)}$ by Corollary \ref{cor9.6}(ii). Altogether we have
\[
\begin{array}{c}
S_{B,D,n}(V)\simeq S_{H,L,n+\mu_{C}(s)}(V)=P(C[s],V)[-(n+\mu_{C}(s))]\\
\simeq P(C[s],W)[-(n+\mu_{C}(s))]= S_{H,L,n+\mu_{C}(s)}(W)\simeq S_{B,D,n}(W),
\end{array}
\]
and so the result follows by parts (ii) and (iii) of Corollary \ref{cor9.6}.
\end{proof}
\section{\label{sec:Evaluation-on-string}Relations on complexes.}
 
\begin{assumption}\label{ass10}
\cite[Definition 2.3.1]{Ben2018} In \S\ref{sec:Evaluation-on-string} we let $C$ be any homotopy $I$-word. We fix a \textit{transversal} $\mathtt{S}$ of $\mathfrak{m}$ by choosing a lift $s\in R$
for each element of $k=R/\mathfrak{m}$ so that $\mathtt{S}\cap(r+\mathfrak{m})$ has precisely one element
for each $r\in R$. We also assume $\mathtt{S}\cap\mathfrak{m}=\{0\}$ and $\mathtt{S}\cap(1+\mathfrak{m})=\{1\}$.
\end{assumption}
\begin{notation}
\label{definition.1.1-1}\cite[Definition 2.3.3]{Ben2018} Let $i\in I$. Recall the symbol $b _{i}$
denotes the coset of $e_{v_{C}(i)}$ in the summand $\Lambda e_{v_{C}(i)}$
of $P^{\mu_{C}(i)}(C)$. Let $\mathbf{P}[i]$ be the set $\mathbf{P}(v_{C}(i)\rightarrow)$ of all non-trivial paths
$\sigma\notin\mathcal{J}$ with tail $v_{C}(i)$. If $x\in\mathbf{A}$ let $\mathbf{P}[x,i]=\{\sigma\in\mathbf{P}[i]\mid\mathrm{l}( \sigma)=x\}$.

Recall $\Lambda\mathfrak{m}\subseteq \rad(\Lambda)$. So, if $m\in P(C)$ then $m=\sum_{i}(\eta_{i} b _{i}+\sum_{\sigma\in\mathbf{P}[i]}r_{\sigma,i}\sigma b _{i})$
where $\eta_{i}\in \mathtt{S}$ and $r_{\sigma,i}\in R$ for each $i$, $\eta_{i}=r_{\sigma,i}=0$
for all but finitely many $i$, and given $i\in I$ fixed we have $r_{\sigma,i}=0$ for all but finitely many $\sigma$. In this case notation $\left\lceil m\right\rceil =\sum_{i}\eta_{i} b _{i}$
and $\left\lfloor m\right\rfloor =\sum_{i,\sigma}r_{\sigma,i}\sigma b _{i}$. Let $t\in I$. Let $\psi_{t}$ denote the $\Lambda$-module
epimorphism $P(C)=\bigoplus_{i}\Lambda e_{v_{C}(i)}\rightarrow\Lambda e_{v_{C}(t)}$
sending $m=\sum_{i}m_{i} b _{i}$ to $m_{t}$. For $m=\sum_{i}(\eta_{i} b _{i}+\sum_{\sigma\in\mathbf{P}[i]}r_{\sigma,i}\sigma b _{i})$
as above this gives $\psi_{t}(\gamma\left\lceil m\right\rceil )=\eta_{t}\gamma$
and $\psi_{t}(\gamma\left\lfloor m\right\rfloor )=\sum_{\sigma}r_{\sigma,t}\gamma\sigma$. 

Let $\left[m\right._{x,t}=\left\lceil m\right._{x,t}+\left\lfloor m\right._{x,t}$
and $\left.m\right]_{\,x,t}=\left.m\right\rceil _{x,t}+\left.m\right\rfloor _{x,t}$
where
\[
\begin{array}{c}
\left.m\right\rceil _{x,t}=\left\{ \begin{array}{cc}
\eta_{t+1}\alpha & (\text{if }t+1\in I\text{ and }l_{t+1}^{-1}r_{t+1}=\alpha^{-1}d_{x})\\
0 & (\text{otherwise})
\end{array}\right.\\
\left\lceil m\right._{x,t}=\left\{ \begin{array}{cc}
\eta_{t-1}\beta & (\text{if }t-1\in I\text{ and }l_{t}^{-1}r_{t}=d_{x}^{-1}\beta)\\
0 & (\text{otherwise})
\end{array}\right.\\
\left.m\right\rfloor _{x,t}=\left\{ \begin{array}{cc}
\sum_{\sigma\in\mathbf{P}[x,t+1]}r_{\sigma,t+1}\sigma\kappa & (\text{if }t+1\in I\text{ and }l_{t+1}^{-1}r_{t+1}=\kappa^{-1}d_{ \mathrm{l}( \kappa)})\\
0 & (\text{otherwise})
\end{array}\right.\\
\left\lfloor m\right._{x,t}=\left\{ \begin{array}{cc}
\sum_{\sigma\in\mathbf{P}[x,t-1]}r_{\sigma,t-1}\sigma\zeta & (\text{if }t-1\in I\text{ and }l_{t}^{-1}r_{t}=d_{ \mathrm{l}( \zeta)}^{-1}\zeta)\\
0 & (\text{otherwise})
\end{array}\right.
\end{array}
\]
For a vertex $v$ recall the sum $\sum_{a\in\mathbf{A}(v\rightarrow)}\Lambda a$
is direct. For any arrow $y$ with tail $v$ let $\theta_{y}:\bigoplus_{a}\Lambda a\rightarrow\Lambda y$
be the canonical projection. For $m'\in P(C)$
write $m'=\sum_{i}\eta'_{i} b _{i}+\sum_{i}\sum_{\sigma\in\mathbf{P}[i]}r'_{\sigma,i}\sigma b _{i}$ where $\eta'_{i}\in \mathtt{S}$ and $r'_{\sigma,i}\in R$ have finite support (as above).
\end{notation}
\begin{lemma}
\label{lemma30}\emph{\cite[Lemma 2.3.5]{Ben2018}} Fix $m,m'\in P(C)$ as in Notation \ref{definition.1.1-1}. Let $x\in\mathbf{A}$ and $t\in I$.
\begin{enumerate}
\item We have $\psi_{t}(d_{x,P(C)}(\left\lceil m\right\rceil ))=\left\lceil m\right._{x,t}+\left.m\right\rceil _{x,t}$
and $\psi_{t}(d_{x,P(C)}(\left\lfloor m\right\rfloor ))=\left\lfloor m\right._{x,t}+\left.m\right\rfloor _{x,t}$.
\item If $\gamma\in\mathbf{P}[x,t]$ and $m\in\gamma^{-1}d_{x}m'$
then $\theta_{ \mathrm{f}( \gamma)}(\psi_{t}(d_{x,P(C)}(m')))-\eta_{t}\gamma\in\gamma\mathrm{rad}(\Lambda)$.
\item If $l_{t+1}^{-1}r_{t+1}=\gamma^{-1}d_{ \mathrm{l}( \gamma)}$ then
$\theta_{ \mathrm{f}( \gamma)}(\left.m\right]_{\, \mathrm{l}( \gamma),t})-\eta_{t+1}\gamma\in\gamma\mathrm{rad}(\Lambda)$
and $\theta_{ \mathrm{f}( \gamma)}(\left[m\right._{\mathrm{l}( \gamma),t})=0$.
\item If $l_{t}^{-1}r_{t}=d_{ \mathrm{l}( \gamma)}^{-1}\gamma$ then $\theta_{ \mathrm{f}( \gamma)}(\left[m\right._{ \mathrm{l}( \gamma),t})-\eta_{t-1}\gamma\in\gamma\mathrm{rad}(\Lambda)$
and $\theta_{ \mathrm{f}( \gamma)}(\left.m\right]_{\,\mathrm{l}( \gamma),t})=0$.
\end{enumerate}
\end{lemma}
\begin{proof}
(i) Let $\iota_{x}:xP(C)\rightarrow\bigoplus yP(C)$
and $\pi_{x}:\bigoplus yP(C)\rightarrow xP(C)$ respectively denote
the natural $R$-module inclusions and projections. Note $d_{x,P(C)}(\sum\eta_{i} b _{i})=\sum_{i}\eta_{i}\iota_{x}(\pi_{x}( b _{i}^{-}+ b _{i}^{+}))$,
and if $\psi_{t}(b _{i}^{\pm})\neq 0$
then $i=t\mp1$. By case analysis $\psi_{t}(d_{x,P(C)}(\left\lceil m\right\rceil ))=\left\lceil m\right._{x,t}+\left.m\right\rceil _{x,t}$, and by Lemma \ref{lemma.1.1}, $d_{x,P(C)}(\sum_{i,\sigma}r_{\sigma,i}\sigma b _{i})=\sum_{i,(\sigma,z,i)}r_{\sigma,i}\sigma d_{z,P(C)}( b _{i})$
for each $i$ where the triples $(\sigma,z,i)$ run through all $\sigma\in\mathbf{P}[x,i]$
and all arrows $z$ with $ \mathrm{f}( \sigma)z\in\mathbf{P}$. For any such triple
$(\sigma,z,i)$ we have $\sigma d_{P(C)}( b _{i})=\sigma d_{z,P(C)}( b _{i})$
and so $\psi_{t}(d_{x,P(C)}(\left\lfloor m\right\rfloor ))=\left\lfloor m\right._{x,t}+\left.m\right\rfloor _{x,t}$. 

(ii) Since $\sum_{\sigma}r_{\sigma,t}\sigma\in\mathrm{rad}(\Lambda)$
we have $\theta_{ \mathrm{f}( \gamma)}(m)-\eta_{t}\gamma\in\gamma\rad(\Lambda)$
and so applying $\theta_{ \mathrm{f}( \gamma)}\psi_{t}$ to either side of $\gamma m=d_{x,P(C)}(m')$
gives $\theta_{ \mathrm{f}( \gamma)}\psi_{t}(d_{x,P(C)}(m'))-\eta_{t}\gamma\in\gamma\mathrm{rad}(\Lambda)$.

(iii) If $l_{t+1}^{-1}r_{t+1}=\gamma^{-1}d_{\mathrm{l}( \gamma)}$ then $\left[m\right._{\mathrm{l}( \gamma),t}=0$
unless $l_{t}^{-1}r_{t}=d_{ \mathrm{l}( \zeta)}^{-1}\zeta$ in which case $ \mathrm{f}( \gamma)\neq  \mathrm{f}( \zeta)$
since $d_{ \mathrm{l}( \zeta)}^{-1}\zeta\gamma^{-1}d_{\mathrm{l}( \gamma)}$ is a homotopy word. Furthermore
$\left.m\right\rceil _{ \mathrm{l}( \gamma),t}=\eta_{t+1}\gamma$ and $\left.m\right\rfloor _{ \mathrm{l}( \gamma),t}$ is equal to the sum over $\sigma\in\mathbf{P}[ \mathrm{l}( \gamma),t+1]$ of the terms $r_{\sigma,t+1}\sigma\gamma$. Hence $\left.m\right\rfloor _{ \mathrm{l}( \gamma),t}\in\rad(\Lambda)\gamma\cap  \mathrm{l}( \gamma)\Lambda$
which is contained in $\gamma\rad(\Lambda)$
by Corollary \ref{corollary.0.2}(iii). For (iv) apply the above to $D=C^{-1}$. 
\end{proof}
\begin{notation}\label{cidelta}\cite[Definition 2.3.7]{Ben2018} For
each $i$ the words $C_{>i}$ and $(C_{\leq i})^{-1}$ have head $v_{C}(i)$
and opposite sign, and we let $C(i,\delta)$ be the one with sign $\delta$. If $C(i,\delta)=C_{>i}$
then let $d_{i}(C,\delta)=1$, and otherwise $C(i,\delta)=(C_{\leq i})^{-1}$
in which case we let $d_{i}(C,\delta)=-1$. If $s\in I_{C(i,\delta)}$
and $s+1\in I_{C(i,\delta)}$ then ($C(i,\delta)_{s}=l_{i+s}^{-1}r_{i+s}$ if $d=1$) and ($C(i,\delta)_{s}=r_{i-s+1}^{-1}l_{i-s+1}$ if $d=-1$).
\end{notation}
\begin{corollary}
\label{corollary23}\emph{\cite[Corollary 2.3.8]{Ben2018}} Fix $m,m'\in P(C)$ as in Notation \ref{definition.1.1-1}. Let $i\in I$ and $\delta=\pm 1$, and let $d=d_{i}(C,\delta)$ as in Notation \ref{cidelta}.
\begin{enumerate}
\item If $n-1,n\in I_{C(i,\delta)}$ and $m\in C(i,\delta)_{n}m'$
then $\eta_{i+d(n-1)}=\eta'_{i+dn}$. 
\item If $I_{C(i,\delta)}=\{0,\dots,h\}$, $C(i,\delta) 1 _{u,\epsilon}=C(i,\delta)$ and $m\in 1 _{u,\epsilon}^{-}(P(C))$
then $\eta_{i+dh}=0$.
\item If $m\in C(i,\delta)^{-}(P(C))$ then $\eta_{i}=0$.
\end{enumerate}
\end{corollary}
\begin{proof}
(i) Let $C(i,\delta)_{n}=\gamma^{-1}d_{ \mathrm{l}( \gamma)}$, $t=i+d(n-1)$ and $x= \mathrm{l}( \gamma)$
so that $\gamma\in\mathbf{P}[x,t]$. By Lemma \ref{lemma30}(i) we have $\theta_{ \mathrm{f}( \gamma)}(\psi_{i+d(n-1)}(d_{x,P(C)}(m')))=\theta_{ \mathrm{f}( \gamma)}(\left[m'\right._{x,i+d(n-1)})+\theta_{ \mathrm{f}( \gamma)}(\left.m'\right]_{\,x,i+d(n-1)})$. By Lemma \ref{lemma30}(ii) this gives $\theta_{ \mathrm{f}( \gamma)}(\psi_{i+d(n-1)}(d_{x,P(C)}(m')))-\eta_{i+d(n-1)}\gamma\in\gamma\mathrm{rad}(\Lambda)$. 

If $d=1$ then $l_{i+n}^{-1}r_{i+n}=\gamma^{-1}d_{ \mathrm{l}( \gamma)}$, and applying Lemma \ref{lemma30}(iii) (where $t=i+n-1$ and $m$ is replaced with $m'$) gives  $\theta_{ \mathrm{f}( \gamma)}(\left.m'\right]_{\, x,i+n-1})-\eta'_{i+n}\gamma\in\gamma\mathrm{rad}(\Lambda)$
and $\theta_{ \mathrm{f}( \gamma)}(\left[m'\right._{x,i+n-1})=0$. This gives $\theta_{ \mathrm{f}( \gamma)}(\psi_{i+n-1}(d_{x,P(C)}(m)))-\eta'_{i+n}\gamma\in\gamma\rad(\Lambda)$, and so $(\eta'_{i+n}-\eta{}_{i+n-1})\gamma\in\gamma\mathrm{rad}(\Lambda)$. In case $d=-1$ we similarly have $(\eta'_{i-n}-\eta{}_{i-(n-1)})\gamma\in\gamma\mathrm{rad}(\Lambda)$, by applying Lemma \ref{lemma30}(iv) where $t=i-n+1$. 

In either case $(\eta'_{i+dn}-\eta{}_{i+d(n-1)})\gamma\in\gamma\mathrm{rad}(\Lambda)$, and if also $\eta'_{i+dn}-\eta{}_{i+d(n-1)}$ lies outside $\mathfrak{m}$
then $\gamma\Lambda=\gamma\rad(\Lambda)$
which contradicts Remark \ref{rem2.6}. Hence $\eta'_{i+dn}-\eta{}_{i+d(n-1)}\in\mathfrak{m}$
and as $\mathtt{S}$ is a transversal in $R$ with respect to $\mathfrak{m}$
this means $\eta'_{i+dn}=\eta{}_{i+d(n-1)}$. For the case where $C(i,\delta)_{n}=d_{ \mathrm{l}( \gamma)}^{-1}\gamma$
the proof is similar: when we use Lemma \ref{lemma30}(ii) we set $t=i+dn$ and swap $m$ and $m'$, and when we use Lemma \ref{lemma30}(iii) and Lemma \ref{lemma30}(iv) we set $t=i-n$ and $t=i+n$ respectively.

(ii) It suffices to prove $\eta_{i+dh}\in\mathfrak{m}$ since $\mathtt{S}\cap\mathfrak{m}=0$.
If there is no $\beta\in\mathbf{P}$ for which $ 1 _{u,\epsilon}\beta^{-1}d_{ \mathrm{l}( \beta)}$
is a word then $ 1 _{u,\epsilon}^{-}(P(C))\subseteq\rad(P(C))$
and so $\psi_{i+dh}(m)\in\rad(\Lambda e_{v_{C}(i+dn)})$.
Since $\eta_{i+dh} e_{v_{C}(i+dh)}=\psi_{i+dh}(m-\left\lfloor m\right\rfloor )$ this
gives $\eta_{i+dh}\in\mathfrak{m}$ as $\Lambda e_{v_{C}(i+dh)}$
is local. Suppose instead there is some $\beta\in\mathbf{P}$ for
which $ 1 _{u,\epsilon}\beta^{-1}d_{ \mathrm{l}( \beta)}$ is a word. By definition
$m\in\gamma^{-1}d_{ \mathrm{l}( \gamma)}m'$ for some $m'\in P(C)$ and some
$\gamma\in\mathbf{P}$ such that $C(i,\delta)\gamma^{-1}d_{ \mathrm{l}( \gamma)}$
is a homotopy word. 

Let $x=\mathrm{l}( \gamma)$. By Lemma \ref{lemma30}(i) $\psi_{i+dn}(d_{x,P(C)}(m'))=\left[m'\right._{x,i+dh}+\left.m'\right]_{\,x,i+dh}$.
Since $i+d(h+1)\notin I$, $d=1$ implies $\left.m'\right]_{\,x,i+h}=0$
and $d=-1$ implies $\left[m'\right._{x,i-h}=0$. If $d=1$ and $\left[m'\right._{x,i+h}\neq0$
then $i+h-1\in I$ and $l_{i-h+1}^{-1}r_{i-h+1}=\tau^{-1}d_{y}$ which
means $d_{y}^{-1}\tau\gamma^{-1}d_{x}$ is a homotopy word and hence
$\theta_{ \mathrm{f}( \gamma)}(\left[m'\right._{x,i+h})=0$. Similarly $\theta_{ \mathrm{f}( \gamma)}(\left.m'\right]_{\,x,i-h})=0$
when $d=-1$ , and altogether this gives $\theta_{ \mathrm{f}( \gamma)}(\psi_{i+dh}(d_{x,P(C)}(m')))=0$
so $\eta_{i+dh}\gamma\in\gamma\mathrm{rad}(\Lambda)$ by Lemma \ref{lemma30}(ii). As above this shows $\eta_{i+dh}\in\mathfrak{m}$, completing the proof of (ii).

(iii) Choose
$h\geq0$ such that ($I_{C(i,\delta)}=\mathbb{N}$ and $m\in C(i,\delta)_{\leq h}0$)
or ($I_{C(i,\delta)}=\{0,\dots,h\}$ and $C(i,\delta) 1 _{u,\epsilon}=C(i,\delta)$
for a vertex $u$ and $\epsilon\in\{\pm1\}$). So we have elements
$m_{j}=\sum_{i}\eta{}_{i,j} b _{i}+\sum_{i}\sum_{\sigma\in\mathbf{P}[i]}r_{\sigma,i,j}\sigma b _{i}$
from $P(C)$ where $m_{0}=m$ and $m_{j}\in C(i,\delta)_{j+1}m_{j+1}$
whenever $j<h$. By assumption when $C(i,\delta)$ is infinite, or
by (ii) when $C(i,\delta)$ is finite, we have $\eta{}_{i+dh,h}=0$.
Applying (i) to each natural number $j\leq h-1$ gives $\eta_{i+dj,j}=\eta_{i+d(j+1),j+1}$
and together this shows $\eta_{i}=\eta_{i,0}=\eta_{i+d,1}=\dots=\eta_{i+d(h-1),h-1}=0$.
\end{proof}
\begin{corollary}
\label{corollary31}\emph{\cite[Corollary 2.3.9]{Ben2018}} We have $ b _{i}\in C(i,\delta)^{+}(P(C))$ for any $i\in I_{C}$.
\end{corollary}
\begin{proof}Let $d=d_{i}(C,\delta)$. It is straightforward to show that, if $n-1,n\in I_{C(i,\delta)}$ then $ b _{i+d(n-1)}\in C(i,\delta)_{n} b _{i+dn}$. In particular, if $I_{C(i,\delta)}=\mathbb{N}$ then $ b _{i}\in C(i,\delta)_{1}b_{i+d}$, $ b _{i+d}\in C(i,\delta)_{2}b_{i+2d}$ and so on. So, if $I_{C(i,\delta)}=\mathbb{N}$ then the existence of the $\mathbb{N}$-sequence $(b_{i}\mid i\in\mathbb{N})$ in $P(C)$ shows that $ b _{i}\in C(i,\delta)^{+}(P(C))$. 

Now consider the case $I_{C(i,\delta)}=\{0,\dots,h\}$ and $C(i,\delta) 1 _{u,\epsilon}=C(i,\delta)$. As above we have $ b _{i}\in C(i,\delta) b _{i+dh}$, and so it suffices to show $b _{i+dh}\in 1 _{u,\epsilon}^{+}(P(C))$. It is enough to assume $C(i,\delta)d_{x}^{-1}x$ is a homotoy word for some arrow $x$, as otherwise $e_{u}P(C)= 1 _{u,\epsilon}^{+}(P(C))$.
Note $d_{P(C)}( b _{i+dh})= b _{i+dh}^{\pm}$ where $d=\mp 1$, and if $l_{i+dh}^{-1}r_{i+dh}=\beta^{-1}d_{ \mathrm{l}( \beta)}$ then
$ \mathrm{l}( \beta)\neq x$. So $d_{x,P(C)}( b _{i+dh})=0$ which
shows $ b _{i+dh}\in d_{x}^{-1}x0\subseteq 1 _{u,\epsilon}^{+}(P(C))$.
\end{proof}
\begin{notation}\label{nastynot}\cite[Definitions 2.3.10 and 2.3.12]{Ben2018} Let $V$ be a $R[T,T^{-1}]$-module with free $R$-basis $(v_{\lambda}\mid\lambda\in\Omega)$, $C={}^{\infty}E^{\,\infty}$ be periodic of period $p$ and $E=l_{1}^{-1}r_{1}\dots l_{p}^{-1}r_{p}$. If $\lambda\in\Omega$ and $0\leq i\leq p-1$ let $ b _{i,\lambda}= b _{i}\otimes v_{\lambda}$.
So $q\in P(C,V)$ gives $q=\sum_{i,\lambda}q_{i,\lambda} b _{i,\lambda}$
where $q_{i,\lambda}=\eta_{i,\lambda}+\sum_{\sigma}r_{\sigma,i,\lambda}\sigma$,
$\eta_{i,\lambda}\in \mathtt{S}$, $r_{\sigma,i,\lambda}\in R$ and $r_{\sigma,i,\lambda}=0$ for all but finitely many $\sigma\in\mathbf{P}[i]$. 

Let $\left\lceil q\right\rceil =\sum_{i,\lambda}\eta_{i,\lambda} b _{i,\lambda}$, $\left\lfloor q\right\rfloor =\sum_{i,\lambda}\sum_{\sigma}r_{\sigma,i,\lambda}\sigma b _{i,\lambda}$. For any $x\in\mathbf{A}$ we define $\left[q\right._{x,t}=\left\lceil q\right._{x,t}+\left\lfloor q\right._{x,t}$
and $\left.q\right]_{\,x,t}=\left.q\right\rceil _{x,t}+\left.q\right\rfloor _{x,t}$ by setting
\[
\begin{array}{c}
\left.q\right\rceil _{x,t}=\left\{ \begin{array}{cc}
\sum_{\lambda}\eta_{t+1,\lambda}\alpha\otimes v_{\lambda} & (\text{if }0\leq t<p-1\text{ and }l_{t+1}^{-1}r_{t+1}=\alpha^{-1}d_{x})\\
\sum_{\lambda}\eta_{0,\lambda}\alpha\otimes Tv_{\lambda} & (\text{if }t=p-1\text{ and }l_{p}^{-1}r_{p}=\alpha^{-1}d_{x})\\
0 & (\text{otherwise})
\end{array}\right.\\
\left\lceil q\right._{x,t}=\left\{ \begin{array}{cc}
\sum_{\lambda}\eta_{t-1,\lambda}\beta\otimes v_{\lambda} & (\text{if }0<t\leq p-1\text{ and }l_{t}^{-1}r_{t}=d_{x}^{-1}\beta)\\
\sum_{\lambda}\eta_{p-1,\lambda}\beta\otimes T^{-1}v_{\lambda} & (\text{if }t=0\text{ and }l_{0}^{-1}r_{0}=d_{x}^{-1}\beta)\\
0 & (\text{otherwise})
\end{array}\right.\\
\left.q\right\rfloor _{x,t}=\left\{ \begin{array}{cc}
\sum_{\lambda}\sum_{\sigma\in\mathbf{P}[x,t+1]}r_{\sigma,t+1,\lambda}\sigma\kappa\otimes v_{\lambda} & (\text{if }0\leq t<p-1\text{ and }l_{t+1}^{-1}r_{t+1}=\kappa^{-1}d_{ \mathrm{l}( \kappa)})\\
\sum_{\lambda}\sum_{\sigma\in\mathbf{P}[x,0]}r_{\sigma,0,\lambda}\sigma\kappa\otimes Tv_{\lambda} & (\text{if }t=p-1\text{ and }l_{p}^{-1}r_{p}=\kappa^{-1}d_{ \mathrm{l}( \kappa)})\\
0 & (\text{otherwise})
\end{array}\right.\\
\left\lfloor q\right._{x,t}=\left\{ \begin{array}{cc}
\sum_{\lambda}\sum_{\sigma\in\mathbf{P}[x,t-1]}r_{\sigma,t-1}\sigma\zeta\otimes v_{\lambda} & (\text{if }0<t\leq p-1\text{ and }l_{t}^{-1}r_{t}=d_{ \mathrm{l}( \zeta)}^{-1}\zeta)\\
\sum_{\lambda}\sum_{\sigma\in\mathbf{P}[x,p-1]}r_{\sigma,p-1,\lambda}\sigma\zeta\otimes T^{-1}v_{\lambda} & (\text{if }t=0\text{ and }l_{0}^{-1}r_{0}=d_{ \mathrm{l}( \zeta)}^{-1}\zeta)\\
0 & (\text{otherwise})
\end{array}\right.
\end{array}
\]
If $0\leq t<p$ let $\varphi_{t}:P(C,V)\rightarrow\Lambda e_{v_{C}(t)}\otimes_{R}V$
be the composition $\omega_{t}\kappa$ where $\kappa:P(C,V)\rightarrow\bigoplus_{i=0}^{p-1}\Lambda e_{v_{C}(i)}\otimes_{R}V$ is the isomorphism from Lemma \ref{lemma.2.3} and $\omega_{t}$
is the canonical projection. If $m\in P(C)$ and $v\in V$ then, in the sense of Notation \ref{definition.1.1-1}, we have $\varphi_{0}(d_{P(C,V)}(m\otimes v))=\psi_{0}(d_{P(C)}(m))\otimes v+\psi_{p}(d_{P(C)}(m))\otimes T^{-1}v$ and $\varphi_{p-1}(d_{P(C,V)}(m\otimes v))=\psi_{-1}(d_{P(C)}(m))\otimes Tv+\psi_{p-1}(d_{P(C)}(m))\otimes v$  \cite[Lemma 2.3.11]{Ben2018}. If $y\in\mathbf{A}(v\rightarrow )$ let $\phi_{y}:(\bigoplus_{a\in\mathbf{A}(v\rightarrow)}\Lambda a)\otimes_{R}V\rightarrow\Lambda a\otimes_{R}V$
be the natural $\Lambda$-module projection. 
\end{notation}
The proofs of Lemma \ref{lemma32} and Lemma \ref{lemma34} are ommited. The proof of Lemma \ref{lemma32} involves a straightforward application of Lemmas \ref{lemma.2.3} and \ref{lemma30}. The proof of Lemma \ref{lemma34} uses Lemma \ref{lemma32}, and is similar to the proof of Corollaries \ref{corollary23} and \ref{corollary31}. See \cite[\S 2.3.1]{Ben2018} for details.
\begin{lemma}\emph{\cite[Lemma 2.3.13]{Ben2018}} 
\label{lemma32} Let $x\in\mathbf{A}$, $0\leq t<p$ and $M=\rad(\Lambda)\otimes_{R} V$. Fix $q\in P(C,V)$ and consider Notation \ref{nastynot}.
\begin{enumerate}
\item We have $\varphi_{t}(d_{x,P(C,V)}(\left\lceil q\right\rceil ))=\left\lceil q\right._{x,t}+\left.q\right\rceil _{x,t}$
and $\varphi_{t}(d_{x,P(C,V)}(\left\lfloor q\right\rfloor ))=\left\lfloor q\right._{x,t}+\left.q\right\rfloor _{x,t}$.
\item If $\gamma\in\mathbf{P}[x,t]$ and $q\in\gamma^{-1}d_{x}q'$
then $\phi_{ \mathrm{f}( \gamma)}(\varphi_{t}(d_{x,P(C,V)}(q')))-\sum_{\lambda}\eta_{t,\lambda}\gamma\otimes v_{\lambda}\in\gamma M$.
\item If $l_{t+1}^{-1}r_{t+1}=\gamma^{-1}d_{x}$ then $\phi_{ \mathrm{f}( \gamma)}(\left[q\right._{x,t})=0$.
\item If $l_{t}^{-1}r_{t}=d_{x}^{-1}\gamma$ then $\phi_{ \mathrm{f}( \gamma)}(\left.q\right]_{\,x,t})=0$.
\item If $l_{t+1}^{-1}r_{t+1}=\gamma^{-1}d_{x}$ and $t<p-1$ then $\phi_{ \mathrm{f}( \gamma)}(\left.q\right]_{\,x,t})-\sum_{\lambda}\eta_{t+1,\lambda}\gamma\otimes v_{\lambda}\in\gamma M$.
\item If $l_{t}^{-1}r_{t}=d_{x}^{-1}\gamma$ and $0<t$ then $\phi_{ \mathrm{f}( \gamma)}(\left[q\right._{x,t})-\sum_{\lambda}\eta_{t-1,\lambda}\gamma\otimes v_{\lambda}\in\gamma M$.
\item If $l_{p}^{-1}r_{p}=\gamma^{-1}d_{x}$ then $\phi_{ \mathrm{f}( \gamma)}(\left.q\right]_{\,x,p-1})-\sum_{\lambda}\eta_{0,\lambda}\gamma\otimes Tv_{\lambda}\in\gamma M$.
\item If $l_{0}^{-1}r_{0}=d_{x}^{-1}\gamma$ then $\phi_{ \mathrm{f}( \gamma)}(\left[q\right._{x,0})-\sum_{\lambda}\eta_{p-1,\lambda}\gamma\otimes T^{-1}v_{\lambda}\in\gamma M$.
\end{enumerate}
\end{lemma}
\begin{lemma}
\label{lemma34}\emph{\cite[Lemmas 2.3.14 and 2.3.15]{Ben2018}} Let $i,n\in\mathbb{Z}$, $0\leq i \leq p-1$, $1\leq n\leq p$, $\delta=\pm 1$, $d=d_{i}(C,\delta)$ and let $\mu\in\Omega$. Fix $q,q'\in P(C,V)$ with $q\in C(i,\delta)_{n}q'$ and consider Notation \ref{nastynot}. 
\begin{enumerate}
\item If \emph{(}$i<p-n$ and $d=1$\emph{)} or \emph{(}$i>n-1$ and $d=-1$) then
$\eta_{i+d(n-1),\mu}=\eta'_{i+dn,\mu}$.
\item If \emph{(}$i>p-n$ and $d=1$\emph{)} or \emph{(}$i<n-1$ and $d=-1$\emph{)},
then $\eta_{i+d(n-p-1),\mu}=\eta'_{i+d(n-p),\mu}$.
\item If $i=p-n$ and $d=1$ then \emph{(}$\eta'_{0,\lambda}=0$ for all $\lambda$
if and only if $\eta_{p-1,\lambda}=0$ for all $\lambda$\emph{)}.
\item If $i=n-1$ and $d=-1$ then \emph{(}$\eta{}_{0,\lambda}=0$ for all $\lambda$
if and only if $\eta'_{p-1,\lambda}=0$ for all $\lambda$\emph{)}.
\end{enumerate}
Consequently $b_{i,\lambda}\in C(i,\delta)^{+}(P(C,V))$,
and if $q\in C(i,\delta)^{-}(P(C,V))$ then $\{\eta_{i,\lambda}\mid\lambda\in\Omega\}=\{0\}$.
\end{lemma}
\section{Refining complexes.}\label{refinecomplex}
\begin{assumption}
As in Assumption \ref{ass10}, in \S\ref{refinecomplex} we let $C$ be a homotopy $I$-word, and we let $\mathtt{S}$ be a transversal of $\mathfrak{m}$ in $R$ such that $\mathtt{S}\cap\mathfrak{m}=\{0\}$ and $\mathtt{S}\cap(1+\mathfrak{m})=\{1\}$.
\end{assumption}
\begin{corollary}\label{corollary25}\label{corollary24} Let $A\in\mathcal{W}_{v,\delta}$ and $(I,A,+)$ (respectively $(I,A,-)$) be the set of $i\in I$
with $v_{C}(i)=v$ and $C(i,\delta)\leq A$ (respectively $C(i,\delta)<A$). 
\begin{enumerate}
\item \emph{\cite[Corollary 2.3.16]{Ben2018} (}see also \emph{\cite[Lemma 8.1]{Cra2018})}. If $C$ is aperiodic then
\[
A^{\pm}(P(C))+e_{v}\rad(P(C))=\sum_{i\in(I,A,\pm)}R b _{i}+e_{v}\rad(P(C)).
\]
\item \emph{\cite[Corollary 2.3.17]{Ben2018} (}see also \emph{\cite[Lemma 8.4]{Cra2018})}. If $C$ is $p$-periodic then
\[
A^{\pm}(P(C,V))+e_{v}\rad(P(C,V))=\sum_{\lambda\in\Omega,\,i\in(\mathbb{Z},A,+)\,\mid\, 0\leq i\leq p-1}R b _{i,\lambda}+e_{v}\rad(P(C,V)).
\]
\end{enumerate}
\end{corollary}
\begin{proof}Recall $m\in P(C)$ satisfies $m=\sum_{i}(\eta_{i} b _{i}+\sum_{\sigma}r_{\sigma,i}\sigma b _{i})$
where $\eta_{i}\in \mathtt{S}$, $r_{\sigma,i}\in R$, $\eta_{i}=r_{\sigma,i}=0$
for all but finitely many $i\in I$, and $r_{\sigma,i}=0$ for all
but finitely many $\sigma\in\mathbf{P}[i]$. By Corollary \ref{corollary31} we have $ b _{i}\in C(i,\delta)^{+}(P(C))$
for each $i\in I$. So if $\sum_{i}\eta_{i} b _{i}$ lies
in the $R$-span of all $ b _{i}$ with $i\in(I,A,+)$ then $\sum_{i}\eta_{i} b _{i}\in A^{+}(P(C))$.
If we also have $\eta_{i}=0$ for $C(i,\delta)=A$ then $\sum_{i}\eta_{i} b _{i}\in A^{-}(P(C))$. This gives the containment $A^{\pm}(P(C))\supseteq\sum_{i\in(I,A,\pm)}R b _{i}$.

If $m\in A^{\pm}(P(C))$ then $m-\sum_{i}\eta_{i} b _{i}\in e_{v}\rad(P(C))$
for some $\eta_{i}\in \mathtt{S}$. For $i\in I$
such that $v_{C}(i)\neq v$ we have $\eta_{i}=0$ since $\sum_{i}\eta_{i} b _{i}\in e_{v}P(C)$.
If $m\in A^{+}(P(C))$ then given any $i\in I$ with $C(i,\delta)>A$
we have $m\in C(i,\delta)^{-}(P(C))$ and so $\eta_{i}=0$ by Corollary \ref{corollary23}. Similarly $m\in A^{-}(P(C))$ implies
$\eta_{i}=0$ given $C(i,\delta)\geq A$. This shows
$A^{\pm}(P(C))\subseteq\sum_{i\in(I,A,\pm)}R b _{i}+e_{v}\rad(P(C))$. This gives (i). The proof of (ii) is similar, uses Lemma \ref{lemma34}, and is omitted.
\end{proof}

\begin{notation}\cite[Definition 2.3.18]{Ben2018}
\label{definition2.3.4} Let  $\Xi:\mathcal{C}_\mathrm{min}(\Lambda\text{-}\boldsymbol{\mathrm{Proj}})\rightarrow\mathcal{K}_\mathrm{min}(\Lambda\text{-}\boldsymbol{\mathrm{Proj}})$ be the restriction of the quotient functor $\mathcal{C}(\Lambda\text{-}\boldsymbol{\mathrm{Proj}})\rightarrow\mathcal{K}(\Lambda\text{-}\boldsymbol{\mathrm{Proj}})$. 
\end{notation}
Recall the set $\Sigma$ of $(B,D,n)$ where $n\in\mathbb{Z}$ and $B$, $D$ and $B^{-1}D$ are homotopy words. Recall the functors $S_{B,D,n}$ and $F_{B,D,n}$ ($(B,D,n)\in\Sigma$). See \S\ref{Refined Functors for Complexes.} and \S\ref{subsec:Natural-Isomorphisms.} for details.
\begin{lemma}
\label{lemma.35}\emph{\cite[Lemma 2.3.19]{Ben2018} (}see also the second lemma of \emph{\cite[p. 26]{Rin1975})}. Let $V$ be a free $R$-module, $P=\Xi (S_{B,D,n}(V))$, $n\in\mathbb{Z}$, and
$B$ and $D$ be homotopy words with $C=B^{-1}D$. 
\begin{enumerate}
\item The map $\Phi_{V}:k\otimes_{R}V\rightarrow\bar{F}_{B,D,n}(P)$, $1\otimes v\mapsto  b _{a_{B,D}}\otimes v+\bar{F}_{B,D,n}^{-}(P)$
is injective.
\item If $C$ is aperiodic then $\Phi$ gives a natural map $k\otimes_{R}-\rightarrow\bar{F}_{B,D,n}\,\Xi\,  S_{B,D,n}$
\item If $C$ is periodic then $\Phi$ induces a natural map
$k[T,T^{-1}]\otimes_{R[T,T^{-1}]}-\rightarrow\bar{F}_{B,D,n}\,\Xi\,  S_{B,D,n}$.
\end{enumerate}
\end{lemma}
\begin{proof}
Let $i=a_{B,D}$. For any $\delta=\pm 1$ let $(i,\delta,\pm)=(I,C(i,\delta),\pm)$ as in the notation from Corollary \ref{corollary25}. Let $(v_{\lambda}\mid\lambda\in\Omega)$ be an $R$-basis
of $V$, and so $((1+\mathfrak{m})\otimes v_{\lambda}\mid\lambda\in\Omega)$ is a $k$-basis of
$k\otimes_{R}V$. Since $\Lambda$ is rad-nilpotent modulo $\mathfrak{m}$ we have $\mathrm{rad}(P(C))\otimes V\subseteq\bar{F}_{B,D,n}^{-}(P)$, and so  $\Phi_{V}$ is well defined. Since $V$ is free, if $\sum_{\lambda}s_{\lambda} b _{i}\otimes v_{\lambda}\in\bar{F}_{B,D,n}^{-}(P)$ for some $s_{\lambda}\in R$ then $s_{\lambda} b _{i}\in\bar{F}_{B,D,n}^{-}(P)$
for each $\lambda$. By Corollary \ref{corollary25} there exists $r_{j,\lambda}\in R$ such that $s_{\lambda} b _{i}-\sum_{j}r_{j,\lambda} b _{j}\in\mathrm{rad}(P(C))$ where $j$ runs through
the union of $(i,1,+)\cap(i,-1,-)$ and $(i,1,-)\cap(i,-1,+)$. Since $i\notin(i,-1,-)\cup(i,1,-)$, if $s_{\lambda}\notin\mathfrak{m}$
then $s_{\lambda}$ is a unit in which case $\Lambda e_{v_{C}(i)}\subseteq\mathrm{rad}(\Lambda e_{v_{C}(i)})$, a contradiction. Thus $\Phi_{V}$ is injective. The proof of part (ii) is straightforward. The proof of part (iii) is a straightforward application of part (ii), noting that there is an natural isomorphism $k[T,T^{-1}]\otimes_{R[T,T^{-1}]}-\simeq k\otimes_{R}-$ of functors $R[T,T^{-1}]\text{-}\boldsymbol{\mathrm{Mod}}_{R\text{-}\boldsymbol{\mathrm{Proj}}}\rightarrow k[T,T^{-1}]\text{-}\boldsymbol{\mathrm{Mod}}$.
\end{proof}
\begin{lemma}
\emph{\cite[Lemma 2.3.20]{Ben2018} (}see also \emph{\cite[Lemma 8.2]{Cra2018})}.\label{lemma.31} Let $(B,D,n),(B',D',n')\in\Sigma$ such that $C=B^{-1}D$ is aperiodic. Let $C'=(B')^{-1}D'$ and $P=P(C)$. 
\begin{enumerate}
\item If $i\in I$ then $\bar{F}_{C(i,1),C(i,-1),n}^{+}(P[\mu_{C}(i)-n])=\bar{F}_{C(i,1),C(i,-1),n}^{-}(P[\mu_{C}(i)-n])+R b _{i}$.
\item If $C'=C$ and $n-n'=\mu_{C}(a_{B,D})-\mu_{C}(a_{B',D'})$
then $k\otimes_{R}-\simeq F_{B',D',n'}\,\Xi\, S_{B,D,n}$.
\item If $(B,D,n)$ is not equivalent to $(B',D',n')$ then
$\bar{F}_{B',D',n'}(P[\mu_{C}(a_{B,D})-n])=0$.
\end{enumerate}
\end{lemma}
\begin{proof}
(i) Let $M=P[\mu_{C}(i)-n]$ and $v=v_{C}(i)$. Let $A=C(i,1)$ and $E=C(i,-1)$. Let $(i,1,\pm)=(I,A,\pm)$ and $(i,-1,\pm)=(I,E,\pm)$ as in the proof of Lemma \ref{lemma.35}. By Corollary \ref{corollary25} $\bar{F}_{A,E,n}^{+}(M)$ contains $\bar{F}_{A,E,n}^{-}(M)+Rb_{i}$. Now let $m\in\bar{F}_{A,E,n}^{+}(M)$. By assumption and
by Corollary \ref{corollary25} we may write $m=\sum_{j}\eta_{j} b _{j}+m_{0}$
for $\eta_{j}\in \mathtt{S}$ and some $m_{0}\in e_{v}\mbox{rad}(M)$, where
$\eta_{j}=0$ whenever $C(j,1)>A$ or $C(j,-1)>E$. Since
$m\in M^{n}=P^{\mu_{C}(i)}(C)$ we have $\eta_{j}=0$ for all
$j\in I$ with $\mu_{C}(j)\neq\mu_{C}(i)$. So $\sum_{j}\eta_{j} b _{j}\in \sum_{t}R b _{t}$ where $t$ runs through $(i,1,+)\cap(i,-1,+)$.

Let $(i,\delta,=)=\{j\in I\mid C(j,\delta)=C(i,\delta)\}$. Then $(i,1,+)\cap(i,-1,+)$ is the union of $(i,1,+)\cap(i,-1,-)$,
$(i,1,-)\cap(i,-1,+)$ and $(i,1,=)\cap(i,-1,=)$. So by Corollary
\ref{corollary25} $\sum_{j}\eta_{j} b _{j}$ lies in $F_{A,E,n}^{-}(M)+\sum_{t}R b _{t}$
where $t$ runs over $j\in(i,1,=)\cap(i,-1,=)$ with $\mu_{C}(j)=\mu_{C}(i)$. It suffices to assume $t\in I$,
$C(t,1)=A$, $C(t,-1)=C(i,-1)$ and $\mu_{C}(t)=\mu_{C}(i)$: and show $t=i$. If $C_{>i}=(C_{\leq t})^{-1}$ and $(C_{\leq i})^{-1}=C_{>t}$ then
$C[t]=C^{-1}[i]$ which contradicts \cite[Lemma 2.2.17]{Ben2018}. Hence $C_{>i}=C_{>t}$ and $C_{\leq i}=C_{\leq t}$
which shows $C=C[t-i]$. Applying   Lemma \ref{lemma.4.1}(iii) twice
yields $\mu_{C}(t-i)=0$, which means $t=i$ since $C$ is aperiodic.

(ii) Any element of $S_{B,D,n}(V)$ may be written as the coset
of a sum of pure tensors $\sum_{t=1}^{n}r_{t} b _{i}\otimes v_{t}+\bar{F}_{B,D,n}^{-}(\left.\Xi (S_{B,D,n}(V)))\right .$
for some $v_{1},\dots,v_{n}\in V$. Hence the $k$-linear embedding
$\Phi_{V}$ from Lemma \ref{lemma.35}(i) is surjective, and so $\Phi$ from Lemma \ref{lemma.35}(ii)
is a natural isomorphism. Since $(B,D,n)\sim(B',D',n')$ we have $F_{B',D',n'}\simeq\bar{F}_{B,D,n}$, and so $F_{B',D',n'}\,\Xi\,S_{B,D,n}\simeq k\otimes_{R}-$ as functors $R\text{-\textbf{Proj}}\rightarrow k\text{-\textbf{Mod}}$.

(iii) It is enough to show $\bar{F}_{B',D',n'}(P(C)[\mu_{C}(a_{B,D})-n])=0$. It suffices to let $s(B')=1$ and
$s(D')=-1$, and show $(B,D,n)\sim(B',D',n')$ assuming $\bar{F}_{B',D',n'}(P)\neq 0$.
By Corollary \ref{corollary25} the $R$-submodules $\bar{G}_{B',D',n'}^{\pm}(P)$
are spanned by sets of elements of the form $b_{i}$ together with
$\mbox{rad}(P)$, and hence, as in the proof of \cite[Lemma 8.2]{Cra2018}, the existence of $ b _{i}$ shows that $B'=C(i,1)$ and $D'= C(i,-1)$, so $B'^{-1}D'=C$. Since $ b _{i}$ lies in both $P^{\mu_{C}(i)}(C)$ and $P^{n'}$ we have $\mu_{C}(i)=n'+\mu_{C}(a_{B,D})-n$ and so $(B,D,n)\sim(B',D',n')$.
\end{proof}
\begin{lemma}
\label{lemma.7.3-1}\emph{\cite[Lemma 2.3.21]{Ben2018} (}see also \emph{\cite[Lemma 8.5]{Cra2018})}. Let $(B,D,n),(B',D',n')\in\Sigma$ where $C=B^{-1}D$ is $p$-periodic, and let $i\in\{0,\dots,p-1\}$. If $P=P(C,V)$ then:
\begin{enumerate}
\item $\bar{F}_{C(i,1),C(i,-1),n}^{+}(P[\mu_{C}(i)-n])=\bar{F}_{C(i,1),C(i,-1),n}^{-}(P[\mu_{C}(i)-n])+\sum_{\lambda}R b _{i,\lambda}$;
\item if $C'=C[m]$ and $n-n'=\mu_{C}(m)$ then $k[T,T^{-1}]\otimes_{R[T,T^{-1}]}-\simeq F_{B',D',n'}\Xi S_{B,D,n}$;
and
\item if $(B,D,n)$ is not equivalent to $(B',D',n')$ then
$\bar{F}_{B',D',n'}(P[-n])=0$.
\end{enumerate}
\end{lemma}
\section{\label{sec:Linear-compactness-and}Compactness and Covering.}
 
\begin{assumption}\label{ass11}
In \S\ref{sec:Linear-compactness-and} we let $v$ be a vertex and fix any homotopically minimal complex $M$ of projective $\Lambda$-modules. Use $M$ to denote the underlying $\Lambda$-module or underlying $R$-module.

The topology we refer to will be the $\mathfrak{m}$-adic topology. Let $\mathfrak{m}^{0}=R$. Recall a base of open sets for an $R$-module $N$ with this topology is the collection of cosets $m+\mathfrak{m}^{n}U$ where $n\in\mathbb{N}$ and $U$ is an $R$-submodule of $N$.
\end{assumption}
Later in this article we apply the results in \S\ref{sec:Linear-compactness-and} in the context where the object $M$ of $\mathcal{K}_{\mathrm{min}}(\Lambda\text{-}\boldsymbol{\mathrm{Proj}})$ lies in $\mathcal{K}_{\mathrm{min}}(\Lambda\text{-}\boldsymbol{\mathrm{proj}})$ (and so $M^{i}$ is finitely generated over $\Lambda$ for each $i\in\mathbb{Z}$).
\begin{remark}\label{fingenrem}
Recall from Remark \ref{rem2.6} the, by Definition \ref{def.1.13}, $\Lambda$ is a finitely generated $R$-module. Consequently, under Assumption \ref{ass11}, if $M$ lies in $\mathcal{K}_{\mathrm{min}}(\Lambda\text{-}\boldsymbol{\mathrm{proj}})$ each $R$-module of the form $e_{v}M^{i}$ is finitely generated. 
\end{remark}
\begin{definition}\cite[\S 1, p. 80, Definition]{Zel1953} Let $L$ be a subset of an $R$-module $N$. We write $L\subseteq_{\mathrm{c}}N$ if and only if $L$ is closed. We say $L$ is a \textit{linear variety} if $L=U+m\subseteq_{\mathrm{c}}N$ for some $R$-submodule $U$ of $N$. We say $N$ is \textit{linearly
compact} if any collection of linear varieties in $N$ with the
finite intersection property must have a non-void intersection. 
\end{definition}
\begin{example}\label{example.2.4.2} For each non-zero $n\in\mathbb{N}$ the $R$-module $R/\mathfrak{m}^{n}$ has the minimum condition on closed submodules, and so it is linearly compact by \cite[Proposition 5]{Zel1953}. Since $R$ is $\mathfrak{m}$-adically complete, $R$ (as a module over itself) is isomorphic to the inverse limit of a system of linearly compact $R$-modules. By \cite[Proposition 4]{Zel1953} this means $R$ is linearly compact. By \cite[Proposition 1]{Zel1953} this means that any finitely generated
free $R$-module is linearly compact. By \cite[Proposition 2]{Zel1953} this means any finitely generated $R$-module is linearly compact.
\end{example}
\begin{lemma}
\label{lemma.5.2-1}\emph{\cite[Lemma 2.4.1]{Ben2018}}. Suppose $M$ lies in $\mathcal{K}_{\mathrm{min}}(\Lambda\text{-}\boldsymbol{\mathrm{proj}})$. Let $r\in\mathbb{Z}$ and $\delta=\pm1$. Let $U$ be an
$R$-submodule of $e_{v}M^{r}$ with $e_{v}\rad(M^{r})\subseteq U$.
\begin{enumerate}
\item \emph{(}See also \emph{\cite[Lemma 10.4]{Cra2018})}. If $H$ is a linear variety in $e_{v}M^{r}$
and $m\in H\setminus U$, then there is a homotopy word $C\in\mathcal{W}_{v,\delta}$
such that $H\cap(U+m)$ meets $C^{+}(M)$ but not $C^{-}(M)$.
\item \emph{(}See also \emph{\cite[Lemma 10.5]{Cra2018})}. If $m\in e_{v}M^{r}\setminus U$ then
there are words $B\in\mathcal{W}_{v,\delta}$ and $D\in\mathcal{W}_{v,-\delta}$
such that $U+m$ meets $G_{B,D,r}^{+}(M)$ but not $G_{B,D,r}^{-}(M)$. 
\end{enumerate}
\end{lemma}
The proof of Lemma \ref{lemma.5.2-1} is given at the end of \S\ref{sec:Linear-compactness-and}. 
\begin{lemma}
\label{lemma.5.3}\emph{\cite[Lemma 2.4.8]{Ben2018} (}see also \emph{\cite[Lemma 10.3]{Cra2018})}. Fix an integer
$r$ and some $\delta\in\{\pm1\}$. For any non-empty subset $S$
of $e_{v}M^{r}$ which does not meet \emph{$\rad(M)$} there
is a homotopy word $C\in\mathcal{W}_{v,\delta}$ such that either:
\begin{enumerate}
\item $C$ is finite and $S$ meets $C^{+}(M)$ but not $C^{-}(M)$; or
\item $C$ is a homotopy $\mathbb{N}$-word and $S$ meets $C_{\leq n}M$
but not $C_{\leq n}\rad(M)$ for each $n\geq0$.
\end{enumerate}
\end{lemma}
In Lemma \ref{lemma.5.3} we do \textit{not} require that $M$ is an object of $\mathcal{K}_{\mathrm{min}}(\Lambda\text{-}\boldsymbol{\mathrm{proj}})$. In Lemma \ref{lemma.5.3-1} we \textit{do} consider this setting, and show $S\cap C^{+}(M)\neq\emptyset =S\cap C^{-}(M)$ in case (ii) of Lemma \ref{lemma.5.3}.
\begin{proof}[of  Lemma \ref{lemma.5.3}]
We assume either $S\cap B^{+}(M)=\emptyset$ or $S\cap B^{-}(M)\neq\emptyset$ for any finite homotopy word $B\in\mathcal{W}_{v,\delta}$. We now construct a homotopy $\mathbb{N}$-word $C$ iteratively from $C_{\leq0}= 1 _{v,\delta}$,
such that $S$ meets $C_{\leq n}M$ but not $C_{\leq n}\mbox{rad}(M)$
for each $n\geq0$. If $n=0$ there is nothing to prove. Assume $S$ meets $C_{\leq m}M$ but not $C_{\leq m}\mbox{rad}(M)$ for some fixed $m\geq0$. It suffices to choose $l_{m+1}$
and $r_{m+1}$ such that $S$ meets $C_{\leq m}l_{m+1}^{-1}r_{m+1}M$
but not $C_{\leq m}l_{m+1}^{-1}r_{m+1}\mbox{rad}(M)$. 

Suppose $S\cap(C_{\leq m})^{-}(M)\neq\emptyset$. We can assume there exists $y\in\mathbf{A}$ where $C_{\leq m}y^{-1}d_{y}$
is a homotopy word, since otherwise $(C_{\leq m})^{-}(M)\cap S\subseteq C_{\leq m}\mathrm{rad}(M)\cap S=\emptyset$. As $S$ meets $(C_{\leq m})^{-}(M)$ there exists $\gamma\in\mathbf{P}$
of minimal length for which $S$ meets $C_{\leq m}\gamma^{-1}d_{ \mathrm{l}( \gamma)}M$.
Let $l_{m+1}=\gamma$ and $r_{m+1}=d_{ \mathrm{l}( \gamma)}$. For case (ii) it suffices to show $S$ does not meet $C_{\leq m}\gamma^{-1}d_{\gamma}\rad(M)$. 

If $\gamma\in\mathbf{A}$ then $\gamma^{-1}d_{\gamma}\rad(M)\subseteq e_{t(\gamma)}\rad(M)$
by Corollary \ref{corollary.3.1-1} and so $S$ does not meet $C_{\leq m}\gamma^{-1}d_{\gamma}\rad(M)$. So we can assume $\gamma= \mathrm{l}( \gamma)\alpha$
for some $\alpha\in\mathbf{P}$. By Corollary \ref{corollary.3.1-1} we also have $\alpha^{-1}d_{\mathrm{l}( \alpha)}M=\alpha^{-1} \mathrm{l}( \gamma)^{-1} \mathrm{l}( \gamma)d_{ \mathrm{l}( \alpha)}M$
and by Lemma \ref{lemma.1.1}(iia) we have $ \mathrm{l}( \gamma)d_{ \mathrm{l}( \alpha)}M=d_{ \mathrm{l}( \gamma)} \mathrm{l}( \gamma)M=d_{ \mathrm{l}( \gamma)}\rad(M)$.
The minimality of the length of $\gamma$ shows that $S$ does not
meet $C_{\leq m}\alpha^{-1}d_{ \mathrm{l}( \alpha)}M$ and altogether this shows
$S$ does not meet $C_{\leq m}\gamma^{-1}d_{ \mathrm{l}( \gamma)}\rad(M)$. 

Suppose instead $S\cap(C_{\leq m})^{-}(M)=\emptyset$. Here $S\cap(C_{\leq m})^{+}(M)=\emptyset$ by the assumption at the beginning of the proof, so there is some arrow $x$ for which $C_{\leq m}d_{x}^{-1}x$
is a homotopy word. By definition $d_{x,M}$ sends elements of $e_{h(x)}M$ to $xM$, and so $S$ meets $C_{\leq m}M=C_{\leq m}d_{x}^{-1}xM$. Suppose the set $L$ of $\lambda\in\mathbf{P}$ where ($C_{\leq m}d_{x}^{-1}\lambda$
is a homotopy word) is infinite. By \cite[Lemma 2.1.19]{Ben2018} we have $\bigcap_{\lambda}C_{\leq m}d_{x}^{-1}\lambda M=(C_{\leq m})^{+}(M)$
which does not meet $S$. By Corollary \ref{corollary.3.2-1} there is some maximal
length $\mu\in L$ for which $S$ meets $C_{\leq m}d_{x}^{-1}\mu M$.
As $L$ is infinite $\mu\eta\in L$ for some arrow $\eta$ in which
case $C_{\leq m}d_{x}^{-1}\mu\rad(M)=C_{\leq m}d_{x}^{-1}\mu\eta M$
which does not meet $S$ by construction. In this case it is sufficient
to let $l_{m+1}=d_{ \mathrm{l}( \mu)}$ and $r_{m+1}=\mu$. Otherwise $L$ is finite
with longest path $\mu'$, in which case it suffices to let $l_{m+1}=d_{ \mathrm{l}( \mu')}$
and $r_{m+1}=\mu'$.
\end{proof}
\begin{lemma}
\label{lemma.5.2}\emph{\cite[Lemma 2.4.4]{Ben2018}} Suppose $M$ lies in $\mathcal{K}_{\mathrm{min}}(\Lambda\text{-}\boldsymbol{\mathrm{proj}})$. Let $i\in\mathbb{Z}$ and $m\in e_{v}M^{i}$. 
\begin{enumerate}
\item The $R$-module $L=e_{v}M^{i}$ is linearly compact.
\item If $U\subseteq_{\mathrm{c}}e_{v}M^{i}$ with $e_{v}\mathfrak{m}{}^{n}M^{i}\subseteq U$
for some $n>0$ then $U+m\subseteq_{\mathrm{c}}e_{v}M^{i}$.
\item For any $m\in e_{v}M^{i}$ we have $\{m\}=0+m\subseteq_{\mathrm{c}}e_{v}M^{i}$.
\end{enumerate}
\end{lemma}
\begin{proof}
(i) By Remark \ref{fingenrem} $L$ is finitely generated. By Example \ref{example.2.4.2} $L$ is linearly compact.

(ii), (iii) Let $l$ be a limit point of $U+m$. So for all $t\geq0$ there exists $u_{t}\in U$ with $u_{t}+m\neq l$
and $u_{t}+m\in l+\mathfrak{m}{}^{t}M$. So there exist $u_{n+1}\in U$ and $x_{n+1}\in\mathfrak{m}{}^{n+1}M$ with $u_{n+1}+m=l+x_{n+1}$.
Since $e_{v}(u_{n+1}+m-l)=u_{n+1}+m-l$ we have $x_{n+1}\in e_{v}\mathfrak{m}{}^{n}M\subseteq U$, and so $l=(u_{n+1}-x_{n+1})+m\in U+m$. This gives (ii). For (iii), any neighborhood of $l$ contains $m$, so $m-l\in e_{v}\mathfrak{m}{}^{t}M$ for all $t\geq0$. Hence $m-l\in\bigcap_{t\geq 0}e_{v}\mathfrak{m}{}^{t}M^{i}=0$
by the last statement in Remark \ref{rem2.6}. 
\end{proof}
\begin{remark}\label{closures}
For any $\gamma\in\mathbf{P}$ and any $\alpha\in\mathbf{A}$ we have that: if $U\subseteq_{\mathrm{c}} e_{t(\gamma)}M^{i}$ then $\gamma U\subseteq_{\mathrm{c}} e_{h(\gamma)}M^{i}$; if $V\subseteq_{\mathrm{c}}e_{h(\alpha)}M^{i}$ then $d_{\alpha}V\subseteq_{\mathrm{c}} e_{h(\alpha)}M^{i+1}$ and $d_{\alpha}^{-1}V\cap e_{h(\alpha)}M^{i-1}\subseteq_{\mathrm{c}} e_{h(\alpha)}M^{i-1}$;
and if $W\subseteq_{\mathrm{c}} e_{h(\gamma)}M^{i}$ then  $\gamma^{-1}W\cap e_{t(\gamma)}M^{i}\subseteq_{\mathrm{c}} e_{t(\gamma)}M^{i}$ (see \cite[Corollary 2.4.5]{Ben2018} for details). 
\end{remark}
\begin{corollary}
\label{cor25}\emph{\cite[Corollary 2.4.6]{Ben2018}} Suppose $M$ lies in $\mathcal{K}_{\mathrm{min}}(\Lambda\text{-}\boldsymbol{\mathrm{proj}})$, $N$
is an $R$-submodule of $M$ and that $C$ is a homotopy $\{0,\dots,t\}$-word. If $M^{i+\mu_{C}(t)}\cap N\subseteq_{\mathrm{c}} M^{i+\mu_{C}(t)}$
then $e_{h(C)}M^{i}\cap CN\subseteq_{\mathrm{c}} e_{h(C)}M^{i}$.
\end{corollary}
\begin{proof}By Corollary \ref{corollary.3.1} we have $e_{h(\gamma)}M^{i}\cap C_{\leq n}N=e_{h(\gamma)}M^{i}\cap C(M^{i+\mu_{C}(n)}\cap N)$ for all $n\leq t$. By Remark \ref{closures} $(d_{ \mathrm{l}( \gamma)}^{-1}\gamma)^{\pm1}(M^{i\pm 1}\cap N)\subseteq _{\mathrm{c}}e_{h(\gamma)}M^{i}$. The claim follows by iteration. 
\end{proof}
\begin{lemma}\label{lemma.5.3-1}\emph{\cite[Lemma 2.4.7]{Ben2018}} Suppose $M$ lies in $\mathcal{K}_{\mathrm{min}}(\Lambda\text{-}\boldsymbol{\mathrm{proj}})$.  Let $\gamma\in\mathbf{P}$, $i\in\mathbb{Z}$, $v(+)=t(\gamma)$ and $v(-)=h(\gamma)$.  Let $C=l_{1}^{-1}r_{1}\dots$ be a homotopy $I$-word where $\emptyset\neq I\subseteq\mathbb{N}$. 

\begin{enumerate}
\item If $m\in e_{h(\gamma)}M^{i}$ then $(\gamma^{-1}d_{ \mathrm{l}( \gamma)})^{\pm1}m\cap e_{v(\pm)}M^{i\pm1}$ is a linear variety if it is non-empty.
\item If $0<t\in I$ and $m\in M^{i+\mu_{C}(t-1)}\cap\bigcap_{n\in I,\,n\geq t}(C_{>t-1})_{\leq n}M$
then $M^{i+\mu_{C}(t)}\cap r_{t}^{-1}l_{t}m$ meets the intersection $\bigcap_{n\in I,\,n\geq t+1}(C_{>t})_{\leq n}M$.
\item If $I=\mathbb{N}$ and $S\subseteq e_{v}M^{i}$
then $S\cap C^{+}(M)=\bigcap_{n\geq0}S\cap C_{\leq n}M$.
\end{enumerate}
\end{lemma}
\begin{proof}
(i) Let $P=\gamma^{-1}d_{ \mathrm{l}( \gamma)}m\cap e_{t(\gamma)}M^{i+1}$. If $x\in P$ then $P'+x\subseteq P$ where $P'=e_{t(\gamma)}M^{i+1}\cap\gamma^{-1}0$. If $y\in P$ then $y-x\in P'$
as $\gamma y=d_{ \mathrm{l}( \gamma),M}(m)=\gamma x$. By Lemma \ref{lemma.5.2}(iii), Remark \ref{closures} and Corollary \ref{cor25}, $P$ is closed. This gives the case $v(\pm)=v(+)$.  The proof for the case $v(\pm)=v(-)$ is similar.

(ii) For each $n\geq t$ we have $m\in(C_{>t-1})_{\leq n}M$ and so
there is some $u_{n}\in M^{i+\mu_{C}(t)}\cap r_{t}^{-1}l_{t}m$ for
which $u_{n}\in(C_{>t})_{\leq n}M$. Let $\Delta$ be the collection of $M^{i+\mu_{C}(t)}\cap(C_{>t})_{\leq n}M$ ($n\geq t$) together with $M^{i+\mu_{C}(t)}\cap r_{t}^{-1}l_{t}m$. Let $V_{n}=M^{i+\mu_{C}(t)}\cap r_{t}^{-1}l_{t}m\cap(C_{>t})_{\leq n}M$
for each $n\geq0$. 

Clearly $\Delta$ has finite intersections. Each member of $\Delta$ lies in $e_{v_{C}(t)}M^{i+\mu_{C}(t)}$, which is  linearly compact by Lemma \ref{lemma.5.2}(i). By (i) we have $M^{i+\mu_{C}(t)}\cap r_{t}^{-1}l_{t}m\subseteq_{\mathrm{c}}e_{v_{C}(t)}M^{i+\mu_{C}(t)}$, and so by Corollary \ref{cor25} the collection $\Delta$ consists of linear varieties. Hence $\bigcap_{n\in\mathbb{N}}V_{n}=0$.

(iii) The argument in the proof of \cite[Lemma 3.1]{BenCra2018} adapts with few complications.
\end{proof}

\begin{proof}[of Lemma \ref{lemma.5.2-1}]
(i) Let $S=H\cap(U+m)$. Note $S\cap\rad(M)=\emptyset$ since $e_{v}\rad(M^{r})\subseteq U$ and $m\notin U$. So by Lemma \ref{lemma.5.3} there is a homotopy word $C$
such that either $C$ is finite and $S\cap C^{+}(M)\neq\emptyset=S\cap C^{-}(M)$,
or $C$ is a homotopy $\mathbb{N}$-word and for all $n\geq0$ we
have $S\cap C_{\leq n}M\neq\emptyset=S\cap C_{\leq n}\rad(M)$.
We assume $I_{C}=\mathbb{N}$ as otherwise there is nothing to prove. Note that the collection $(S\cap C_{\leq n}M\mid n\geq0)$ has the finite intersection property. By  Lemma \ref{lemma.5.2}(ii) and Corollary \ref{cor25} each $S\cap C_{\leq n}M$ is a linear variety in $e_{v}M^{r}$. By Lemma \ref{lemma.5.2}(i) $e_{v}M^{r}$ is linearly compact, and so $\bigcap_{n\geq0}S\cap C_{\leq n}M\neq\emptyset$.
This shows $S\cap C^{+}(M)\neq\emptyset$ by   Lemma \ref{lemma.5.3-1}(iii).
Since $S\cap C_{\leq n}\rad(M)=\emptyset$ for all $n$, $S\cap C^{-}(M)\subseteq\bigcup S\cap C_{\leq n}\rad(M)=\emptyset.$

(ii) The argument in the proof of \cite[Lemma 10.5]{Cra2018} adapts with few complications.
\end{proof}

\section{Local and global mapping properties.}\label{Local Mapping Properties.}
 
\begin{assumption}
In \S\ref{Local Mapping Properties.} fix an
object $M$ of $\mathcal{K}_{\mathrm{min}}(\Lambda\text{-}\boldsymbol{\mathrm{Proj}})$.
\end{assumption}
Recall: $\Sigma$ is the set of triples $(B,D,n)$ where $B^{-1}D$ is a homotopy word and $n\in\mathbb{Z}$; $\Sigma(s)$ (respectively $\Sigma(b)$) is the set of $(B,D,n)\in\Sigma$ where $B^{-1}D$ is aperiodic  (respectively periodic); $\mathcal{I}(s)\subseteq\Sigma(s)$ (respectively $\mathcal{I}(b)\subseteq\Sigma(b)$) is a collection of representatives $(B,D,n)$, one for each
equivalence class of $\Sigma(s)$ (respectively $\Sigma(b)$); and $\mathcal{I}=\mathcal{I}(s)\sqcup\mathcal{I}(b)$ (see Definition \ref{9.444}).
\begin{assumption}\label{ass12.2}
In Lemmas \ref{lemma.7.4} and \ref{lemma.7.4-1} fix $(B,D,n)\in\Sigma$ and let $C=B^{-1}D$. If $j\in I_{C}$, $i=a_{B,D}$ and $t=i-j$ then $v_{C}(j)=v_{B}(t)$ and $\mu_{C}(j)=\mu_{C}(i)+\mu_{B}(t)$ if $t\geq0$, and $v_{C}(j)=v_{D}(-t)$ and $\mu_{C}(j)=\mu_{C}(i)+\mu_{D}(-t)$ if $t<0$ \cite[Lemma 2.5.1]{Ben2018}. If $I_{B}\neq\{0\}$ let $B=l_{1}^{-1}r_{1}\dots$ and if $I_{D}\neq\{0\}$ let $D={l'}_{1}^{-1}{r'}_{1}\dots$ 
\end{assumption}
\begin{lemma}
\label{lemma.7.4}\emph{\cite[Lemma 2.5.2]{Ben2018} (}see also \emph{\cite[Lemma 8.3]{Cra2018})}. If $(B,D,n)\in\mathcal{I}(s)$ and $\mathcal{B}=(\bar{u}_{\lambda}\mid\lambda\in\Omega)$
is a $k$-basis of $F_{B,D,n}(M)$ then there is a morphism of complexes $\theta_{B,D,n,M}:\bigoplus_{\lambda}P(C)[\mu_{C}(a_{B,D})-n]\rightarrow M$
such that $F_{B,D,n}(\theta_{B,D,n,M})$ is an isomorphism. 
\end{lemma}
\begin{proof}
Let $i=a_{B,D}$. For each $\lambda$ choose a lift $u_{\lambda}\in F_{B,D,n}^{+}(M)\setminus F_{B,D,n}^{-}(M)$
of $\bar{u}_{\lambda}$. Since $u_{\lambda}\in B^{+}(M)$ we have, by Corollary \ref{corollary.3.1}, that for all $s\in I_{B}$ there exists $u''_{s,\lambda}\in e_{v_{B}(s)}M$ in degree $n+\mu_{B}(s)$ where $u''_{0,\lambda}=u_{\lambda}$ and
$u''_{s,\lambda}\in l_{s+1}^{-1}r_{s+1}u''_{s+1,\lambda}$ whenever
 $s+1\in I_{B}$. Similarly there are $u_{t,\lambda}'\in e_{v_{D}(t)}M$ in degree $n+\mu_{D}(t)$ where $u_{0,\lambda}'=u_{\lambda}$ and $u_{t,\lambda}'\in {l'}_{t+1}^{-1}{r'}_{t+1}u_{t+1,\lambda}'$ whenever $t+1\in I_{D}$. Let $u_{j,\lambda}=u''_{i-j,\lambda}$ if
$j\leq i$ and $u_{j,\lambda}=u_{j-i,\lambda}'$ if $j\geq i$.

By \cite[Lemma 2.5.1]{Ben2018} (see Assumption \ref{ass12.2}) setting $\theta_{B,D,n,M}( b _{j,\lambda})=u_{j,\lambda}$ (for all $j$ and $\lambda$)
defines a degree $0$ graded $\Lambda$-module map into $M$.
 Note $d_{M}(u_{j,\lambda})=u_{j,\lambda}^{+1}+u_{j,\lambda}^{-1}$
where $u_{j,\lambda}^{\pm1}=\sum_{\sigma^{\pm1}}d_{\sigma^{\pm1},M}(u_{j,\lambda})$
and $\sigma^{\pm1}$ runs through the arrows with head $v_{C}(j)$
and sign $\pm1$. It is straightforward to show $u_{j,\lambda}^{\pm 1}=\theta_{B,D,n,M}( b _{j,\lambda}^{\pm})$ by separating the cases $j\in I$ and $j\notin I$.

Together this gives $\theta_{B,D,n,M}( b _{j,\lambda}^{+}+ b _{j,\lambda}^{-})=d_{M}(u_{j,\lambda})$ and so $\theta_{B,D,n,M}$ is a morphism of complexes. By  
Lemma \ref{lemma.31}(i) the elements $c _{i,\lambda}= b _{i,\lambda}+F_{B,D,n}^{-}(P(C)[\mu_{C}(i)-n])$ ($\lambda\in \mathcal{B}$) together define a basis of the
$k$-vector space $F_{B,D,n}(\bigoplus_{\mathcal{\lambda}}P(C)[\mu_{C}(i)-n])$.
Since $\theta_{B,D,n,M}( b _{i,\lambda})=u_{i,\lambda}$
we have that $F_{B,D,n}(\theta_{B,D,n,M})(c _{i,\lambda})=\bar{u}_{i,\lambda}=\bar{u}_{\lambda}$
so $F_{B,D,n}(\theta_{B,D,n,M})$ is an isomorphism. 
\end{proof}
Recall $E(n)$ consists of pairs $(m,m')$ with $m,m'\in e_{v}M^{n}$ and $m\in Em'$, and so $E(n)^{\sharp}=F_{B,D,n}^{+}(M)$
and $E(n)^{\flat}=F_{B,D,n}^{-}(M)$ (see Remark \ref{remreffun}). 
\begin{lemma}
\label{lemma.7.4-1}\emph{\cite[Lemma 2.5.4]{Ben2018} (}see also \emph{\cite[Lemma 8.6]{Cra2018})}. Suppose $(B,D,n)\in\mathcal{I}(b)$, say where $C={}^{\infty}E{}^{\,\infty}$ is periodic of period $p>0$ and $E=l_{1}^{-1}r_{1}\dots l_{p}^{-1}r_{p}$. If there is a reduction $(U,f)$ of the relation $E(n)$ on $e_{v}M^{n}$ then there is a morphism $\theta_{B,D,n,M}:P(C,U)[-n]\rightarrow M$
of complexes such that $F_{B,D,n}(\theta_{B,D,n,M})$ is an isomorphism.
\end{lemma}
\begin{proof}
Let $F_{B,D,n}(M)=V$. Choose an $R$-basis $(u_{\lambda}\mid \lambda\in\Omega)$
of $U$. Since $(U,f)$ is a reduction, $\im(f)\subseteq E(n)^{\sharp}$
and so  $f(u_{\lambda})\in F_{B,D,n}^{+}(M)$ for each $\lambda$.  Similarly $f(Tu_{\lambda})\in E f(u_{\lambda})$, and so there are elements $u_{0,\lambda},\dots u_{p,\lambda}\in M$
where $u_{p,\lambda}=g(u_{\lambda})$,
$u_{0,\lambda}=g(Tu_{\lambda})$, and $u_{j-1,\lambda}\in l_{j}^{-1}r_{j}u_{j,\lambda}$ given $j>0$. By Lemma \ref{lemma.2.3}, to define a $\Lambda$-module
map $\theta_{B,D,n,M}:P(C,U)[-n]\rightarrow M$ it is enough to extend
$\theta_{B,D,n,M}( b _{j}\otimes u_{\lambda})=u_{j,\lambda}$ linearly
over $\Lambda$ for each $\lambda$ and each $j$ with $0\leq j \leq p-1$. By Corollary \ref{corollary.3.1} we have that $u_{j,\lambda}\in e_{v_{E}(j)}M^{n+\mu_{E}(j)}$ for each $j$, and so $\theta_{B,D,n,M}$ is homogeneous
of degree $0$. As in the proofs of \cite[Lemma 8.6]{Cra2018} and \cite[\S 5, Proposition]{Rin1975} the morphism $P(C)\times U\to M$ given by $((b_{j},u_{\lambda})\mapsto u_{j,\lambda})$ is $R[T,T^{-1}]$-balanced. To show $\theta_{B,D,n,M}$ is a morphism
of complexes one uses Lemma \ref{lemma.1.1}(ii), separating the cases $j=0$, $j=p-1$ and $p-1\neq j\neq0$. 

We now show that $F_{B,D,n}(\theta_{B,D,n,M})$ is an isomorphism. Write $e _{p,\lambda}$
for the coset $  b _{p}\otimes u_{\lambda}+F_{B,D,n}^{-}(P(C,U)[-n])$. By Lemma \ref{lemma.7.3-1}(ii) the elements $(e _{p,\lambda}\mid \lambda\in\Omega)$
give a $k$-basis of $F_{B,D,n}(P(C,U)[-n])\simeq k\otimes_{R}V$. Let $u_{p,\lambda}+F_{B,D,n}^{-}(P(C,U)[-n])=\bar{u}_{p,\lambda}$ for each $\lambda\in\Omega$. 

Since $F_{B,D,n}(\theta_{B,D,n,M})(e _{p,\lambda})=\bar{u}_{p,\lambda}$,
to prove $F_{B,D,n}(\theta_{B,D,n,M})$ is an isomorphism we need
only show $(\bar{u}_{p,\lambda}\mid\lambda\in\Omega)$ is a $k$-linearly
independent subset of $V=F_{B,D,n}(M)$. If we have $\sum_{\lambda}(r_{\lambda}+\mathfrak{m})\bar{u}_{p,\lambda}=0$
in $V$ for some $(r_{\lambda}\mid\lambda\in\Omega)\in \bigoplus _{\lambda} R$ then $f( 
\sum_{\lambda}r_{\lambda}u_{\lambda})\in E(n)^{\flat}$. Since the reduction $(U,f)$ meets in $\mathfrak{m}$, we have $\sum_{\lambda}r_{\lambda}u_{\lambda}\in\bigoplus_{\lambda}\mathfrak{m}u_{l}$, as required.
\end{proof}
\begin{lemma}
\label{lemma.7.5}\emph{\cite[Lemma 2.5.5]{Ben2018} (}see also \emph{\cite[Lemma 10.5]{Cra2018}} and \emph{\cite[p. 163]{ButRin1987})}. Let $\theta:P\rightarrow M$ be a morphism in \emph{$\mathcal{K}_{\mathrm{min}}(\Lambda\text{-}\boldsymbol{\mathrm{Proj}})$},
and suppose $M^{i}$ is finitely generated for each $i$. If $F_{B,D,n}(\theta)$ is surjective for each $(B,D,n)\in\Sigma$ then $\theta^{i}$
is surjective for each $i$.
\end{lemma}
\begin{proof}
For a contradiction suppose that $\theta^{i}$ is not surjective for
some $i\in\mathbb{Z}$. Since $\mbox{rad}(M^{i})$ is a superfluous
submodule of $M^{i}$, $e_{v}\mbox{im}(\theta^{i})+e_{v}\mbox{rad}(M^{i})$
is contained in a maximal $R$-submodule $U$ of $e_{v}M^{i}$. Since
$e_{v}\mbox{rad}(M^{i})\subseteq U$ and $U\neq e_{v}M^{i}$, by Lemma \ref{lemma.5.2-1}(ii)
for some element $m\in e_{v}M^{i}\setminus U$ there are homotopy words $B\in\mathcal{W}_{v,\delta}$
and $D\in\mathcal{W}_{v,-\delta}$ for which ($B^{-1}D$ is a homotopy word)
and $U+m$ meets $G_{B,D,i}^{+}(M)$ but not $G_{B,D,i}^{-}(M)$. From here one can show $F_{B,D,i}(\theta)$ is not surjective by adapting the argument from the proof of \cite[Lemma 10.6]{Cra2018}.
\end{proof}
\begin{assumption}\label{ass12.5}
For the remainder of \S\ref{Local Mapping Properties.} we fix  a direct sum $N$ of shifts of string and band complexes as follows.
Let $\mathcal{S}$ and $\mathcal{B}$ be sets, $\{t(\sigma),s(\beta)\mid\sigma\in\mathcal{S},\beta\in\mathcal{B}\}$
be a set of integers, $\{V^{\beta}\mid\beta\in\mathcal{B}\}$
be a set of objects from $R[T,T^{-1}]\text{-\textbf{Mod}}_{R\text{-\textbf{Proj}}}$
and $\{A(\sigma),E(\beta)\mid\sigma\in\mathcal{S},\beta\in\mathcal{B}\}$
be a set of homotopy words, where each $A(\sigma)$ is aperiodic
and each $E(\beta)$ is $p_{\beta}$-periodic. Let 
\[
N=\Bigl( \bigoplus_{\sigma\in\mathcal{S}}P(A(\sigma))[-t(\sigma)]\Bigr)\oplus\Bigl( \bigoplus_{\beta\in\mathcal{B}}P(E(\beta),V^{\beta})[-s(\beta)]\Bigr)
\]
\end{assumption} 
\begin{definition}\label{nota}The homotopy words $(A(\sigma)_{\leq 0})^{-1}$ and $A(\sigma)_{>0}$ (respectively $(E(\beta)_{\leq 0})^{-1}$ and $E(\beta)_{>0}$) have the same head and opposite sign, and we let $A(\sigma,\pm)$ (respectively $E(\beta,\pm)$) be the one with sign $\pm1$.
\end{definition}
\begin{lemma}\emph{\cite[Lemma 2.5.6]{Ben2018} (}see also \emph{\cite[Lemma 9.4]{Cra2018})}.
\label{lemma.7.1-1} Let $N$ be the direct sum of string and band complexes from Assumption \ref{ass12.5}. Let $\theta:N\rightarrow M$ be a map in $\mathcal{K}_{\mathrm{min}}(\Lambda\text{-}\boldsymbol{\mathrm{Proj}})$
where $\bar{F}_{B,D,n}(\theta)$ is injective for all $(B,D,n)\in\mathcal{I}$.
Then each $\theta^{i}$ is injective.
\end{lemma}
\begin{proof}
Assume there is some $h\in\mathbb{Z}$ for which $\theta^{h}$ is not
injective. Since $N^{h}$ and $M^{h}$ are projective $\Lambda$-modules we have $\mathrm{rad}(N^{h})=\mathrm{rad}(\Lambda)N^{h}$ and $\mathrm{rad}(M^{h})=\mathrm{rad}(\Lambda)M^{h}$ (see for example \cite[Theorem 24.7]{Lam1991}). Since $\Lambda/\mathrm{rad}(\Lambda)$
is semisimple, $N^{h}/\rad(N^{h})$ is an injective $\Lambda/\rad(\Lambda)$-module.
Hence the induced map $\bar{\theta}^{h}:N^{h}/\mathrm{rad}(N^{h})\rightarrow M^{h}/\mathrm{rad}(M^{h})$
is not injective, as otherwise it must be a section, which would
mean $\theta^{h}$ was injective by \cite[Lemma 2.2]{Jon1976}.

So there is a vertex
$v$ and non-zero $n\in e_{v}N^{h}\setminus e_{v}\rad(N^{h})$ with $\theta^{h}(n)\in e_{v}\rad(M^{h})$. Since $n$ has finite support over the summands of $N$, we can assume $\mathcal{S}=\{1,\dots,m\}$ and $\mathcal{B}=\{1,\dots,q\}$.  For each $\sigma\in\mathcal{S}$, and each $i\in I_{A(\sigma)}$, let $b_{i}^{\sigma}=b_{i,A(\sigma)}$. For each $\beta\in\mathcal{B}$, fix an $R$-basis $(v^{\beta}_{\lambda}\mid\lambda\in\Omega(\beta))$ of $V^{\beta}$, and let $b_{j,\lambda}^{\beta}=b_{j,E(\beta)}\otimes v_{\lambda}^{\beta}$. Fix a transversal $\mathtt{S}$ of $\mathfrak{m}$ in $R$. Altogether 
\[
n=x+\sum_{\sigma=1}^{m}\sum_{i\in I_{A(\sigma)} }\eta_{\sigma,i} b _{i}^{\sigma}+\sum_{\beta=1}^{q}\sum_{j\in I_{E(\beta)}}\sum_{\lambda\in\Omega(\beta)}\eta_{\beta,j,\lambda} b _{j,\lambda}^{\beta}
\]
for some $x\in e_{v}\rad(N^{h})$ and some $\eta_{\sigma,i},\eta_{\beta,j,\lambda}\in \mathtt{S}$. After reordering
we can assume $A(\sigma,1)\leq A(\sigma',1)$ for $\sigma\leq\sigma'$
and $E(\beta,1)\leq E(\beta',1)$ for $\beta\leq\beta'$. Let $B$
be the largest of $A(m,1)$ and $E(q,1)$. If $B=A(m,1)$
let $D$ be the largest homotopy word $A(\sigma',-1)$ among $A(1,-1),\dots,A(m,-1)$
for which $A(\sigma',1)=A(m,1)$. In this case there is some such $i$ with $v_{A(\sigma)}(i)=v$ and $\mu_{A(\sigma)}(i)=h-t(\sigma)$: and so $n_{\sigma,i}b^{\sigma}_{i}-x\in G_{B,D,h}^{+}(N)\setminus G_{B,D,h}^{-}(N)$. If instead $B=E(q,1)$ we let $D=E(\beta',-1)$,
the largest among $E(1,-1),\dots,E(q,-1)$
for which $E(\beta',1)=E(q,1)$. As above $n_{\beta,j,\lambda}b^{\beta}_{j,\lambda}-x\in G_{B,D,h}^{+}(N)\setminus G_{B,D,h}^{-}(N)$ for some $j$. As in the proof of \cite[Lemma 9.4]{Cra2018}, in either case this shows $\bar{G}_{B,D,h}(\theta)$ is not injective.
\end{proof}
\section{\label{sec:Proofs-of-the}Completing the proofs of the main theorems.}
We now verify the hypotheses of Lemma \ref{lemma} in our setting. 
Recall, from Definition \ref{definition.7.14}, the equivalence relation on the triples $(B,D,n)$ where $B^{-1}D$ is a homotopy word and $n\in\mathbb{Z}$. Recall, from Definition \ref{9.444}, that $\mathcal{I}=\mathcal{I}(s)\sqcup \mathcal{I}(b)$ where $\mathcal{I}(s)$ (respectively $\mathcal{I}(b)$) denotes a chosen set of representatives $(B,D,n)$ such that $B^{-1}D$ is aperiodic (respectively periodic). 

Recall that if $(B,D,n)$ lies in $\mathcal{I}(s)$ (respectively $\mathcal{I}(b)$) then the functor $S_{B,D,n}$ has the form $R\text{-}\mathbf{Proj}\rightarrow \mathcal{C}_{\mathrm{min}}(\Lambda\text{-}\mathbf{Proj})$ (respectively $R[T,T^{-1}]\text{-}\mathbf{Mod}_{R\text{-}\mathbf{Proj}}\rightarrow\mathcal{C}_{\mathrm{min}}(\Lambda\text{-}\mathbf{Proj})$), and the functor $F_{B,D,n}$ has the form $\mathcal{K}_{\mathrm{min}}(\Lambda\text{-}\mathbf{Proj})\rightarrow k\text{-}\mathbf{Mod}$ (respectively $\mathcal{K}_{\mathrm{min}}(\Lambda\text{-}\mathbf{Proj})\rightarrow k[T,T^{-1}]\text{-}\mathbf{Mod}$). Recall Definition \ref{detect}.
\begin{proposition}
\label{proposition6}\emph{\cite[Proposition 2.6.6]{Ben2018} (}see also  \emph{\cite[p. 163, Proposition]{ButRin1987})}. Let $\mathcal{M}=\Lambda\text{-}\boldsymbol{\mathrm{Mod}}$ and $\mathfrak{I}=\mathcal{I}(s)\sqcup\mathcal{I}(b)$;
and for each $i=(B,D,n)\in\mathfrak{I}$ let
\[
(\mathfrak{A}_{i},\mathfrak{X}_{i})=\left\{ \begin{array}{cc}
(R\text{-}\boldsymbol{\mathrm{Proj}},k\text{-}\boldsymbol{\mathrm{Mod}}) & (\text{if }B^{-1}D\text{ is aperiodic}),\\
(R[T,T^{-1}]\text{-}\boldsymbol{\mathrm{Mod}}_{R\text{-}\boldsymbol{\mathrm{Proj}}},k[T,T^{-1}]\text{-}\boldsymbol{\mathrm{Mod}}) & (\text{if }B^{-1}D\text{ is periodic}).
\end{array}\right.
\]
\begin{enumerate}
    \item The category $\mathcal{N}=\Lambda\text{-}\boldsymbol{\mathrm{mod}}$  has enough projective covers.
    \item The collection $\{(S_{B,D,n},F_{B,D,n})\mid(B,D,n)\in\mathcal{I}\}$
detects the objects in $\mathcal{K}_{\mathrm{min}}(\Lambda\text{-}\boldsymbol{\mathrm{proj}})$.
\end{enumerate} 
\end{proposition}
In the proof of Proposition \ref{proposition6}(ii) we verify the conditions FFI, FFII, FFIII, FFIV, FFV and FFVI from Definition \ref{detect}. Later we use Proposition \ref{proposition6} in the context of Lemma \ref{lemma}.
\begin{proof}(i) By Corollary \ref{corollary.0.1}(vi) $\Lambda$ is a semiperfect ring, and so every finitely generated $\Lambda$-module has a projective cover (see for example \cite[Proposition 24.12]{Lam1991}).

(ii) FFI) Recall $k[T,T^{-1}]\otimes_{R[T,T^{-1}]}-\simeq k\otimes_{R}-$ as functors $R[T,T^{-1}]\text{-}\boldsymbol{\mathrm{Mod}}_{R\text{-}\boldsymbol{\mathrm{Proj}}}\rightarrow k[T,T^{-1}]\text{-}\boldsymbol{\mathrm{Mod}}$. Let $(B,D,n)\in\Sigma$ and $B^{-1}D=C$. If $C$ is aperiodic (respectively periodic) then by Lemma \ref{lemma.31}(ii) (respectively  Lemma \ref{lemma.7.3-1}(ii)) we have $F_{B,D,n}\,\Xi\, S_{B,D,n}\simeq k\otimes_{R}-$ as functors $\mathfrak{A}_{B,D,n}\to\mathfrak{X}_{B,D,n}$. We show $k\otimes_{R}-:\mathfrak{A}_{B,D,n}\to\mathfrak{X}_{B,D,n}$ is dense and reflects isomorphisms. 

Recall that any local ring (such as $R$) is a semiperfect ring. If $g$ is a morphism in $\mathfrak{A}_{B,D,n}$ then $g$ is homomorphism of free $R$-modules, and $g$ is an isomorphism in $\mathfrak{A}_{B,D,n}$ if and only if $g$ is bijective. Similarly any isomorphism in $\mathfrak{X}_{B,D,n}$ is an isomorphism of vector spaces. Altogether, by Remark \ref{projiso}, if $g$ is a morphism in $\mathfrak{A}_{B,D,n}$ such that $k\otimes g$ is an isomorphism in $\mathfrak{X}_{B,D,n}$, then $g$ was an isomorphism. Hence $k\otimes_{R}-$ reflects isomorphisms.

Note that for any vector space $V$ with basis $(v_{\lambda}\mid\lambda\in\Omega)$ over $k$ there is an isomorphism $f:V\to k\otimes_{R}F$ of $k$-vector spaces where $F=\bigoplus _{\lambda} Ra_{\lambda}$ is the free $R$-module with $R$-basis $(a_{\lambda}\mid\lambda\in\Omega)$, such that $f(v_{\lambda})=(1+\mathfrak{m})\otimes a_{\lambda}$. This shows $k\otimes_{R}-$ is dense when $(B,D,n)$ lies in $\mathcal{I}(s)$. Suppose instead that $(B,D,n)$ lies in $\mathcal{I}(b)$, and that $V$ is a $k[T,T^{-1}]$-module. There is an $R$-module endomorphism of $F$ given by $a_{\lambda}\mapsto\sum_{\mu}r_{\mu\lambda}a_{\mu}$ for each $\lambda$ where $Tv_{\lambda}=\sum_{\mu}(r_{\mu\lambda}+\mathfrak{m})v_{\mu}$ for some $r_{\mu\lambda}\in R$. As above, by Remark \ref{projiso} this endomorphism is an isomorphism, and so $F$ is an $R[T,T^{-1}]$-module where $Ta_{\lambda}=\sum_{\mu}r_{\mu\lambda}a_{\mu}$. By construction $f$ is $k[T,T^{-1}]$-linear.

FFII) Let $(B',D',n')\in\mathcal{I}$. If $(B',D',n')\neq (B,D,n)\in\mathcal{I}(s)$ then $\bar{F}_{B',D',n'}(P(C)[\mu_{C}(a_{B,D})-n])=0$ by Lemma \ref{lemma.31}(iii) where $C=B^{-1}D$. This shows $F_{B',D',n'}\Xi S_{B,D,n}=0$ since $\bar{F}_{B',D',n'}\simeq F_{B',D',n'}$.  If $(B,D,n)\in\mathcal{I}(b)$ then the proof is similar and uses Lemma \ref{lemma.7.3-1}(iii).

FFIII)  By Corollary \ref{corollary.3.4} each of the subfunctors $C^{\pm}$ of the forgetful functor $\mathcal{C}_{\mathrm{min}}(\Lambda\text{-}\boldsymbol{\mathrm{Proj}})\rightarrow R\text{-}\boldsymbol{\mathrm{Mod}}$ commutes with arbitrary direct sums. It follows that $F^{\pm}_{B,D,n}$ commutes with direct sums of objects in $\mathcal{K}_{\mathrm{min}}(\Lambda\text{-}\boldsymbol{\mathrm{Proj}})$ (see \cite[Lemma 2.1.21]{Ben2018} for details).

FFIV) Let $M$ be an object in $\mathcal{K}_{\mathrm{min}}(\Lambda\text{-}\boldsymbol{\mathrm{proj}})$, which means $F_{B,D,n}(M)$ is  finite-dimensional as a $k$-vector space. If $(B,D,n)$ lies in $\mathcal{I}(s)$, by Lemma \ref{lemma.7.4} there is a free $R$-module $U$ of $R$-rank $\mathrm{dim}_{k}(F_{B,D,n}(M))$ and a morphism $\theta_{B,D,n}:S_{B,D,n}(U)\rightarrow M$ for which $F_{B,D,n}(\theta_{B,D,n})$ is an isomorphism. Now suppose instead $(B,D,n)$ lies in $\mathcal{I}(b)$, say where $B^{-1}D={}^{\infty}E{}^{\,\infty}$ is periodic of period $p>0$, and where $E=l_{1}^{-1}r_{1}\dots l_{p}^{-1}r_{p}$. Note $F_{B,D,n}^{+}(M)=E(n)^{\sharp}$ and $F_{B,D,n}^{-}(M)=E(n)^{\flat}$, so $E(n)^{\sharp}/E(n)^{\flat}=F_{B,D,n}(M)$. By Lemma \ref{lemma.2.11} there is a reduction $(U,f)$ of
$E(n)$ which meets in $\mathfrak{m}$. The required morphism exists by Lemma \ref{lemma.7.4-1}.

FFV, FFVI) Let $\theta:N\rightarrow M$ be an morphism in the category $\mathcal{K}_{\mathrm{min}}(\Lambda\text{-}\boldsymbol{\mathrm{Proj}})$. If $M$ is a complex in $\mathcal{K}_{\mathrm{min}}(\Lambda\text{-}\boldsymbol{\mathrm{proj}})$ and $F_{B,D,n}(\theta)$ is epic for all $(B,D,n)\in\mathcal{I}$ then $\theta^{n}$ is epic for each $n\in\mathbb{N}$ by Lemma \ref{lemma.7.5}. This shows FFV holds, and similarly FFVI holds by Lemma \ref{lemma.7.1-1}.
\end{proof}
\begin{proof}[of Theorem \ref{theorem.1.1}]
Parts (i) and (ii) of Theorem \ref{theorem.1.1} are precisely
parts (i) and (iii) of  Lemma \ref{lemma}, after applying Proposition
\ref{proposition6}.\end{proof}
Note that, in the context of the proof of Theorem \ref{theorem.1.1}, Lemma \ref{lemma}(ii) says that any indecomposable object in $\mathcal{K}_{\mathrm{min}}(\Lambda\text{-}\boldsymbol{\mathrm{proj}})$ is the shift of a string or band complex.
\begin{assumption}\label{finalass}
In what remains we recall and fix notation from Assumption \ref{ass12.5} and Definition \ref{nota}. Let $\mathcal{S}$, $\mathcal{S}'$, $\mathcal{B}$ and $\mathcal{B}'$ be index sets. For each $\sigma\in\mathcal{S}$ and each $\sigma'\in\mathcal{S}'$ let $t(\sigma),w(\sigma')\in\mathbb{Z}$ and let $A(\sigma)$ and $B(\sigma')$ be aperiodic homotopy words. For each $\beta\in\mathcal{B}$ and each $\beta'\in\mathcal{B}'$ let $s(\beta),r(\beta')\in\mathbb{Z}$, let $E(\beta)$ and $D(\beta')$ be periodic homotopy words and let $V^{\beta}$ and $W^{\beta'}$ be indecomposable objects of $R[T,T^{-1}]\text{-}\boldsymbol{\mathrm{Mod}}_{R\text{-}\mathbf{Proj}}$. Define the objects $N$ and $N'$ in $\mathcal{K}_{\mathrm{min}}(\Lambda\text{-}\mathbf{Proj})$ by
\[\begin{array}{c}
N = \Bigl( \bigoplus_{\sigma\in\mathcal{S}}P(A(\sigma))[-t(\sigma)]\Bigr)\oplus\Bigl( \bigoplus_{\beta\in\mathcal{B}}P(E(\beta),V^{\beta})[-s(\beta)]\Bigr),\\
N'=\Bigl( \bigoplus_{\sigma'\in\mathcal{S}'}P(B(\sigma'))[-w(\sigma')]\Bigr)\oplus\Bigl( \bigoplus_{\beta'\in\mathcal{B}'}P(D(\beta'),W^{\beta'})[-r(\beta')]\Bigr).
\end{array}
\]
Let $(B,D,n)\in\Sigma$.  If $(B,D,n)\in\mathcal{I}(s)$ let $\mathcal{S}(B,D,n)$ be the set of $\sigma$ such that $t(\sigma)-n=r(B,D;C(1),C(-1))$ and $A(\sigma)\sim B^{-1}D$. If $(B,D,n)\in\mathcal{I}(b)$ let $\mathcal{B}^{\pm}(B,D,n)$ be the set of $\beta$ such
that $s(\beta)-n=\mu_{B^{-1}D}(\pm m)$ and $E(\beta)=(B^{-1}D)^{\pm1}[m]$. For $\beta\in\mathcal{B}(B,D,n)^{\pm}$ and $\beta'\in\mathcal{B}'(B,D,n)^{\pm}$ define $V_{\pm}(\beta)$, $\bar{V}_{\pm}(\beta)=k\otimes_{R}V_{\pm}^{\beta}$, $W_{\pm}(\beta')$ and $\bar{W}_{\pm}(\beta')=k\otimes_{R}W_{\pm}^{\beta}$ by
\[
\begin{array}{cccc}
V_{+}(\beta)=V^{\beta}, & V_{-}(\beta)=\mathrm{res}_{\iota,R}\,V^{\beta}, & W_{+}(\beta')=W^{\beta'}, & W_{-}(\beta')=\mathrm{res}_{\iota,R} \,W^{\beta'}.
\end{array}
\] 
\end{assumption}
\begin{lemma}
\label{lemma.7.1}\emph{\cite[Lemma 2.5.8]{Ben2018}} Let $B$, $B'$, $D$ and $D'$ be homotopy words such that $C=B^{\,-1}D$ and $C'=B'^{\,-1}D'$ are homotopy words. Let $n$ and $n'$ be integers. 
\begin{enumerate}
\item If $n\neq n'$ then $\mathcal{S}(B,D,n)\cap\mathcal{S}(B,D,n')=\emptyset=\mathcal{B}(B,D,n)^{\pm}\cap\mathcal{B}(B,D,n')^{\pm}$.
\item If $B^{-1}D\nsim B'^{-1}D'$ then $\mathcal{S}(B,D,n)\cap\mathcal{S}(B',D',n)=\emptyset=\mathcal{B}(B,D,n)^{\pm}\cap\mathcal{B}(B',D',n)^{\pm}$.
\item We have $\bigcup_{t\in\mathbb{Z}}\mathcal{S}(B,D,t)=\{\sigma\in\mathcal{S}\mid A(\sigma)\sim C\}$.
\item We have $\bigcup_{t\in\mathbb{Z}}\mathcal{B}(B,D,t)^{+}\cup\mathcal{B}(B,D,t)^{-}=\{\beta\in\mathcal{B}\mid E(\beta)\sim C\}$.
\end{enumerate}
\end{lemma}
\begin{theorem} \emph{\cite[Theorem 2.5.9]{Ben2018} (}see also \emph{\cite[Theorem 9.1]{Cra2018})}.
\label{theorem.7.1}Let $(B,D,m)\in\Sigma$. 
\begin{enumerate}
\item If $C=B^{-1}D$ is aperiodic then in $k\text{-}\mathbf{Mod}$ we have isomorphisms 
\[
F_{B,D,m}(N)\simeq\bigoplus_{\sigma\in \mathcal{S}(B,D,m)}k\text{, and }\bigoplus_{n\in\mathbb{Z}}F_{B,D,n}(N)\simeq \bigoplus_{\sigma\in \mathcal{S}\, \mid\,A(\sigma)=C^{\pm 1}}k.
\]
\item If $C$ is periodic then in $k[T,T^{-1}]\text{-}\mathbf{Mod}$ we have isomorphisms
\[
F_{B,D,m}(N)\simeq\bigoplus_{\beta\in \mathcal{B}(B,D,m)^{+}}\bar{V}_{+}(\beta)\oplus\bigoplus_{\beta'\in \mathcal{B}(B,D,m)^{-}}\bar{V}_{-}(\beta')\text{, and }
\]
\[
\bigoplus_{n\in\mathbb{Z}}F_{B,D,n}(N)\simeq\bigoplus_{\beta\in \mathcal{B}\,\mid\,E(\beta)= C[t]}\bar{V}_{+}(\beta)\oplus \bigoplus_{\beta'\in \mathcal{B}\,\mid\,E(\beta')= C[t]^{-1}}\bar{V}_{-}(\beta').
\]
\end{enumerate}
\end{theorem}
\begin{proof}
(i) By Corollary \ref{corollary.3.4} the functor $F_{B,D,m}$ preserves small coproducts. This together with 
Lemma \ref{lemma.31}(iii) shows $F_{B,D,m}(N)\simeq\bigoplus_{\sigma\in\mathcal{S}}F_{B,D,m}(P(A(\sigma))[-t(\sigma)])$
as $C$ is aperiodic. If $\sigma\in\mathcal{S}(B,D,m)$
then $(A(\sigma,1),A(\sigma,-1),t(\sigma))\sim(B,D,m)$
and so by Lemma \ref{corollary.6.6} and Lemma \ref{lemma.31}(ii)
we have $F_{B,D,m}(P(A(\sigma))[-t(\sigma)])\simeq\bar{F}_{B,D,m}(S_{B,D,m}(R))\simeq k$. Otherwise $\sigma\notin\mathcal{S}(B,D,m)$ and so as above
$F_{B,D,m}(P(A(\sigma))[-t(\sigma)])=0$. By  Lemma \ref{lemma.7.1}(i) $\sum_{n\in\mathbb{Z}}\#\,\mathcal{S}(B,D,n)=\#\,\bigcup_{n\in\mathbb{Z}}\mathcal{S}(B,D,n)$. By Lemma \ref{lemma.7.1}(iii) this completes the proof of (i).

(ii) Similar to the above, $F_{B,D,m}(N)$  is isomorphic to $\bigoplus_{\beta\in\mathcal{B}}F_{B,D,m}(P(E(\beta),V^{\beta})[-s(\beta)])$ by Lemma  \ref{lemma.7.3-1}(iii).
If $\beta\in\mathcal{B}(B,D,m)^{\pm}$ then by Lemma \ref{corollary.6.6} and Lemma \ref{lemma.7.3-1}(ii) we have that $F_{B,D,m}(P(E(\beta),V^{\beta})[-s(\beta)])$ is isomorphic to $\bar{V}_{\pm}(\beta)$
as above. If $\beta\notin\mathcal{B}(B,D,m)^{+}\cup\mathcal{B}(B,D,m)^{-}$
then $F_{B,D,m}(P(E(\beta),V^{\beta})[-s(\beta)])=0$ by 
Lemma \ref{lemma.7.3-1}(iii). By Lemma \ref{lemma.7.1}(ii)  this shows $F_{B,D,m}(N)\simeq(\bigoplus_{\beta_{+}}\bar{V}_{+}(\beta_{+}))\oplus(\bigoplus_{\beta_{-}}\bar{V}_{-}(\beta_{-}))$
where $\beta_{\pm}$ runs through $\mathcal{B}(B,D,m)^{\pm}$. By Lemma \ref{lemma.7.1}(iv) this completes the proof of (ii).
\end{proof}
\begin{proof}[of Theorem \ref{theorem.1.2}]
Recall Definition \ref{nota}. Let $C$ and $E$ be homotopy words and let $n\in\mathbb{Z}$. We assume $s((E_{\leq 0})^{-1})=1$ and $s(E_{>0})=-1$. The case where  $s((E_{\leq 0})^{-1})=-1$ and $s(E_{>0})=1$ is similar, and omitted. Recall that homotopy words $(C_{\leq 0})^{-1}$ and $C_{>0}$ have the same head and opposite sign. We let $C(\pm 1)$ be the one with sign $\pm1$. Choose $\delta\in\{1,-1\}$ such that  
$C(-\delta)=(C_{\leq 0})^{-1}$ and $C(\delta)=C_{>0}$. Let $V$ and $W$ be objects of  $R[T,T^{-1}]\text{-}\boldsymbol{\mathrm{Mod}}_{R\text{-}\boldsymbol{\mathrm{Proj}}}$.

(i) Let $C$ and $E$ be aperiodic. Assume momentarilly that: $I_{C}=\{0,\dots,m\}$ and ($(I_{E},E,n)=(I_{C},C,0)$ or $(I_{E},E,n)=(I_{C},C^{-1},\mu_{C}(m))$); or that $I_{C}=\pm\mathbb{N}$ and ($(I_{E},E,n)=(\pm\mathbb{N},C,0)$ or $(I_{E},E,n)=(\mp\mathbb{N},C^{-1},0)$); or that $I_{C}=\mathbb{Z}$ and $(I_{E},E,n)=(\mathbb{Z}, C^{\pm1}[t],\mu_{C}(\pm t))$ for some $t\in\mathbb{Z}$. In each of the cases above (respectively) there is an isomorphism $P(C)[n]\simeq P(E)$ in $\mathcal{K}(\Lambda\text{-}\boldsymbol{\mathrm{Proj}})$ given by parts (i), (ii) and (iii) of Lemma \ref{lemma.4.1}.

For the converse suppose $P(C)[n]\simeq P(E)$ in $\mathcal{K}(\Lambda\text{-}\boldsymbol{\mathrm{Proj}})$. In the notation of Assumption \ref{finalass} we have $N=P(C)[n]$ and $N'=P(E)$ where $\mathcal{B}=\mathcal{B}'=\emptyset$, $\mathcal{S}=\{\sigma\}$, $\mathcal{S}'=\{\sigma'\}$, $A(\sigma)=C$, $B(\sigma')=E$,  $t(\sigma)=-n$ and $w(\sigma')=0$. By Theorem \ref{theorem.7.1}(i) we have that, for any $(H,L,q)\in\Sigma$, $F_{H,L,q}(N)$ (respectively $F_{H,L,q}(N')$) is isomorphic to $\bigoplus k$ where the direct sum runs through all $\sigma\in \mathcal{S}$ (respectively $\sigma'\in\mathcal{S}'$) such that $n-q=r(H,L;C(1),C(-1))$  (respectively  $0-q=r(H,L;(E_{\leq 0})^{-1},E_{>0})$) and $C\sim H^{-1}L$ (respectively $E\sim H^{-1}L$). 

Applying $F_{H,L,0}$ to $N\simeq N'$ in case $(H,L)=((E_{\leq 0})^{-1},E_{>0})$ gives $F_{H,L,0}(N)\simeq k$, which shows that $n=r((E_{\leq 0})^{-1},E_{>0};C(1),C(-1))$ and $C\sim E$. By definition this means that: $C=E$ which is not a homotopy $\mathbb{Z}$-word, and so $E_{>0}=C$ and $s(C)=1$, and so $n=\mu_{E}(0)=0$; or $C=E^{-1}$ which is not a homotopy $\mathbb{Z}$-word, and so $E_{\leq m}=C^{-1}$ and $s(C)=-1$, and so $n=\mu_{E}(m)$ where $I_{C}=\{0,\dots,m\}$ or ($I_{C}=\pm \mathbb{N}$ and $m=0$); or $C=E^{\pm1}[t]$ is a homotopy $\mathbb{Z}$-word and $n=\mu_{E}(\pm t)$. This shows one of the conditions (a), (b) or (c) must hold.

(ii) Let $C$ and $E$ be periodic. Assume momentarilly, for some $t\in\mathbb{Z}$, that: $E=C[t]$ and $n=\mu_{C}(t)$; or that $E=C^{-1}[t]$ and $n=\mu_{C}(-t)$. By definition this means $(H,L,0)\sim (C(1),C(-1),-n)$ where $H=(E_{\leq 0})^{-1}$ and $L=E_{>0}$. By Corollary \ref{allbandsgood} we have that: if ($E=C[t]$, $k\otimes_{R}V\simeq k\otimes_{R}W$ and $n=\mu_{C}(t)$) then $S_{C(1),C(-1),-n}(V)\simeq S_{H,L,0}(W)$; and if ($E=C^{-1}[t]$, $k\otimes_{R}V\simeq k\otimes_{R}\mathrm{res}_{\iota,R} \,W$ and $n=\mu_{C}(-t)$) then $S_{C(1),C(-1),-n}(V)\simeq S_{H,L,0}(\mathrm{res}_{\iota,R}(\mathrm{res}_{\iota,R}(W)))$. Since $\mathrm{res}^{2}_{\iota,R}=\mathrm{id}$ we have $P(C,V)[n]\simeq P(E,W)$ in $\mathcal{K}(\Lambda\text{-}\boldsymbol{\mathrm{Proj}})$.

For the converse suppose $P(C,V)[n]\simeq P(E,W)$ in $\mathcal{K}(\Lambda\text{-}\boldsymbol{\mathrm{Proj}})$. As above we have $N'=P(C,V)[n]$ and $N=P(E,W)$ where $\mathcal{S}=\mathcal{S}'=\emptyset$, $\mathcal{B}=\{\beta\}$, $\mathcal{B}'=\{\beta'\}$, $E(\beta)=E$, $D(\beta')=C$,  $s(\beta)=0$ and $r(\beta')=-n$. By Theorem \ref{theorem.7.1}(ii) we have that, for any $(H,L,q)\in\Sigma$, $F_{H,L,q}(N)\simeq \bigoplus W$ (respectively $F_{H,L,q}(N')\simeq \bigoplus V$) in $k[T,T^{-1}]\text{-}\boldsymbol{\mathrm{Mod}}$) where the sum runs through all $\beta\in \mathcal{B}$ (respectively $\beta'\in\mathcal{B}'$) with $-q=r(H,L;(E_{\leq 0})^{-1},E_{>0})$ (respectively  $n-q=r(H,L;C(1),C(-1))$) and $E\sim H^{-1}L$ (respectively $C\sim H^{-1}L$). 

Applying $F_{H,L,n}$ to the isomorphism $N\simeq N'$ in case $(H,L)=(C(1),C(-1)))$ gives $F_{H,L,n}(N)\simeq k\otimes_{R} V$, which shows that $n=r((E_{\leq 0})^{-1},E_{>0};C(1),C(-1))$ and $C\sim E$. By definition this means that: $E=C[t]$ and $n=\mu_{C}(t)$, in which case $F_{H,L,n}(N)\simeq k\otimes_{R} W$; or that $E=C^{-1}[t]$ and $n=\mu_{C}(-t)$, in which case $F_{H,L,n}(N)\simeq k\otimes_{R} \mathrm{res}_{\iota,R}(W)$.

(iii) Using similar ideas to the above, one can show that if $P(C)[n]\simeq P(E,V)$ then we must have $C\sim E$, which is impossible if $C$ is not periodic, yet $E$ is periodic.
\end{proof}
\begin{proof}[of Theorem \ref{theorem.1.3}] Let $(B,D,n)\in\mathcal{I}$. Suppose $N\simeq N'$ in the notation of Assumption \ref{finalass}. We define the bijection $\chi:\mathcal{S}\cup\mathcal{B}\rightarrow\mathcal{S}'\cup\mathcal{B}'$ as follows. 

Suppose firstly that $(B,D,n)\in\mathcal{I}(s)$. As in the proof of Theorem \ref{theorem.1.2}, by Theorem \ref{theorem.7.1}(i) there is a bijection $\varphi:\mathcal{S}(B,D,n)\rightarrow \mathcal{S}'(B,D,n)$. By Lemma \ref{lemma.4.1}, if $\varphi(\sigma)=\sigma'$ then both $P(A(\sigma))[-t(\sigma)]$ and $P(B(\sigma'))[-u(\sigma')]$ are isomorphic to $P(C)[-n]$. If instead $(B,D,n)\in\mathcal{I}(b)$ then, as $F_{B,D,n}(N)\simeq F_{B,D,n}(N')$, by Theorem \ref{theorem.7.1}(ii) we have
\[
\begin{array}{c}
\bigoplus_{\beta_{+}}\bar{V}_{+}(\beta_{+})\oplus\bigoplus_{\beta_{-}}\bar{V}_{-}(\beta_{-})\simeq \bigoplus_{\beta'_{+}}\bar{W}_{+}(\beta'_{+})\oplus\bigoplus_{\beta'_{-}}\bar{W}_{-}(\beta'_{-})
\end{array}
\]
where $\beta_{\pm}$ (respectively $\beta'_{\pm}$) runs through $\mathcal{B}(B,D,n)^{\pm}$ (respectively $\mathcal{B}'(B,D,n)^{\pm}$). As in the proof of Theorem \ref{theorem.1.2}, by Theorem \ref{theorem.7.1}(ii) and by the Krull-Remak-Schmidt
property for $k[T,T^{-1}]\text{-}\mathbf{Mod}_{k\text{-}\mathbf{mod}}$ there is an isomorphism class preserving bijection
\[
\psi:\mathcal{B}(B,D,n)^{-}\cup\mathcal{B}(B,D,n)^{+}\to \mathcal{B}'(B,D,n)^{-}\cup\mathcal{B}'(B,D,n)^{+}.
\]
Provided  $\psi(\beta)=\beta'$ for some $\beta\in\mathcal{B}(B,D,n)^{\pm}$, note that: if $\beta'\in\mathcal{B}'(B,D,n)^{\pm}$ then $\bar{V}_{\pm}(\beta)\simeq \bar{W}_{\pm}(\beta')$; and if $\beta'\in\mathcal{B}'(B,D,n)^{\mp}$ then $\bar{V}_{\pm}(\beta)\simeq \bar{W}_{\mp}(\beta')$. By Corollary \ref{allbandsgood}, if $\psi(\beta)=\beta'$ then 
\[
P(E(\beta),V^{\beta})[-s(\beta)] \simeq S_{B,D,n}(V_{\pm}(\beta))
\simeq P(D(\beta'),W^{\beta'})[-r(\beta')]
\] 
Define $\chi:\mathcal{S}\cup\mathcal{B}\rightarrow\mathcal{S}'\cup\mathcal{B}'$ by $\chi(\alpha)=\alpha'$ if and only if (($\alpha\in\mathcal{S}(B,D,n)$ and $\alpha'=\varphi=(\alpha)$) or ($\alpha\in\mathcal{B}(B,D,n)^{+}\cup\mathcal{B}(B,D,n)^{-}$ and $\alpha'=\psi(\alpha)$)) for some $(B,D,n)\in\mathcal{I}$. Note that $\chi$ is well defined and bijective by Lemma \ref{lemma.7.1}.\end{proof}

\bibliography{biblio}
\bibliographystyle{abbrv}

%
\end{document}